%% file: main.tex
\documentclass[USenglish]{article}

\usepackage[USenglish]{babel}
\usepackage{color}
\usepackage{amsmath,amsthm,amssymb,amsfonts}
\usepackage[margin=2.5 cm]{geometry}
\usepackage[unicode,colorlinks=true,linkcolor=blue,citecolor=green,urlcolor=black,pagebackref=false]{hyperref}
\usepackage{stmaryrd}
\usepackage{mathabx}
\usepackage{custom_macros_general, custom_macros_spaceops, custom_macros_funcs, custom_macros_freqops,mathtools}

\let\oldmu\mu
\let\oldeps\varepsilon
\let\oldnu\nu
\let\oldlambda\lambda
\let\oldzeta\zeta

\newcommand{\revisioncolorred}{}

\renewcommand{\mu}{\boldsymbol{\oldmu}}
\renewcommand{\varepsilon}{\boldsymbol{\oldeps}}
\renewcommand{\nu}{\boldsymbol{\oldnu}}
\renewcommand{\lambda}{\boldsymbol{\oldlambda}}
\renewcommand{\zeta}{\boldsymbol{\oldzeta}}

\def\revision#1{{#1}}

\def\revisionc#1{{#1}}
\numberwithin{equation}{section}

\DeclareMathAlphabet\mathbfcal{OMS}{cmsy}{b}{n}
\makeatletter
\let\@fnsymbol\@arabic
\makeatother

\title{Wavenumber-explicit analytic regularity of the heterogeneous Maxwell equations with impedance boundary conditions}

\author{Jens Markus Melenk\thanks{(melenk@tuwien.ac.at), Institut f\"{u}r Analysis und Scientific Computing, Technische Universit\"{a}t Wien, Wiedner Hauptstrasse 8-10, A--1040 Wien, Austria.}
 \and 
David W\"{o}rg\"{o}tter\thanks{(david.woergoetter@tuwien.ac.at), Institut f\"{u}r Analysis und Scientific Computing, Technische Universit\"{a}t Wien, Wiedner Hauptstrasse 8-10, A--1040 Wien, Austria.}}

\date{\today}

\begin{document}

\maketitle
\begin{abstract}
		We consider the time-harmonic Maxwell equations at a wavenumber $k\in\Co\setminus\{0\}$ on a bounded and simply connected Lipschitz domain $\Omega$ with an analytic boundary $\Gamma$, on which we impose impedance boundary conditions. We suppose that the (possibly complex-valued) permeability and permittivity tensor fields $\mu^{-1}$ and $\varepsilon$ are piecewise analytic in $\Omega$ and discontinuous only across certain mutually disjoint analytic surfaces inside of $\Omega$. We show that under these circumstances, any weak solution of Maxwell's equations is piecewise analytic in $\Omega$ and that the growth of its derivatives can be controlled explicitly in the wavenumber $k$.
\end{abstract}

\input{intro.tex}

\input{notation/notation.tex}

	\input{analytic/analytic.tex}

\section*{Acknowledgement}
    JMM acknowledges funding by the Austrian Science Fund (FWF) under
    grant F65 ``taming complexity in partial differential systems'' (\href{https://doi.org/10.55776/F65}{DOI:10.55776/F65}).

	\appendix
	\input{appendix/appendix1.tex}

	\input{appendix/appendix2.tex}

	\bibliographystyle{amsplain}
	\bibliography{bibliog}
\end{document}

%% file: intro.tex
\section{Introduction}

We consider the heterogeneous time-harmonic Maxwell equations at wavenumber $k\in\Co\setminus\{0\}$ posed on a bounded and simply connected Lipschitz domain $\Omega$ with analytic boundary $\Gamma:=\partial\Omega$. 

These equations read as follows: For a given right-hand side $\solf$ and a tangent field $\solgi$ on $\Gamma$, find a vector field $\solu:\Omega\rightarrow\Co^3$ satisfying

\begin{equation}\label{Maxwellorig}
\begin{alignedat}{2}
				\curl\mu^{-1}\curl\solu-k^2\varepsilon\solu &= \solf &&\quad {\rm in}\ \Omega, \\
				\left(\mu^{-1}\curl\solu\right)\times\soln-ik\zeta\solu_T &= \solgi &&\quad {\rm on}\ \Gamma,
\end{alignedat}
\end{equation}
where $\mu^{-1}$ and $\varepsilon$ are complex-valued tensor fields which satisfy \revision{certain coercivity conditions (see Assumption~\ref{assumptionsepsmu} below)} and are piecewise analytic in $\Omega$, but may be discontinuous across certain mutually disjoint analytic surfaces $\interf_1,\ldots,\interf_r$ inside of $\Omega$. 

Moreover, $i:=\sqrt{-1}$ is the imaginary unit, $\soln$ is the outer unit normal to $\Gamma$, $\solu_T := \soln\times(\solu\times\soln)$ is the tangential component of a vector field $\solu$ and $\zeta:\Gamma\rightarrow\Co^{3\times 3}$ is a given tensor field which is coercive and analytic on $\Gamma$ and satisfies certain admissibility conditions described in Assumption~\ref{assumptionsepsmu} below. 

\revisionc{Maxwell's equations are fundamental in the theory of electromagnetism and they have numerous applications in physics and engineering. As such they have been extensively studied in the literature both from an analytical (see e.g. \cite{BookCessenat, BookColton,BookNedelec}) and from a numerical point of view, see for example \cite{Hiptmair1,Hiptmair2,BookMonk}. While the Sobolev regularity properties of solutions of Maxwell's equations \eqref{Maxwellorig} are by now well-understood both in smooth domains \cite{MaxwellImpedanceMelenk,MaxwellMyself,RegularityWeber} and in polyhedral domains \cite{CostabelDauge}, the case of wavenumber-explicit analytic regularity in the case of piecewise analytic input data appears to be open.} 

\revisionc{	Rather recently, however, wavenumber-explicit analytic regularity estimates for solutions $\solu$ of \eqref{Maxwellorig} turned out to be a powerful tool in the analysis of numerical schemes for Maxwell's equations \cite{MaxwellImpedanceMelenk, MaxwellTomezyk}, hence proving such estimates is of great practical interest. 
		So far, wavenumber-explicit analytic regularity estimates for Maxwell's equations have been derived for the case of homogeneous coefficients $\mu^{-1}=\varepsilon=\zeta=1$, see \cite{MaxwellImpedanceMelenk, MaxwellTomezyk}, but for heterogeneous coefficients such estimates have been lacking so far. Applications include physical situations where several materials are present in the computational domain. 
}

\bigskip

The aim of this work is \revisionc{to consider the case of heterogeneous and piecewise analytic coefficients and} to prove a) local analyticity of $\solu$ in regions where $\mu^{-1}$ and $\varepsilon$ are analytic, and b) piecewise analyticity of $\solu$ on the whole domain $\Omega$, provided that $\Gamma$ and the surfaces of discontinuity $\interf_1,\ldots,\interf_r$ are analytic. On top of that, we aim for $k$-explicit estimates for the growth of the derivatives of $\solu$. \revisionc{We highlight that in this work we restrict ourselves to the case of an analytic boundary $\Gamma$ and analytic surfaces of discontinuity $\interf_i$. Considering more complicated domains (like domains with only piecewise analytic boundary or polyhedral domains) is beyond the scope of the present paper.}

\revisionc{For constant scalar coefficients $\mu^{-1},\varepsilon$ the approach of }\cite{MaxwellTomezyk, MaxwellImpedanceMelenk} is to exploit the identity 
\begin{align}\label{curlcurldelta}
\curl\curl\solu = \nabla(\diverg\solu)-\Delta\solu
\end{align}
to transform \eqref{Maxwellorig} to a set of three Helmholtz equations. Simply speaking, this allows \revision{one} to reduce the problem of analyticity of Maxwell's equations to that of the analyticity of the vector Laplacian; this reduces the complexity of the problem since the vector Laplacian allows for a direct application of the well-understood regularity theory of elliptic operators \cite{ADN1, ADN2, MorreyNirenberg, BookMorrey, BabuskaGuo1, Babuskaguo2}.

In the setting of this work, however, the identity \eqref{curlcurldelta} cannot be used any more. Instead, we rely on shift results for vector fields with piecewise regular curl and divergence from \cite{MaxwellMyself} and nested ball techniques also used in e.g. \cite{GrandLivre, MaxwellTomezyk} to control the growth of derivatives of both $\solu$ and $\mu^{-1}\curl\solu$ in the interior of $\Omega$. This will provide analyticity of $\solu$ on regions inside of $\Omega$ where $\mu^{-1}$ and $\varepsilon$ are analytic. \revisionc{Let us mention that the techniques that we use in this work are inspired by and very similar to the techniques used in \cite{GrandLivre} to prove analyticity of solutions of elliptic problems. One major novel aspect, however, is that in this work we include a wavenumber-dependence in the nested ball techniques; this leads to wavenumber-explicit analytic estimates for solutions of Maxwell's equations.}

\medskip

In order to infer analyticity near the boundary $\Gamma$ or near subdomain interfaces $\interf_1,\ldots,\interf_r$ we encounter the problem of $\revisionc{\solu\mapsto}\curl\mu^{-1}\curl\solu$ not being an elliptic operator. Hence, although the aforementioned methods allow us to control tangential derivatives of $\solu$ and $\mu^{-1}\curl\solu$, we cannot exploit ellipticity of the differential operator at hand to infer estimates on normal derivatives. The remedy to this problem is the observation that \eqref{Maxwellorig} implies $-k^2\diverg\varepsilon\solu = \diverg\solf$; from this and the original equation \eqref{Maxwellorig} we are able to control the normal derivatives of $\solu$ on $\Gamma$ or on $\interf_1,\ldots,\interf_r$ in terms of tangential derivatives of $\solu$ and $\solf$; for further details we refer to Lemma~\ref{linklemmalem} below. \revision{Let us highlight that the techniques that we use in this work can also be applied to discuss other boundary conditions. For example, it is possible to extend the results of this work to Maxwell's equations \eqref{Maxwellorig} subject to natural boundary conditions $\mu^{-1}\curl\solu\times\soln = \solgn$ on $\Gamma$ for some given $\solgn$ or essential boundary conditions $\solu_T=0$ on $\Gamma$.} \revisionc{For the sake of simplicity of exposition we refrain from discussing other boundary conditions in this work, however, let us mention that an extension of the results from this paper to Maxwell's equations with natural or essential boundary conditions will be provided in \cite{MyDiss}.}

\medskip

We mention that this work is the second paper in a series of research articles. The first paper \cite{MaxwellMyself} provides piecewise regularity shift results that are fundamental for the analytic regularity properties of $\solu$ derived here in this work. The present paper is the key stepping stone for the design and analysis of wavenumber-explicit high order finite element methods for Maxwell's equations in heterogeneous media. Indeed, in \cite{MaxwellMyself3} we will provide a wavenumber-explicit analysis of the $hp$-finite element method applied to Maxwell's equations \eqref{Maxwellorig} in the presence of piecewise analytic coefficients and impedance boundary conditions. In particular, if $k$ denotes the wavenumber, $h$ is the mesh width and $p$ denotes the polynomial degree in the $hp$-FEM approximation, we will show in \cite{MaxwellMyself3} that the $hp$-FEM is quasi-optimal under the scale-resolution conditions a) that $h|k|/p$ is sufficiently small and b) that $p/\log|k|$ is sufficiently large. Additionally, the $k$-explicit regularity theory of the present paper allows one to quantify best approximation results for the solution of Maxwell's equations that are explicit in the mesh size $h$ and the approximation order $p$. Maxwell's equations with piecewise constant coefficients can also be addressed by boundary integral equation (BIE) methods. By taking traces, the present PDE regularity could be used to infer regularity results for the solution of such BIEs and develop a $k$-explicit analysis of their discretizations.

%The PDE results obtained here may also pave the way to obtain regularity properties of solutions of integral equations and a subsequent analysis of boundary element methods for Maxwell's equations with piecewise constant coefficients.

\medskip 

The outline of this work is as follows: In Section~\ref{notation} we introduce the general notation used throughout this work. Among others, this includes the concepts of piecewise Sobolev spaces, piecewise analytic vector- and tensor fields and $\ana$-partitions. We conclude Section~\ref{notation} by stating the main results of this work, namely Theorem~\ref{Mainresult1} and Theorem~\ref{Mainresult2}.

\medskip 

Subsequently, in Section~\ref{diffops} we recall the surface curl and the surface divergence operators from, e.g., \cite{BookMonk,BookNedelec}. These operators are needed for the variational formulation of Maxwell's equations \eqref{Maxwellorig}, which is stated at the end of Section~\ref{diffops}.

In Section~\ref{analyticinteriorsec} we define and develop techniques and methods to prove Theorem~\ref{Mainresult1}. We define Sobolev-Morrey seminorms and commutator fields, present commutator estimates that are inspired by \cite[Lem.~1.6.2]{GrandLivre} and recall useful results from \cite{MaxwellMyself}. We then conclude Section~\ref{analyticinteriorsec} by giving a proof of Theorem~\ref{Mainresult1}.

\medskip

Following that, we turn our attention to the process of locally flattening the boundary $\Gamma$ or the surfaces of discontinuity $\interf_1,\ldots,\interf_r$ to the plane $\{x_3=0\}$. We describe this flattening process in great detail in Section~\ref{flatteningsec}, as it is of utmost importance for the proof of Theorem~\ref{Mainresult2}. In fact, we restrict ourselves to a certain type of transformations which we call {\it local flattenings} and after defining these special transformations we show that Maxwell's equations \eqref{Maxwellorig} are invariant under transformations by local flattenings.

\medskip 

As we shall see in Section~\ref{analyticinteriorsec}, the piecewise regularity results from \cite{MaxwellMyself} are fundamental for the proof of Theorem~\ref{Mainresult1}. The proof of Theorem~\ref{Mainresult2} follows the same principles as the proof of Theorem~\ref{Mainresult1}, however, as it involves locally flattening $\Gamma$ and $\interf_1,\ldots,\interf_r$, the technical details are much more involved. In particular, locally flattening $\Gamma$ and $\interf_1,\ldots,\interf_r$ changes the underlying geometry; instead of $\Omega$ and subdomains of $\Omega$ we have to consider balls and half-balls centered at the origin. Hence, we have to adapt the piecewise regularity results from \cite{MaxwellMyself} to this new geometric setting. In addition, we need to study the behaviour of vector- and tensor fields under local flattenings. Performing these tasks is the aim of Section~\ref{adaptnotation}.

\medskip

Afterwards, in Section~\ref{auxiliarsec} we study the interplay between normal and tangential derivatives of a solution $\solu$ of Maxwell's equations \eqref{Maxwellorig} near $\Gamma$ or near a surface of discontinuity $\interf_1,\ldots,\interf_r$. Indeed, Lemma~\ref{linklemmalem} provides a link between normal and tangential derivatives; more precisely, it allows us to control the normal derivatives of $\solu$ by its tangential derivatives, $\solu$, $\curl\solu$ and $\solf$. Subsequently we provide estimates for tangential derivatives of $\solu$, and in combination with Lemma~\ref{linklemmalem}, this provides combined estimates for the tangential and the normal derivatives of $\solu$ near $\Gamma$ and near $\interf_1,\ldots,\interf_r$. 

\medskip

Finally, in Section~\ref{finalsec} we exploit the results from Section~\ref{auxiliarsec} to give a proof of Theorem~\ref{Mainresult2}, which concludes this work except for the proofs of some rather technical results; these proofs are postponed to the subsequent Appendix~\ref{appendix:commest}.

%% file: notation/notation.tex
\section{General notation and main results}\label{notation}

For any two vectors $\solw,\solz\in\Co^3$ with $\solw = (w_1, w_2, w_3)^T$ and $\solz=(z_1,z_2,z_3)^T$ we set $\solw\cdot\solz:=\sum_{i=1}^3w_iz_i$ and write $\SCP{\solw}{\solz}{}:=\solw\cdot\overline{\solz}$ for the scalar product between $\solw$ and $\solz$, where $\overline{\solz}:=(\overline{z_1},\overline{z_2},\overline{z_3})^T$ denotes the complex conjugate of $\solz$. 
Furthermore, the cross product between the vectors $\solw$ and $\solz$ is defined in the usual way as $\solw\times\solz:=(w_2z_3-w_3z_2, w_3z_1-w_1z_3, w_1z_2-w_2z_1)^T$.
  
\medskip 

  As usual, let $\SL(\Omega)$ denote the Lebesgue space of complex-valued square integrable functions, and define its vector-valued version as $\VSL(\Omega):=(\SL(\Omega))^3$.
  For a bounded Lipschitz domain $\Omega\subseteq \R^3$ and \revision{$\ell\in\N_0$, the space $\SHM{\ell}(\Omega)$ is the usual Sobolev space of order $\ell$, see \cite[Ch.~3]{BookMcLean}, and $\SHM{\ell}_0(\Omega)$ denotes the closure of $\CM{\infty}_0(\Omega)$ in $\SHM{\ell}(\Omega)$. In order to deal with vector fields, we define the vector-valued space $\VSHM{\ell}(\Omega):=(\SHM{\ell}(\Omega))^3$.}
  In addition, we define the spaces
  \begin{align*}
  \Hcurl := \{\solu\in\VSL(\Omega)\ |\ \curl \solu\in\VSL(\Omega)\}\quad {\rm and}\quad \Hdiv:=\{\solu\in\VSL(\Omega)\ |\ \diverg \solu\in\SL(\Omega)\}.
  \end{align*}
  %which are equipped with the norms
  %\begin{align*}
  %\norm{\solu}{\Hcurl}^2:=\norm{\solu}{\VSL(\Omega)}^2+\norm{\curl\solu}{\VSL(\Omega)}^2\ {\rm and}\ \norm{\solu}{\Hdiv}^2:=\norm{\solu}{\VSL(\Omega)}^2+\norm{\diverg\solu}{\SL(\Omega)}^2,
  %\end{align*}
  %respectively.

\medskip

Furthermore, if $\ell\in\N_0$ and if $\Sigma\subseteq\R^3$ is the (smooth) boundary of a bounded and smooth Lipschitz domain $\Omega_{\Sigma}\subseteq\R^3$, let $\SHM{\ell}(\Sigma)$ be the Sobolev space of order $\ell$ with dual space $\SHM{-\ell}(\Sigma)$, see \cite[Ch.~3]{BookMcLean}. Moreover, for $\ell\in\N_0$ we define the Sobolev spaces $\SHM{\ell+1/2}(\Sigma)$ as trace spaces
		\begin{align*}
				\SHM{\ell+1/2}(\Sigma):=\{U\vert_{\Sigma}\ \big|\ U\in\SHM{\ell+1}(\Omega_{\Sigma})\}
		\end{align*}
		with norm
		\begin{align*}
				\norm{u}{\SHM{\ell+1/2}(\Sigma)}:=\inf\left\{\norm{U}{\SHM{\ell+1}(\Omega_{\Sigma})}\ \big| \ U\vert_{\Sigma}=u\right\}.
		\end{align*}
		For $s\in\R$ with $2s\in\N_0$ we then set $\VSHM{\ell}(\Sigma):=(\SHM{\ell}(\Sigma))^3$. Furthermore, with $\soln$ denoting the outer unit normal to $\Omega_{\Sigma}$ we define the space of square-integrable tangent fields by

\begin{align*}
\VSL_T(\Sigma):=\{\solv\in\VSL(\Sigma)\ |\ \solv\cdot\soln = 0\}.
\end{align*}
\revisionc{For} \revision{$s\in\R$ satisfying $2s\in\N_0$ we set} 
\begin{align*}
\VSHM{s}_T(\Sigma):=\VSL_T(\Sigma)\cap\VSHM{s}(\Sigma),\ \ {\rm as\ well\ as}\ \ \VSHM{-s}_T(\Sigma):=(\VSHM{s}_T(\Sigma))'.
\end{align*}

\revision{From Section~\ref{adaptnotation} onwards we will also work with Sobolev spaces on discs of the form
		\begin{align}\label{PZMdef}
		\PZ_R := \left\{x,y\in\R\ \big|\ x^2+y^2 < R^2\right\}, \quad R>0.
		\end{align}
		Note that $\PZ_R$ is not the boundary of a Lipschitz domain, hence the previous definitions do not apply. We postpone the precise definitions of (fractional) Sobolev spaces on $\PZ_R$ to the beginning of Section~\ref{adaptnotation}.
}
%%\revision{
%%		From Section~\ref{adaptnotation} onwards we will frequently work with fractional Sobolev spaces on discs of the form
%%		$$\PZ_R := \left\{x,y\in\R\ \big|\ x^2+y^2\leq R^2\right\}, \quad R>0.$$
%%		We define $\SHM{1/2}(\R^2)$ in terms of the Fourier transform cf. \cite[Chapter 3]{BookMcLean} and set
%%		\begin{align*}
%%				\SHM{1/2}(\PZ_R):=\left\{u\in\SL(\PZ_R)\ \big| \ u = U\vert_{\PZ_R}\ {\rm for\ some}\ U\in\SHM{1/2}\left(\R^2\right)\right\}.
%%		\end{align*}
%%		In addition, we define
%%		\begin{align*}
%%				\SHMZ{1/2}(\PZ_R):=\left\{u\in\SHM{1/2}(\PZ_R)\ \big| \ \operatorname{supp} u\in\PZ_R\right\}
%%		\end{align*}
%%		and we write $\SHM{-1/2}(\PZ_R)$ for the dual of $\SHMZ{1/2}(\PZ_R)$. Finally, we set 
%%		\begin{align*}
%%				\VSHM{1/2}_T(\PZ_R):=\left\{\solv\in\left(\SHM{1/2}(\PZ_R)\right)^3\ \big|\ \solv\cdot\sole_3 = 0\right\},\ {\rm where}\ \sole_3 = (0,0,1)^T.
%%		\end{align*}
%%}
%%

\medskip

For a bounded Lipschitz domain $\Omega\subseteq\R^3$ we say that a function $u$ or a vector field $\solu$ is smooth in $\Omega$ if $u\in\CM{\infty}(\Omega)$ or $\solu\in\left(\CM{\infty}(\Omega)\right)^3$, respectively. In addition, a function $u$ or a vector field $\solu$ is said to be smooth up to the boundary of $\Omega$ if $u$ or $\solu$ can be extended to a function $\widetilde{u}\in\CM{\infty}(\R^3)$ or a vector field $\widetilde{\solu}\in\left(\CM{\infty}(\R^3)\right)^3$, respectively. In this case we write $u\in\CM{\infty}(\overline{\Omega})$ or $\solu\in\VCM{\infty}(\overline{\Omega})$, respectively.
In general, $\overline{\Omega}$ denotes the closure of $\Omega$, and for two bounded Lipschitz domains $\Omega$ and $\Omega_2$ we write $\Omega_2\subset\subset\Omega$ if $\overline{\Omega}_2\subseteq\Omega$.

Concerning tensor fields, a tensor field $\nu:\Omega\rightarrow\Co^{3\times 3}$ is called {\bf coercive}, if there exists a constant $c>0$ and a complex number $\alpha$ with $|\alpha|=1$ such that
\begin{align}\label{coercivenudomain}
		\revision{\forall\solz\in\Co^3\ \forall x\in\Omega:}\	\realpart\SCP{\alpha\nu\solz}{\solz}{}\geq c\norm{\solz}{}^2.
\end{align}

Analogously, if $\nu:\Gamma\rightarrow\Co^{3\times 3}$ is a tensor field on the boundary $\Gamma$ of a bounded Lipschitz domain, $\nu$ is called {\bf coercive} if there exist $c>0$ and $\alpha\in\Co$ with $|\alpha|=1$ such that
\begin{align}\label{coercivenubdr}
		\revision{\forall\solz\in\Co^3\ \forall x\in\Gamma :}\	\realpart\SCP{\alpha\nu\solz}{\solz}{}\geq c\norm{\solz}{}^2.
\end{align}

Finally, we need to define classes of analytic functions on $\Omega$. Since we aim for wavenumber-explicit analytic regularity estimates for solutions of \eqref{Maxwellorig}, it turns out to be beneficial to incorporate the wavenumber $k$ in the subsequent definitions of analyticity classes.
 Let us mention that for a multiindex $\alpha\in\N_0^3$ and a vector field $\solu = (u_1, u_2, u_3)^T$ we set $\D^{\alpha}(\solu):=\left(\D^{\alpha}(u_1),\D^{\alpha}(u_2), \D^{\alpha}(u_3)\right)^T$. 

\begin{definition}\label{analytictensors}
		For a bounded Lipschitz domain $\Omega\subseteq\R^3$, a wavenumber $k\in\Co\setminus\{0\}$ and two constants $\omega\geq 0, M>0$, the symbols $\anaf{\Omega}$ and $\anavec{\Omega}$ denote the sets of all functions $u:\Omega\rightarrow\Co$ and vector fields $\solu:\Omega\rightarrow\Co^3$ which satisfy
\begin{align*}
		\revision{\forall\ell\in\N_0:}\ \sum_{|\alpha|=\ell}\norm{\D^{\alpha}u}{\SL(\Omega)}\leq \omega M^{\ell}(|k|+\ell)^{\ell}\quad {\rm and}\quad 
		\revision{\forall\ell\in\N_0:}\ \sum_{|\alpha|=\ell}\norm{\D^{\alpha}(\solu)}{\VSL(\Omega)}\leq \omega M^{\ell}(|k|+\ell)^{\ell}.
\end{align*}

Finally, the symbol $\anatens{\Omega}$ denotes the space of all tensor fields $\nu:\Omega\rightarrow\Co^{3\times 3}$ for which there exist constants $\omega\geq 0,M>0$ such that 
\begin{align*}
		\revision{\forall\ell\in\N_0:}\ \sum_{i,j=1}^3\sum_{|\alpha|=\ell}\norm{\D^{\alpha}(\nu_{i,j})}{\SL(\Omega)}\leq \omega M^{\ell}\ell^{\ell}.
\end{align*}

%We say that a tensor field $\nu\in\anatens{\Omega}$ is coercive if there exist constants $\alpha\in\Co$ with $|\alpha|=1$ and $c>0$ such that 

%\begin{align*}
%		\realpart\SCP{\alpha\nu\solz}{\solz}{}\geq c\norm{\solz}{}^2
%\end{align*}
%for all $\solz\in\Co^3$ uniformly in $\Omega$.
\end{definition}
\begin{remark}\label{kmimpl}
		\revision{Let us highlight that we specifically mentioned $\omega, M$ and $k$ in the notions of $\anaf{\Omega}$ and $\anavec{\Omega}$ but not in $\anatens{\Omega}$. This is due to the fact that we are primarily interested in the influence of the given data $\solf$ and $\solgi$ on the solution $\solu$, whereas the influence of the tensor fields $\mu^{-1}$, $\varepsilon$ and $\zeta$ is left implicit.} \revisionc{This is done for simplicity of exposition: It allows us to track the dependence on the analytic data $\solf, \solgi$ but does not precisely track the dependence on $\mu^{-1}, \varepsilon$ and $\zeta$.} 
\end{remark}

In addition to analytic functions, vector fields and tensor fields in the domain $\Omega$ we will also need a notion for the corresponding entities on the boundary $\Gamma$ of a bounded Lipschitz domain $\Omega$ with analytic boundary. To that end, we introduce the following definition:

\begin{definition}\label{assumptionimpedance}
		Let $\Omega\subseteq\R^3$ be a bounded Lipschitz domain with analytic boundary $\Gamma$ and suppose $k\in\Co\setminus\{0\}$ as well as $\omega\geq 0, M>0$. The sets $\anagammaf$ and $\anagammavec{\Gamma}$ consist of all functions $g:\Gamma\rightarrow\Co$ and tangent fields $\solg_T:\Gamma\rightarrow\Co^3$ for which there exist continuations $g^*$ and $\solg^*$ to \revision{an open neighbourhood $\mathcal{U}_{\Gamma}$} of $\Gamma$ which satisfy
		\begin{align*}
				\revision{\forall\ell\in\N_0:}\	\sum_{|\alpha|=\ell}\norm{\D^{\alpha}g^*}{\SL(\mathcal{U}_{\Gamma})}\leq \omega M^{\ell}(|k|+\ell)^{\ell}\quad{\rm and}\quad
				\revision{\forall\ell\in\N_0:}\ \sum_{|\alpha|=\ell}\norm{\D^{\alpha}(\solg^*)}{\VSL(\mathcal{U}_{\Gamma})}\leq \omega M^{\ell}(|k|+\ell)^{\ell}.
		\end{align*}

		Finally, the space $\anagamma$ consists of all tensor fields $\zeta:\Gamma\rightarrow\Co^{3\times 3}$ for which there exist constants $\omega\geq 0, M>0$ and a continuation $\zeta^*$ on an open neighbourhood $\mathcal{U}_{\Gamma}$ which satisfies
\begin{align*}
		\revision{\forall\ell\in\N_0:}\ \sum_{i,j=1}^3\sum_{|\alpha|=\ell}\norm{\D^{\alpha}(\zeta^*_{i,j})}{\SL(\mathcal{U}_{\Gamma})}\leq \omega M^{\ell}\ell^{\ell}.
\end{align*}

%A tensor field $\zeta\in\anagamma$ is called coercive, if there exist $\alpha\in\Co$ with $|\alpha|=1$ and $c>0$ such that 
%\begin{align*}
%	\realpart\SCP{\alpha\zeta\solz}{\solz}{}\geq c\norm{\solz}{}^2
%\end{align*}
%for all $\solz\in\Co^3$ uniformly on $\Gamma$.
\end{definition}

\subsection{$\ana$-partitions and classes of piecewise analytic functions}

Let $\Omega\subseteq\R^3$ be a bounded Lipschitz domain with boundary $\Gamma$. 
Throughout most of this work we suppose that $\Omega$ is partitioned into subdomains $\Gp_1,\ldots,\Gp_n$ such that the boundaries of these subdomains form a family of mutually disjoint closed\footnote{By {\it closed surface} we denote a compact surface without boundary.} and analytic surfaces inside of $\Omega$. 
The following definition makes this precise.

\begin{definition}\label{partitiondef}
		An $\ana$-partition $\Gp$ is a tuple $\Gp = \geom$ which consists of domains $\Omega, \Gp_1,\ldots, \Gp_n\subseteq\R^3$ satisfying 
		\begin{itemize}
				\item[(i)] The domains $\Omega, \Gp_1,\ldots,\Gp_n$ are bounded Lipschitz domains in $\R^3$ with $\Omega$ being simply connected and with the domains $\Gp_1,\ldots\Gp_n$ being mutually disjoint and satisfying $\overline{\Omega} = \overline{\Gp_1}\cup\ldots\cup\overline{\Gp_n}$. 
				\item[(ii)] The boundary $\Gamma:=\partial\Omega$ is analytic and simply connected.
				%\item[(ii)] The boundaries $\Gamma:=\partial\Omega$ and $\partial\Gp_1,\ldots,\partial\Gp_n$ are analytic and each of these boundaries is either simply connected or consists of simply connected components. 
				\item[(iii)] If $n>1$ there exist bounded Lipschitz domains $V_1,\ldots,V_r$ each having a simply connected and analytic boundary $\interf_i:=\partial V_i$, such that the $\interf_i$ are mutually disjoint and disjoint to $\Gamma$ and
						\begin{enumerate}
						\item For every $i\in\{1,\ldots,m\}$ there exist distinct indices $j_1,j_2\in\{1,\ldots,n\}$ such that $\interf_i=\partial\Gp_{j_1}\cap\partial\Gp_{j_2}$.
						\item There holds $$
		\Gamma\cup\bigcup_{j=1}^r\interf_j = \bigcup_{j=1}^n\partial\Gp_j.$$
						\end{enumerate}
						
				%\item[(iii)] If $n>1$, there exist closed and simply connected analytic surfaces $\interf_1,\ldots, \interf_r$ such that $\Gamma, \interf_1,\ldots,\interf_r$ are mutually disjoint and such that 
\end{itemize}

If $n>1$, the surfaces $\interf_1,\ldots,\interf_r$ are called subdomain interfaces and their union $\interf := \interf_1\cup\ldots\cup\interf_r$ is referred to as the union of all subdomain interfaces.
\end{definition}

The point of introducing subdomains in the preceding fashion is to incorporate the location of (possible) discontinuities of piecewise regular coefficients $\mu^{-1}$ and $\varepsilon$ into the geometry of the problem. If $n=1$, then $\Omega$ has only one subdomain, namely itself, hence in this case there is no subdomain interface present.

\begin{remark}
We notice that requirement (iii) in Definition~\ref{partitiondef} implies that every connected component of $\partial\Gp_i$ either coincides with $\Gamma$ or with an interface component $\interf_j$. Moreover, for every $\interf_j$ there are precisely two subdomains $\Gp_i$ and $\Gp_h$ such that $\interf_j = \Gp_i\cap\Gp_h$. 
In particular, there may be no point in $\Omega$ where three or more subdomains meet.
\end{remark}

\begin{remark}
		On a union of subdomain interfaces $\interf:=\interf_1\cup\ldots,\cup \interf_r$ we define
		$$\VSL_T(\interf):=\left\{\solv\in\VSL(\interf)\ |\ {\solv|}_{\interf_i}\in\VSL_T(\interf_i)\ \text{for\ all}\ i=1,\ldots, r\right\},$$
where we use that $\interf_i$ is the boundary of a smooth Lipschitz domain thus $\VSL_T(\interf_i)$ being well-defined.

Furthermore, for $s\in\R$ with $2s\in\N_0$ we set
\begin{align*}
		\SHM{s}(\interf)&:=\left\{u\in\SL(\interf)\ |\ {u|}_{\interf_i}\in\SHM{s}(\interf_i)\ \text{for\ all}\ i=1,\ldots,r\right\}, \\
		\VSHM{s}(\interf)&:=\left(\SHM{s}(\interf)\right)^3, \\
		\VSHM{s}_T(\interf)&:=\VSL_T(\interf)\cap\VSHM{s}(\interf),\\
		\VSHM{-s}_T(\interf)&:=\left(\VSHM{s}_T(\interf)\right)'.
\end{align*}
\end{remark}

Sometimes the requirement of $\Gamma$ and $\partial\Gp_1,\ldots,\partial\Gp_n$ being analytic is too strong and can be relaxed to $\Gamma$ and $\partial\Gp_1,\ldots,\partial\Gp_n$ being merely smooth. We formalize this in the following definition:
\begin{definition}\label{partitiondefsmooth}
		A $\CM{\infty}$-partition $\Gp$ is a tuple $\Gp=\geom$ which consists of domains $\Omega,\Gp_1,\ldots,\Gp_n\subseteq\R^3$ satisfying conditions $(i)-(iii)$ from Definition~\ref{partitiondef} with ``analytic'' replaced by ``smooth'' in every instance.	
As before, the surfaces $\interf_1,\ldots,\interf_r$ are called subdomain interfaces and their union $\interf := \interf_1\cup\ldots\cup\interf_r$ is referred to as the union of all subdomain interfaces.
\end{definition}

Clearly, every $\ana$-partition is also a $\CM{\infty}$-partition, whereas not every $\CM{\infty}$-partition is an $\ana$-partition. 

\medskip

The notion of $\CM{\infty}$-partition is the basis for the definition of broken Sobolev spaces, which play a crucial role in this work. Let $\Gp=\geom$ be a $\CM{\infty}$-partition and assume $m\in\N_0$. Then, the broken Sobolev spaces of order $m$ are defined as
\begin{align*}
\PSHM{m}(\Gp):=\{u\in\SL(\Omega)\ \big|\ u\vert_{\Gp_i}\in\SHM{m}(\Gp_i)\ {\rm for}\ i=1,\ldots,n\} \quad{\rm and}\quad \PVSHM{m}(\Gp):=(\PSHM{m}(\Gp))^3,
\end{align*}
with their respective canonical norms
\begin{align*}
		\norm{u}{\PSHM{m}(\Gp)}^2:=\sum_{j=1}^n\norm{u}{\SHM{m}(\Gp_j)}^2\quad{\rm and} \quad \norm{\solu}{\PVSHM{m}(\Gp)}^2:=\sum_{i=1}^3\norm{u_i}{\PSHM{m}(\Gp)}^2
\end{align*}
for $\solu = (u_1,u_2,u_3)^T$. 

In addition, for $m\in\N_0$ we define the spaces of vector fields with piecewise regular curl and divergence by
\begin{align*}
		\PVHcurl{m}:=\{\solv\in\Hcurl\ \big|\ \curl\solv\in\PVSHM{m}(\Gp)\}
\end{align*}
and
\begin{align*}
\PVHdiv{m}:=\{\solv\in\Hdiv\ \big|\ \diverg\solv\in\PSHM{m}(\Gp)\},
\end{align*}
respectively.

\revision{
\begin{remark}
		The reader will immediately notice that we defined $\PVHcurl{m}$ and $\PVHdiv{m}$ in an unconventional way. Canonically, one would rather write $\solv\in\PVHcurl{m}$ if $\solv\in\PVSHM{m}(\Gp)$ and $\curl\solv\in\PVSHM{m}(\Gp)$ and analogously for $\PVHdiv{m}$. We deliberately chose this unorthodox way in order to be able to talk about (piecewise) regularity of $\curl\solu$ and $\diverg\solu$ independently of the (piecewise) regularity of $\solu$.
\end{remark}
}

%which are equipped with the norms 
%\begin{align*}
%\norm{\solu}{\PVHcurl{m}}^2:=\norm{\solu}{\VSL(\Omega)}^2+\sum_{j=1}^n\norm{\curl\solu}{\VSHM{m}(\Gp_j)}^2\quad{\rm and}\quad \norm{\solu}{\PVHdiv{m}}^2:=\norm{\solu}{\VSL(\Omega)}^2+\sum_{j=1}^n\norm{\diverg\solu}{\SHM{m}(\Gp_j)}^2,
%\end{align*}
%respectively.

The following definition clarifies the notion of piecewise smooth functions on a $\CM{\infty}$-partition.
\begin{definition}\label{pwsmooth}
		Let $\Gp=\geom$ be a $\CM{\infty}$-partition. A function $v:\Omega\rightarrow\Co$ or a vector field $\solv:\Omega\rightarrow\Co^3$ is called piecewise smooth if for all $\ell\in\N_0$ we have $v\in\PSHM{\ell}(\Gp)$ or $\solv\in\PVSHM{\ell}(\Gp)$, respectively.
\end{definition}
We mention that the terminology {\it piecewise smooth} is justified by the Sobolev embedding theorem, according to which a piecewise smooth function or vector field is indeed smooth up to the boundary of every subdomain.

\bigskip

In addition to piecewise smooth functions we will also need the concept of piecewise analytic functions. The following definition is the analog to Definition~\ref{analytictensors} in the setting of $\ana$-partitions.
\begin{definition}\label{pwanalyticdef}
		For an $\ana$-partition $\Gp=\geom$, a wavenumber $k\in\Co\setminus\{0\}$ and given constants $\omega\geq 0, M>0$, the symbols $\anapw{\Gp}$ and $\anavecpw{\Gp}$ denote the sets of all functions $u:\Omega\rightarrow\Co$ and vector fields $\solu:\Omega\rightarrow\Co^3$ which satisfy 
\begin{align*}
		&\revision{\forall\ell\in\N_0:}\ \sum_{i=1}^n\sum_{|\alpha|=\ell}\norm{\D^{\alpha}(u\vert_{\Gp_i})}{\SL(\Gp_i)}\leq \omega M^{\ell}(|k|+\ell)^{\ell},\\
		&\revision{\forall\ell\in\N_0:}\ \sum_{i=1}^n\sum_{|\alpha|=\ell}\norm{\D^{\alpha}(\solu\vert_{\Gp_i})}{\VSL(\Gp_i)}\leq \omega M^{\ell}(|k|+\ell)^{\ell},
\end{align*}
respectively.

Finally, the symbol $\anatenspw{\Gp}$ denotes the space of all tensor fields $\nu:\Omega\rightarrow\Co^{3\times 3}$ for which exist constants $\omega\geq 0, M>0$ such that 
\begin{align*}
		\revision{\forall\ell\in\N_0:}\ \sum_{m=1}^n\sum_{i,j=1}^3\sum_{|\alpha|=\ell}\norm{\D^{\alpha}(\nu_{i,j}\vert_{\Gp_m})}{\SL(\Gp_m)}\leq \omega M^{\ell}\ell^{\ell}.
\end{align*}

%We say that a tensor field $\nu\in\anatens{\Omega}$ is coercive if there exist constants $\alpha\in\Co$ with $|\alpha|=1$ and $c>0$ such that 

%\begin{align*}
%		\realpart\SCP{\alpha\nu\solz}{\solz}{}\geq c\norm{\solz}{}^2
%\end{align*}
%for all $\solz\in\Co^3$ uniformly in $\Omega$.
\end{definition}

\begin{remark}
		As in Definition~\ref{analytictensors} we track $k,\omega$ and $M$ only in the case of functions and vector fields, in the case of tensor fields we leave the definition of $\anatenspw{\Gp}$ independent of these constants, cf. Remark~\ref{kmimpl}.
\end{remark}

\subsection{General assumptions on the coefficients and main results}\label{mainresultssec}

Throughout the rest of this work we make the following coercivity assumptions concerning the given coefficients $\mu^{-1}$, $\varepsilon$ and $\zeta$:
\begin{assumption}\label{assumptionsepsmu}
		Henceforth, we assume that the coefficients $\mu^{-1}, \varepsilon:\Omega\rightarrow\Co^{3\times 3}$ and $\zeta:\Gamma\rightarrow\Co^{3\times 3}$ are coercive on $\Omega$ or $\Gamma$, respectively, and that $\zeta$ only acts in the tangent plane to $\Gamma$. More precisely, we assume that there exists constants $c>0$ and complex numbers $\alpha_{\mu^{-1}}$, $\alpha_{\varepsilon}$ and $\alpha_{\zeta}$ with $|\alpha_{\mu^{-1}}| = |\alpha_{\varepsilon}| = |\alpha_{\zeta}| = 1$ such that
		\begin{itemize}
				\item there holds
						\begin{align}\label{coercivedomain}
								\realpart\SCP{\alpha_{\mu^{-1}}\mu^{-1}\solw}{\solw}{}+\realpart\SCP{\alpha_{\varepsilon}\varepsilon\solz}{\solz}{}\geq c\norm{\solw}{}^2+c\norm{\solz}{}^2
						\end{align}
						for all $\solw,\solz\in\Co^3$ uniformly in $\Omega$,
				\item there holds 
						\begin{align}\label{coercivebdr}
								\realpart\SCP{\alpha_{\zeta}\zeta\solz}{\solz}{}\geq c\norm{\solz}{}^2
						\end{align}
						for all $\solz\in\Co^3$ uniformly on $\Gamma$,
				\item there holds 
						\begin{align}\label{zetaproperty}
(\zeta\solv)_T = \zeta\solv_T 
\end{align}
for all \revision{continuous} vector fields $\solv:\Gamma\rightarrow\Co^3$, where $\solz_T:=\soln\times(\solz\times\soln)$ denotes the tangential component of a vector field $\solz$ on $\Gamma$.
		\end{itemize}
\end{assumption}

We observe that the condition $(\zeta\solv)_T=\zeta\solv_T$ implies that multiplication with $\zeta$ maps tangent fields again to tangent fields, and for the outer normal field $\soln$ to $\Gamma$ there holds $\zeta\soln = \oldlambda\soln$ for a function $\oldlambda:\Gamma\rightarrow\Co$. If $\zeta$ is smooth or analytic, then $\oldlambda$ is also smooth or analytic.

\bigskip 

Finally, we have gathered all the tools and notations that are necessary to present the two main results of this work. 
The first main result deals with analyticity of a solution $\solu$ of \eqref{Maxwellorig} in regions where the coefficients $\mu^{-1}$ and $\varepsilon$ are analytic. For a bounded Lipschitz domain $\Omega$ we recall the spaces $\anavec{\Omega}$ and $\anatens{\Omega}$ from Definition~\ref{analytictensors}, and remember that $\Omega_2\subset\subset\Omega$ means that the closure of $\Omega_2$ is contained in $\Omega$.

\begin{theorem}\label{Mainresult1}
		Let $\Omega, \Omega_2\subseteq\R^3$ be bounded Lipschitz domains with $\Omega_2\subset\subset\Omega$ and let $\mu^{-1}, \varepsilon\in\anatens{\Omega}$ satisfy \eqref{coercivedomain}. 
		Under these assumptions, consider a vector field $\solu\in\Hcurl$ satisfying 
		\begin{align*}
				\SCP{\mu^{-1}\curl\solu}{\curl\solv}{\VSL(\Omega)}-k^2\SCP{\varepsilon\solu}{\solv}{\VSL(\Omega)}=\SCP{\solf}{\solv}{\VSL(\Omega)}
		\end{align*}
		for all $\solv\in\Hcurl$ with compact support in $\Omega$, where the wavenumber $k\in\Co\setminus\{0\}$ and the right-hand side $\solf\in\VSL(\Omega)$ are given. If $\solf\in\anaomegavec{\Omega}{\omega_{\solf}}{M}$ with $\diverg\solf\in\anaomega{\Omega}{\omega_{\diverg\solf}}{M}$ for some $\omega_{\solf},\omega_{\diverg\solf}\geq 0, M>0$, then $\solu$ satisfies 
		\begin{align}\label{Mainestimate1}
				\revision{\forall\ell\in\N_0:}\ \norm{\solu}{\VSHM{\ell}(\Omega_2)}\leq C\left(\norm{\solu}{\VSL(\Omega)}+|k|^{-1}\norm{\curl\solu}{\VSL(\Omega)}+|k|^{-2}\omega_{\solf}+|k|^{-3}\omega_{\diverg\solf}\right)L^{\ell}(\ell+|k|)^{\ell},
		\end{align}
		where $C>0$ depends only on $\Omega$ and $\Omega_2$, and  $L>0$ depends only on $\Omega$, $\Omega_2$, $\mu^{-1}$, $\varepsilon$, and $M$. In particular, $\solu$ is analytic in $\Omega_2$. 
\end{theorem}

The second main result of this work shows piecewise analyticity of a solution $\solu$ of \eqref{Maxwellorig} on a domain $\Omega$ consisting of multiple subdomains, provided that the coefficients and the given data are piecewise analytic and that the underlying geometry $\Gp=\geom$ forms an $\ana$-partition. 
We recall the spaces $\anavecpw{\Gp}$ and $\anatenspw{\Gp}$ from Definition~\ref{pwanalyticdef} as well as $\anagammavec{\Gamma}$ and $\anagamma$ from Definition~\ref{assumptionimpedance}. 

%Furthermore, in order to appropriately represent the given data in the $k$-explicit estimate~\eqref{mainres2eq} below, we abbreviate, depending on the imposed boundary conditions,
%\begin{align*}
%		&\solB_{k}:= |k|^{-3/2}C_{\solgi}\quad && {\rm in\ the\ case\ of\ } \eqref{Maxwellimpedance}, \\
%		&\solB_{k}:= |k|^{-3/2}C_{\solgn}+|k|^{-5/2}C_{\diverg\solgn}+|k|^{-5/2}C_{\solf\cdot\soln} \quad && {\rm in\ the\ case\ of\ }\eqref{Maxwellnatural},
% \\
%		&\solB_{k}:= 0 \quad && {\rm in\ the\ case\ of\ }\eqref{Maxwellessential}, 
%\end{align*}
%where the constants $C_{\solgi}, C_{\solgn}, C_{\diverg\solgn}, C_{\solf\cdot\soln}>0$ are specified in the subsequent Theorem~\ref{Mainresult2}, and where $\diverg_{\Gamman}\solgn$ denotes the surface divergence of $\solgn$ from Section~\ref{diffops}.
%With these definitions, there holds the following result:
%

\begin{theorem}\label{Mainresult2}
		Let $\Gp=\geom$ be an $\ana$-partition and suppose that $\mu^{-1},\varepsilon\in\anatenspw{\Gp}$ satisfy \eqref{coercivedomain} and that $\zeta\in\anagamma$ satisfies \eqref{coercivebdr} and \eqref{zetaproperty}. 

		Under these assumptions, consider a weak solution $\solu$ of Maxwell's equations \eqref{Maxwellorig} in the sense of Section~\ref{secweakform} below. If $\solf\in\Hdiv\cap\anaomegavecpw{\Gp}{\omega_{\solf}}{M}$ with $\diverg\solf\in\anaomegapw{\Gp}{\omega_{\diverg\solf}}{M}$ and $\solgi\in\anaomega{\Gamma}{\omega_{\solgi}}{M}$ for some $\omega_{\solf}, \omega_{\diverg\solf},\omega_{\solgi}\geq 0, M> 0$, then $\solu$ satisfies
\begin{align}\label{mainres2eq}
		\revision{\forall\ell\in\N_0:}\ \norm{\solu}{\PVSHM{\ell}(\Gp)}\leq C\left(\norm{\solu}{\VSL(\Omega)}+|k|^{-1}\norm{\curl\solu}{\VSL(\Omega)}+|k|^{-2}\omega_{\solf}+|k|^{-3}\omega_{\diverg\solf}+|k|^{-3/2}\omega_{\solgi}\right)L^{\ell}(\ell+|k|)^{\ell},
		\end{align}
		where $C>0$ depends only on $\Gp$, and $L>0$ depends only on $\Gp$, $\mu^{-1}$, $\varepsilon$, $\zeta$ and $M$.
		In particular, for any subdomain $\Gp_i$, the solution $\solu$ restricted to $\Gp_i$ is analytic up to the boundary of $\Gp_i$.

		%In the case of impedance boundary conditions \eqref{Maxwellimpedance} or natural boundary conditions \eqref{Maxwellnatural} we further suppose that 
		%\begin{itemize}
		%		\item there holds $\solgi\in\anagammavecdiscCM{\Gp}{C_{\solgi}}{M_{\solgi}}$ and that $\zeta\in\anagamma$ is coercive in the case of impedance boundary conditions \eqref{Maxwellimpedance},
		%		\item or $\solgn\in\anagammavecdiscCM{\Gp}{C_{\solgn}}{M_{\solgn}}$ in the case of natural boundary conditions \eqref{Maxwellnatural}.
		%\end{itemize}
	%Moreover, let $C_{\diverg\solf}, M_{\diverg\solf}>0$ and, in the case of natural boundary conditions, $C_{\diverg\solgn}, M_{\diverg\solgn}>0$ and $C_{\solf\cdot\soln}, M_{\solf\cdot\soln}>0$ be such that there holds $\diverg\solf\in\anapwCM{\Gp}{C_{\diverg\solf}}{M_{\diverg\solf}}$ as well as $\diverg_{\Gamma}(\solgn)\in\anafCM{\Gamma}{C_{\diverg\solgn}}{M_{\diverg\solgn}}$ and $\solf\cdot\soln\in\anafCM{\Gamma}{C_{\solf\cdot\soln}}{M_{\solf\cdot\soln}}$, respectively, where $\soln$ is the outer unit normal to $\Gamma$.

    %Then, the weak solution $\solu$ satisfies
		\end{theorem}

		We highlight that under the assumptions on $\Gp$, $\mu^{-1}$, $\varepsilon$ and $\zeta$ made by Theorem~\ref{Mainresult2}, there might not exist a weak solution of Maxwell's equations \eqref{Maxwellorig} at all. However, if there exists a (possible non-unique) solution, then Theorem~\ref{Mainresult2} asserts that any solution $\solu$ of \eqref{Maxwellorig} is piecewise analytic in $\Omega$ and satisfies \eqref{mainres2eq}.

\section{Surface differential operators and traces of $\Hcurl$}\label{diffops}

Let us assume that $\Sigma$ is a closed and orientable smooth surface and suppose that $\soln:\Sigma\rightarrow\mathbb{S}_2$ denotes a unit vector field normal to $\Sigma$. 

We follow \cite{BookMonk, BookNedelec, MaxwellMyself} and suppose that \revision{$\mathcal{U}$} is a sufficiently small tubular neighborhood around $\Sigma$. Following the notation from \cite{MaxwellImpedanceMelenk,BookNedelec}, the constant extensions (in normal direction) of a sufficiently smooth scalar function $u$ defined on $\Sigma$ to \revision{$\mathcal{U}$} is denoted by $u^*$; the surface gradient $\nabla_{\Sigma}$ and the tangential curl operator  $\overrightarrow{\curl_{\Sigma}}$ are then defined as 
\begin{align*}
\nabla_{\Sigma}u:= (\nabla u^*)\vert_{\Sigma}\quad {\rm and}\ \ \overrightarrow{\curl_{\Sigma}}u := \nabla_{\Sigma}u\times \soln
\end{align*} 
with adjoint operators $\diverg_{\Sigma}$ and $\curl_{\Sigma}$ defined by
\begin{align}\label{surfaceopdef}
	\SCP{\nabla_{\Sigma} u}{\solv}{\VSL(\Sigma)} = -\SCP{u}{\diverg_{\Sigma}\solv}{\SL(\Sigma)}\quad {\rm and}\quad \SCP{\overrightarrow{\curl_{\Sigma}}u}{\solv}{\VSL(\Sigma)}= \SCP{u}{\curl_{\Sigma}\solv}{\SL(\Sigma)}
\end{align}
for all smooth scalar functions $u$ and smooth tangent fields $\solv$. 
These operators extend to bounded operators on half-order Sobolev spaces, see \cite{TracesHcurl,MaxwellMyself}.
\begin{proposition}\label{surfaceopmappingprop}
		For any closed and orientable smooth surface $\Sigma$ and any $\ell\in\N_0$, the surface differential operators $\nabla_{\Sigma}$ and $\overrightarrow{\curl_{\Sigma}}$ extend to bounded linear operators
	$\nabla_{\Sigma}:\SHM{\ell+1/2}(\Sigma)\rightarrow \VSHM{\ell-1/2}_T(\Sigma)$ and $\overrightarrow{\curl_{\Sigma}}:\SHM{\ell+1/2}(\Sigma)\rightarrow \VSHM{\ell-1/2}_T(\Sigma)$. Moreover, the operators $\diverg_{\Sigma}$ and $\curl_{\Sigma}$ extend to bounded linear operators $\diverg_{\Sigma}:\VSHM{\ell+1/2}_T(\Sigma)\rightarrow\SHM{\ell-1/2}(\Sigma)$ and $\curl_{\Sigma}:\VSHM{\ell+1/2}_T(\Sigma)\rightarrow\SHM{\ell-1/2}(\Sigma)$.
\end{proposition}

The subsequent result \cite{MaxwellMyself} is often useful as it relates surface differential operators on the boundary of a domain $\Omega$ to traces of vector fields on $\Omega$.

\begin{proposition}
		Let $\Omega\subseteq\R^3$ be a bounded Lipschitz domain with smooth boundary $\Gamma$ and outer unit normal $\soln$. 
		Then, for any sufficiently smooth function $\psi$ on $\Omega$ there holds
	\begin{align*}
		\nabla_{\Gamma}(\psi\vert_{\Gamma}) = \soln\times(\nabla\psi\vert_{\Gamma}\times\soln) \quad{\rm and}\quad \overrightarrow{\curl_{\Gamma}}(\psi\vert_{\Gamma}) = \nabla\psi\vert_{\Gamma}\times\soln.
	\end{align*}
	
	Moreover, for any sufficiently smooth vector field $\solv$ there holds
	\begin{align}\label{surfacecurlconnect}
		\diverg_{\Gamma}\solv_t = \curl_{\Gamma}\solv_T = (\curl\solv)\vert_{\Gamma}\cdot\soln,
	\end{align}
	where $\solv_t:=\solv\vert_{\Gamma}\times\soln$ and $\solv_T:=\soln\times(\solv\vert_{\Gamma}\times\soln)$.
\end{proposition}

The next proposition \cite[Sec.~3]{TracesHcurl} turns out to be fundamental. It implies an exact sequence property of surface differential operators and asserts the existence of a Hodge decomposition for tangent fields.
\begin{proposition}\label{setting:Hodgedecomp}
		Let $\Omega\subseteq\R^3$ be a bounded Lipschitz domain with a simply connected smooth boundary $\Gamma$. Then, there holds
	\begin{align*}
			\operatorname{ker}(\curl_{\Gamma})\cap\VSL_T(\Gamma) &= \nabla_{\Gamma}\SHM{1}(\Gamma), \\
			\operatorname{ker}(\diverg_{\Gamma})\cap\VSL_T(\Gamma) &= \overrightarrow{\curl_{\Gamma}}\SHM{1}(\Gamma),\\ 
			\operatorname{ker}(\curl_{\Gamma})\cap\VSHM{-1/2}_T(\Gamma) &= \nabla_{\Gamma}\SHM{1/2}(\Gamma), \\
			\operatorname{ker}(\diverg_{\Gamma})\cap\VSHM{-1/2}_T(\Gamma) &= \overrightarrow{\curl_{\Gamma}}\SHM{1/2}(\Gamma),
	\end{align*}
	as well as 
	\begin{align}\label{setting:kersurfacegrad}
	\operatorname{ker}(\nabla_{\Gamma})\cap \SHM{1/2}(\Gamma) = \operatorname{ker}(\overrightarrow{\curl_{\Gamma}})\cap \SHM{1/2}(\Gamma) = \R.
	\end{align}
	Moreover, there holds $\VSL_T(\Gamma) = \nabla_{\Gamma}\SHM{1}(\Gamma)\oplus\overrightarrow{\curl_{\Gamma}}\SHM{1}(\Gamma)$ and this decomposition is $\VSL_T$-orthogonal, that is, $\nabla_{\Gamma}\SHM{1}(\Gamma)\perp\overrightarrow{\curl_{\Gamma}}\SHM{1}(\Gamma)$ in $\VSL_T(\Gamma)$.
\end{proposition}

\revision{
\begin{remark}
		From Section~\ref{adaptnotation} onwards we will also need surface differential operators on discs $\PZ_R$ defined by \eqref{PZMdef}. Since these discs do not form closed surfaces, the above definitions do not apply. We give a proper definition of surface differential operators on a disc $\PZ_R$ in Section~\ref{adaptnotation} below.
\end{remark}
}

\subsection{Trace spaces of $\Hcurl$}

We assume that $\Omega$ is a bounded Lipschitz domain with smooth boundary $\Gamma$ and define the auxiliary spaces 

\begin{align*}
		\Hcurlgamma{\Gamma}&:=\left\{\solu\in\VSHM{-1/2}_T(\Gamma)\ \big|\ \curl_{\Gamma}\solu\in\SHM{-1/2}(\Gamma)\right\},\\
		\Hdivgamma{\Gamma}&:=\left\{\solu\in\VSHM{-1/2}_T(\Gamma)\ \big|\ \diverg_{\Gamma}\solu\in\SHM{-1/2}(\Gamma)\right\} 
\end{align*}
with associated norms
\begin{align*}
		\norm{\solu}{\Hcurlgamma{\Gamma}}^2:=\norm{\solu}{\VSHM{-1/2}_T(\Gamma)}^2+\norm{\curl_{\Gamma}\solu}{\SHM{-1/2}(\Gamma)}^2, \\
\norm{\solu}{\Hdivgamma{\Gamma}}^2:=\norm{\solu}{\VSHM{-1/2}_T(\Gamma)}^2+\norm{\diverg_{\Gamma}\solu}{\SHM{-1/2}(\Gamma)}^2.
\end{align*}

There holds the following trace result, see \cite[Thm.~3.29]{BookMonk} or \cite[Thm.~5.4.2]{BookNedelec}.

\begin{proposition}\label{traceprop}
		Let $\Omega\subseteq\R^3$ be a bounded Lipschitz domain with boundary $\Gamma$ and outer normal vector $\soln$. 
		\revision{For vector fields $\solv\in\VCM{\infty}(\overline{\Omega})$ we consider the maps $\Pi_{T}$ and $\Pi_{t}$ defined by} 
	\begin{align*}
			\Pi_{T}\solv := \soln\times(\solv\vert_{\Gamma}\times\soln)\quad {\rm and}\quad \Pi_{t}\solv := \solv\vert_{\Gamma}\times\soln.
	\end{align*}
	These maps extend to bounded  and surjective operators $\Pi_{T}:\Hcurl\rightarrow\Hcurlgamma{\Gamma}$ and $\Pi_{t}:\Hcurl\rightarrow\Hdivgamma{\Gamma}$. Moreover, there exist bounded lifting operators $\liftcurl:\Hcurlgamma{\Gamma}\rightarrow\Hcurl$ and $\liftdiv:\Hdivgamma{\Gamma}\rightarrow\Hcurl$. 
\end{proposition}

\begin{remark}
In order to shorten notation, we abbreviate $\Pi_T\solv$ by $\solv_T$ and $\Pi_t\solv$ by $\solv_t$, respectively. 
\end{remark}

As a direct consequence of Proposition~\ref{traceprop} we have that the spaces $\Hdivgamma{\Gamma}$ and $\Hcurlgamma{\Gamma}$ are dual to each other: Indeed, for $\solv\in\Hdivgamma{\Gamma}$ and $\solw\in\Hcurlgamma{\Gamma}$ we may define a duality pairing $\dualpair{\solv}{\solw}{\Zdiv}{\Ycurl}$ by
\begin{align*}
\dualpair{\solv}{\solw}{\Zdiv}{\Ycurl}:=\SCP{\curl\liftdiv\solv}{\liftcurl\solw}{\VSL(\Omega)}-\SCP{\liftdiv\solv}{\curl\liftcurl\solw}{\VSL(\Omega)},
\end{align*}
see, e.g., \cite[Sec. 3.5.3]{BookMonk}. 
Moreover, if there holds $\solv\in\Hdivgamma{\Gamma}\cap\VSL_T(\Gamma)$ as well as $\solw\in\Hcurlgamma{\Gamma}\cap\VSL_T(\Gamma)$, the duality pairing $\dualpair{\solv}{\solw}{\Zdiv}{\Ycurl}$ coincides with $\SCP{\solv}{\solw}{\VSL_T(\Gamma)}$.

%The subsequent result \cite{MaxwellMyself} shows that higher regularity of $\solu$ yields higher regularity of its tangential trace $\solu_T$.

%\begin{lemma}
%		Let $\Omega$ be a bounded domain with smooth boundary $\Gamma$ and let $\ell\in\N_0$. Moreover, let $\solu\in\VSHM{\ell+1}(\Omega)$ satisfy $\curl\solu\in\VSHM{\ell+1}(\Omega)$. Then, there holds $\solu_T\in\VSHM{\ell+1/2}_T(\Gamma)$ as well as $\curl_{\Gamma}\solu_T\in\SHM{\ell+1/2}_T(\Gamma)$ with 
%	\begin{align*}
%		\norm{\solu_T}{\VSHM{\ell+1/2}_T(\Gamma)}+\norm{\curl_{\Gamma}\solu_T}{\SHM{\ell+1/2}_T(\Gamma)}\leq C \left(\norm{\solu}{\VSHM{\ell+1}(\Omega)}+\norm{\curl\solu}{\VSHM{\ell+1}(\Omega)}\right)
%	\end{align*}
%	for a constant $C>0$ depending only on $\ell$ and $\Omega$.
%\end{lemma}

\subsection{Variational formulation of Maxwell's equations}\label{secweakform}

We still have to clarify the notion of a weak solution of Maxwell's equations \eqref{Maxwellorig}. To that end, we assume that $\Gp=\geom$ is a $\CM{\infty}$-partition and suppose that $\mu^{-1}, \varepsilon\in\anatenspw{\Gp}$ satisfy \eqref{coercivedomain}, and that $\zeta$ is a smooth tensor field which satisfies \eqref{coercivebdr} and \eqref{zetaproperty}. 

Furthermore, let $\HXI$ denote the energy space
\begin{align*}
		\HXI:=\{\solu\in\Hcurl\ |\ \solu_T\in\VSL(\Gamma)\}.
\end{align*}
Under these assumptions, let a wavenumber $k\in\Co\setminus\{0\}$, a right-hand side $\solf\in\VSL(\Omega)$ and boundary data $\solgi\in\VSL_T(\Gamma)$ be given. Then, a vector field $\solu\in\HXI$ is called a weak solution of Maxwell's equations \eqref{Maxwellorig} if
				\begin{align*}
					\SCP{\mu^{-1}\curl \solu}{\curl\solv}{\VSL(\Omega)}-k^2\SCP{\varepsilon\solu}{\solv}{\VSL(\Omega)}-ik\SCP{\zeta\solu_T}{\solv_T}{\VSL(\Gamma)}
					= \SCP{\solf}{\solv}{\VSL(\Omega)}+\SCP{\solgi}{\solv_T}{\VSL(\Gamma)}
				\end{align*}
				for all $\solv\in\HXI$.

\bigskip

The following proposition \cite[Cor.~2.12]{MaxwellMyself} asserts that under the present assumptions on $\mu^{-1}$, $\varepsilon$ and $\zeta$, any weak solution $\solu$ of Maxwell's equations is already piecewise smooth in $\Omega$ as long as the input data $\solf$ and $\solgi$ are piecewise smooth in $\Omega$ and smooth on $\Gamma$, respectively.

\begin{proposition}\label{SmoothnessU}
		Let $\Gp=\geom$ be a $\CM{\infty}$-partition and suppose that $\mu^{-1},\varepsilon\in\anatenspw{\Gp}$ satisfy \eqref{coercivedomain}, and that $\zeta$ is smooth and satisfies \eqref{coercivebdr} and \eqref{zetaproperty}. Under these assumptions, consider a weak solution $\solu$ of Maxwell's equations \eqref{Maxwellorig} and suppose that the right-hand side $\solf\in\Hdiv$ is piecewise smooth in the sense of Definition~\ref{pwsmooth}, and that $\solgi$ is a smooth tangent field on $\Gamma$. 
		Then, $\solu$ is piecewise smooth in the sense of Definition~\ref{pwsmooth}.
\end{proposition}

In fact, the assumptions of Proposition~\ref{SmoothnessU} can be relaxed, see \cite[Cor.~2.12]{MaxwellMyself}. 
Proposition~\ref{SmoothnessU} shows that under the hypotheses of Theorem~\ref{Mainresult2}, any weak solution $\solu$ is already piecewise smooth. Hence, to infer piecewise analyticity of $\solu$ we only need to control the growth of its derivatives; a task which is carried out in the subsequent sections.

%% file: analytic/analytic.tex
\input{analytic/analytic_interior.tex}

\input{analytic/analytic_bdr.tex}

\input{analytic/analytic_normal.tex}

\input{analytic/analytic_tangential.tex}

\input{analytic/analytic_tangential_2.tex}

\input{analytic/analytic_tangential_3.tex}

\section{Proof of Theorem~\ref{Mainresult2}}\label{finalsec}

We conclude this work by giving a proof of Theorem~\ref{Mainresult2}. To that end we consider a weak solution $\solu$ of Maxwell's equations \eqref{Maxwellorig} and discuss its analyticity properties locally. 
In the interior of subdomains of an $\ana$-partition $\Gp=\geom$, analyticity of $\solu$ follows from Theorem~\ref{Mainresult1}. Near the boundary $\Gamma$ or near an interface component $\interf_j$, we employ a local flattening $\Upsilon$ to transform $\solu$ to $\PM_R$ or $\PW_R$, respectively. This leads to a transformed solution $\wu$ to which we may apply Lemma~\ref{bdryanalyticprefinal}.
The subsequent result asserts that $\solu$ and $\wu$ have the same analyticity properties. Its proof is postponed to Appendix~\ref{appendix:commest}.

\begin{lemma}\label{analyticchangevars}
		Let $\Gp=\geom$ be an $\ana$-partition in the sense of Definition~\ref{partitiondef} and suppose that $\solv$ is a piecewise smooth vector field on $\Omega$. Furthermore, let $0<\Rmin\leq R\leq\Rmax$ and assume that $\Upsilon:\PW_{2R}\rightarrow\R^3$ is a local flattening of the boundary $\Gamma$ or of an interface component $\I_{j}$. Henceforth, we write $\dom_R:=\Upsilon(\PW_R)$ for the image of $\PW_R$ under $\Upsilon$.

		In addition, let $\ell,n\in\N_0$ and let $\gamma\in\N_0^3$ be a multiindex which can be written as $\gamma=\gamma'+\gamma^*$ for multiindices $\gamma'\in\N_0^2\times\{0\}$ and $\gamma^*\in\{0\}\times\{0\}\times\N_0$ satisfying $|\gamma'|=\ell$ and $|\gamma^*|=n$.
		Then, for all $k\in\Co\setminus\{0\}$ there holds 
\begin{align*}
		\frac{1}{(\ell+n+|k|)^{\ell+n}}	\norm{\D^{\gamma}(\solv)}{\VSL\left(\dom_R\cap\Omega\right)}\leq C^{\ell+n+1}\sum_{d=0}^{\ell}\sum_{b=0}^nA^{\ell+n-d-b}\morrey{\wv}_{d,b,\PM_{2R}}
\end{align*}
if $\Upsilon$ is a local flattening of $\Gamma$, where $\wv$ is the covariant transform of $\solv$ from Definition~\ref{trafoquantitiesdef}. 

Analogously, if $\Upsilon$ is a local flattening of a subdomain interface $\I_{j}=\partial\Gp_{j_1}\cap\partial\Gp_{j_2}$, there holds 
\begin{align*}
		\frac{1}{(\ell+n+|k|)^{\ell+n}}	\left[\norm{\D^{\gamma}(\solv)}{\VSL\left(\dom_R\cap\Gp_{j_1}\right)}+\norm{\D^{\gamma}(\solv)}{\VSL\left(\dom_R\cap\Gp_{j_2}\right)}\right]\leq C^{\ell+n+1}\sum_{d=0}^{\ell}\sum_{b=0}^nA^{\ell+n-d-b}\morrey{\wv}_{d,b,\PS_{2R}}.
\end{align*}
In both cases, the constants $C,A\geq 1$ depend only on $\Rmin$, $\Rmax$, $\Gp$ and $\Upsilon$.
\end{lemma}

\bigskip 

We finish this work by providing a proof of Theorem~\ref{Mainresult2}.

\medskip

\begin{fatproofmod}{Theorem~\ref{Mainresult2}}
		First, we notice that due to Proposition~\ref{SmoothnessU} the solution $\solu$ is piecewise smooth. It remains to find wavenumber-explicit growth estimates for $\solu$ and its derivatives.

		We cover $\overline{\Omega}$ with $p\geq 1$ open balls $\PW_{r_1}(x_1),...,\PW_{r_p}(x_p)$ such that every ball $\PW_{r_m}(x_m)$ falls into one of the following three categories: 

		\begin{enumerate}
				\item 	Either $\PW_{2r_m}(x_m)\subseteq \Gp_i$ for a subdomain $\Gp_i$,
				\item   or $\PW_{2r_m}(x_m)\subseteq\Upsilon(\PW_{R_m})$ for some $R_m>0$, where $\Upsilon:\PW_{2R_m}\rightarrow\R^3$ is a local flattening of an interface component $\I_{j}$, 
				\item   or $\PW_{2r_m}(x_m)\cap \Omega\subseteq \Upsilon(\PM_{R_m})$ for some $R_m>0$, where $\Upsilon:\PW_{2R_m}\rightarrow\R^3$ is a local flattening of the boundary $\Gamma$.
		\end{enumerate}

		If $\PW_{2r_m}(x_m)\subseteq\Gp_i$ for a subdomain $\Gp_i$, then Theorem~\ref{Mainresult1} proves 
		\begin{align*}
				\norm{\D^{\gamma}(\solu)}{\VSL\left(\PW_{r_m}(x_m)\right)}\leq C\left(\norm{\solu}{\VSL(\Omega)}+|k|^{-1}\norm{\curl\solu}{\VSL(\Omega)}+|k|^{-2}\omega_{\solf}+|k|^{-3}\omega_{\diverg\solf}\right)M^{|\gamma|}\left(|\gamma|+|k|\right)^{|\gamma|}
		\end{align*}
		for all multiindices $\gamma\in\N_0^3$. 
		
		If, instead, $\Upsilon:\PW_{2R_m}\rightarrow\R^3$ is a local flattening of an interface component $\I_{j}=\partial\Gp_{j_1}\cap\partial\Gp_{j_2}$ and there holds $\PW_{2r_m}(x_m)\subseteq \Upsilon(\PW_{R_m})$, we proceed as follows: For $\ell,n\in\N_0$ assume that $\gamma\in\N_0^3$ is a multiindex such that $\gamma=\gamma'+\gamma^*$ for $\gamma'\in\N_0^2\times\{0\}$ and $\gamma^*\in\{0\}\times\{0\}\times \N$ with $|\gamma'|=\ell$ and $|\gamma^*|=n$. 
Under these assumptions and with $\mathcal{U}_{j_i}^m:=\PW_{r_m}(x_m)\cap\Gp_{j_i}$, Lemma~\ref{analyticchangevars} and Lemma~\ref{bdryanalyticprefinal} imply
		\begin{align*}
\sum_{i=1}^2\norm{\D^{\gamma}(\solu)}{\VSL\left(\mathcal{U}_{j_i}^m\right)}&\leq C^{\ell+n+1}(\ell+n+|k|)^{\ell+n}\sum_{d=0}^{\ell}\sum_{b=0}^nB^{\ell+n-d-b}\morrey{\wu}_{d,b,\PS_{R_m}}\\ 
																																			   &\leq  \left(\norm{\solu}{\VSL(\Omega)}+|k|^{-1}\norm{\curl\solu}{\VSL(\Omega)}+|k|^{-2}\omega_{\solf}+|k|^{-3}\omega_{\diverg\solf}
\right)M^{|\gamma|}\left(|\gamma|+|k|\right)^{|\gamma|}
		\end{align*}
		for appropriate constants $M,B,C> 0$.

		Analogous arguments show that if $\Upsilon$ is a local flattening of the boundary $\Gamma$ and $\PW_{2r_m}(x_m)\cap\Omega\subseteq\Upsilon(\PM_{R_m})$, then  	
\begin{align*}
		\norm{\D^{\gamma}(\solu)}{\VSL\left(\mathcal{U}^m\right)} 
		  \leq \left(\norm{\solu}{\VSL(\Omega)}+|k|^{-1}\norm{\curl\solu}{\VSL(\Omega)}+|k|^{-2}\omega_{\solf}+|k|^{-3}\omega_{\diverg\solf}+|k|^{-3/2}\omega_{\solgi}
\right)M^{|\gamma|}\left(|\gamma|+|k|\right)^{|\gamma|},
\end{align*}
where $\mathcal{U}^m:=\PW_{r_m}(x_m)\cap\Omega$.

In total, we obtain
\begin{align*}
		\sum_{i=1}^n\norm{\D^{\gamma}(\solu)}{\VSL(\Gp_i)}&\leq \sum_{i=1}^n\sum_{m=1}^p\norm{\D^{\gamma}(\solu)}{\VSL\left(\Gp_i\cap\PW_{r_m}(x_m)\right)}\\
														  &\leq C\left(\norm{\solu}{\VSL(\Omega)}+|k|^{-1}\norm{\curl\solu}{\VSL(\Omega)}+|k|^{-2}\omega_{\solf}+|k|^{-3}\omega_{\diverg\solf}+|k|^{-3/2}\omega_{\solgi}
\right)M^{|\gamma|}\left(|\gamma|+|k|\right)^{|\gamma|},
\end{align*}
for a constant $C>0$ depending only on $\Gp$ and $M>0$ depending only on $\Gp$, $\mu$, $\varepsilon$, $\zeta$ and $M$,
\revision{which} is equivalent to \eqref{mainres2eq}. Therefore, the proof is complete. 

\end{fatproofmod}

%% file: analytic/analytic_interior.tex
\section{Analytic regularity in the interior of subdomains}\label{analyticinteriorsec}

The aim of this section is to prove wavenumber-explicit analyticity of weak solutions of Maxwell's equations in the interior of subdomains. 

The outline of this Section is as follows: In \ref{commfields} we introduce wavenumber-dependent Sobolev-Morrey seminorms which will be used throughout the subsequent sections. Furthermore, inspired by \cite{GrandLivre} we provide commutator estimates that are crucial for this work as they allow us to treat non-constant coefficients $\mu^{-1}$ and $\varepsilon$. 
Subsequently, in Section~\ref{secauxest} we prove some technical estimates which are needed in Section~\ref{Mainresult1proof}, where we provide a proof for Theorem~\ref{Mainresult1}. 

\medskip

For the rest of this work we stick to the following convention: 
For $R>0$, the symbol $\PW_R(x_0)$ shall denote the ball of radius $R$ centered at $x_0$, and when $x_0$ coincides with the origin we set $\PW_R:=\PW_R(0)$ as well as $\norm{\solv}{R} := \norm{\solv}{\VSL(\PW_R(0))}$ for all smooth vector fields $\solv$. 

\begin{definition}\label{cutoffdef}
For $0<R\leq \Rmax$ and $\delta\in (0,R)$, the symbol $\chi_{R,\delta}:\R^3\rightarrow [0,1]$ denotes a smooth function with compact support in $\PW_R$ which satisfies $\chi_{R,\delta}=1$ on $\PW_{R-\delta}$ as well as 
\begin{align}\label{cutoffestimate}
		\max_{x\in\R^3}\max_{|\alpha|=1}|\D^{\alpha}\chi_{R,\delta}(x)|\leq \frac{C}{\delta}
\end{align}
for a constant $C>0$ depending only on $\Rmax$. 
\end{definition}

Let us mention that smooth cutoff functions with the properties required in Definition~\ref{cutoffdef} indeed exist; see for example \cite[Sec.~1.4]{GrandLivre} for the construction of such a function. 

%The appearance of the upper bound $R^*$ in Definition~\ref{cutoffdef} requires some comments: In the following we will often employ cutoff functions $\chi_{R_1,\delta_1}, \chi_{R_2,\delta_2},\ldots,\chi_{R_q,\delta_q}$ on balls $\PW_{R_1}, \PW_{R_2},\ldots,\PW_{R_q}$. In such a case we implicitly set $R^*:=R_q$ to ensure that the constant $C>0$ in \eqref{cutoffestimate} can be chosen uniformly for all $\chi_{R_1,\delta_1}, \chi_{R_2,\delta_2},\ldots,\chi_{R_q,\delta_q}$. 

\subsection{Commutator fields and Sobolev-Morrey seminorms}\label{commfields}

In order to ease the presentation we adopt the following notation for derivatives: For any vector field $\solv = (v_1,v_2,v_3)^T$ and any multiindex $\beta\in\N_0^3$ we define the derivative $\D^{\beta}(\solv) := (\D^{\beta}(v_1),\D^{\beta}(v_2),\D^{\beta}(v_3))^T$. With this convention, the differential operator $\D^{\beta}$ commutes with the $\curl$-operator, i.e., $\curl\D^{\beta}(\solv) = \D^{\beta}(\curl\solv)$ for all smooth vector fields $\solv$.

Motivated by \cite[Sec.~1.6]{GrandLivre}, \cite{MaxwellTomezyk} we define wavenumber-dependent Sobolev-Morrey seminorms as follows:

\begin{definition}\label{defmorrey}
		Let $0<R\leq\Rmax$, $\ell\in\N_0$ and $k\in\Co\setminus\{0\}$ be given and suppose that $\solv\in\VCM{\infty}(\overline{\PW_{\Rmax}})$. Then, with $R(\ell,k) := \frac{R}{2(\ell+|k|)}$ we set
	\begin{align*}
	 \morrey{\solv}_{\ell,R}:=\max_{0<\rho\leq R(\ell,k)}\max_{|\delta|=\ell}\rho^{\ell}\norm{\D^{\delta}\solv}{{R-(\ell+|k|)\rho}},
	\end{align*}
	and
	\begin{align*}
	\morreyone{\solv}_{\ell,R}:=\max_{0<\rho\leq R(\ell,k)}\max_{|\delta|=\ell}\rho^{\ell+1}\norm{\D^{\delta}\solv}{{R-(\ell+|k|)\rho}}.
	\end{align*}
	
	Moreover, let us suppose that $\nu:\PW_{\Rmax}\rightarrow\Co^{3\times 3}$ is a smooth tensor field on $\PW_{\Rmax}$. For all $\ell\in\N_0$ we define the wavenumber-dependent cumulative quantity $\summorreynu{\solv}{\ell,R}$ by
	\begin{align*}
	 \summorreynu{\solv}{\ell,R}:=\morrey{\solv}_{\ell,R}+|k|^{-1}\ \morrey{\nu\curl\solv}_{\ell,R}.
	\end{align*}
\end{definition}

For a smooth tensor field $\nu$, a smooth vector field $\solv$ and a multiindex $\beta\in\N_0^3$ we define the commutator field $\commute(\nu,\solv,\beta)$ by
\begin{align*}
\commute(\nu,\solv,\beta) := \D^{\beta}(\nu\solv)-\nu\D^{\beta}(\solv).
\end{align*}

The following commutator estimate will be crucial in our analysis; its proof is postponed to Appendix \ref{appendix:commest}. In essence it is a variant of \cite[Lem.~1.6.2]{GrandLivre}.

\begin{lemma}\label{commutatorestimate}
		Let $0<R\leq\Rmax$ and assume that $\nu\in\anatens{\PW_R}$, that is, there exist constants $\omega, M>0$ such that
	\begin{align*}
			\revision{\forall \ell\in\N_0:\ }	\sum_{i,j=1}^3\sum_{|\alpha|=\ell}\norm{\D^{\alpha}(\nu_{i,j})}{\SL(\PW_R)}\leq \omega M^{\ell}\ell^{\ell}.
	\end{align*}
Moreover, let $\alpha,\beta\in\N_0^3$ be arbitrary, set $b:=|\alpha|+|\beta|$ and suppose that $k\in\Co\setminus\{0\}$ is given. Then, for all $0<\rho\leq \frac{R}{2(b+|k|)}$ there holds 
	\begin{align}\label{commutone}
		\rho^b\norm{\D^{\alpha}(\commute(\nu,\solv,\beta))}{R-(b-1+|k|)\rho} \leq  C\sum_{d=0}^{b-1}K^{b-1-d}\morrey{\solv}_{d,R}
\end{align}
	for all smooth vector fields $\solv$, where the constants $C,K\geq 1$ depend only on $\omega, M$ and $\Rmax$. 
\end{lemma}

\subsection{Regularity of vector fields with piecewise regular curl and divergence}

The workhorse of this paper is a nested ball technique similar to the ones employed in \cite{GrandLivre,MaxwellTomezyk}. This technique allows us to bound higher derivatives of $\solu$ by lower derivatives on a slightly larger set. The basic requirement of such a technique is always a regularity shift result for the considered differential operator. 

\medskip 

For the Maxwell equations considered in this work, a suitable regularity shift result has been provided in \cite{MaxwellMyself}. Let us remark that similar regularity shift properties have been proved by \cite{RegularityWeber} and, very recently, also in \cite{MaxwellSpence}. 
However, in \cite{RegularityWeber,MaxwellSpence} only homogeneous boundary conditions are considered, hence we stick with the following result from \cite{MaxwellMyself}, which holds true for inhomogeneous boundary conditions on $\Gamma$ as well.
We recall that $\solu_T:=\soln\times(\solu\times\soln)$ denotes the tangential component of $\solu$, and $\solu_t:=\solu\times\soln$ is the tangential trace of $\solu$.

\begin{proposition}\label{finiteshiftprop}
		Let $\Gp$ be a $\CM{\infty}$-partition and suppose that $\nu\in\anatenspw{\Gp}$ satisfies \eqref{coercivenudomain}. Let $\solv\in\Hcurl$ satisfy $\nu\solv\in\Hdiv$ and suppose that either
		\begin{enumerate}
				\item the tangential component $\solv_T$ of $\solv$ satisfies $\solv_T=\solg_T$ for a tangent field $\solg_T\in\VSHM{1/2}_T(\Gamma)$ on $\Gamma$,
				\item or the tangential trace $\solv_t$ of $\solv$ satisfies $\solv_t = \solg_T$ for a tangent field $\solg_T\in\VSHM{1/2}_T(\Gamma)$ on $\Gamma$,
				\item or the normal trace $\nu\solv\cdot\soln$ of $\nu\solv$ satisfies $\nu\solv\cdot\soln = h$ for a function $h\in\SHM{1/2}(\Gamma)$ on $\Gamma$. 
		\end{enumerate}
		Then, $\solv\in\PVSHM{1}(\Gp)$ and 
		\begin{align*}
				\norm{\solv}{\PVSHM{1}(\Gp)}\leq C\left(\norm{\curl\solv}{\VSL(\Omega)}+\norm{\diverg\nu\solv}{\SL(\Omega)}+\norm{\solg_T}{\VSHM{1/2}_T(\Gamma)}\right)
		\end{align*}
\revision{in the cases $1.$ and $2.$, and}
		\begin{align*}
				\norm{\solv}{\PVSHM{1}(\Gp)}\leq C\left(\norm{\curl\solv}{\VSL(\Omega)}+\norm{\diverg\nu\solv}{\SL(\Omega)}+\norm{h}{\SHM{1/2}(\Gamma)}\right)
		\end{align*}
		\revision{in case $3.$ In all three cases,} the constant $C>0$ depends only on $\Gp$ and $\nu$.

		If $\nu$ even satisfies $\nu\in\anatens{\Omega}$, then $\solv\in\VSHM{1}(\Omega)$ and the above estimates hold with $\norm{\solv}{\PVSHM{1}(\Gp)}$ replaced by $\norm{\solv}{\VSHM{1}(\Omega)}$ in both instances.

\end{proposition}

We notice that Proposition~\ref{finiteshiftprop} considers only vector fields on a domain $\Omega\subseteq\R^3$. However, there is a similar result for tangent fields defined on an analytic surface \cite[Lem.~3.6]{MaxwellMyself}.

\begin{proposition}\label{shiftproptangent}
		Let $\Omega\subseteq\R^3$ be a bounded Lipschitz domain with simply connected smooth boundary $\Gamma$, and let $\lambda:\Gamma\rightarrow\Co^{3\times 3}$ be a smooth tensor field which satisfies \eqref{coercivenubdr} as well as $(\oldlambda\solz)_T = \lambda \solz_T$ for all smooth $\solz:\Gamma\rightarrow\Co^3$. 
	Then, for any tangent field $\solv\in\VSL_T(\Gamma)$ satisfying $\diverg_{\Gamma}\solv\in\SHM{-1/2}(\Gamma)$ and $\curl_{\Gamma}(\lambda\solv)\in\SHM{-1/2}(\Gamma)$, there holds $\solv\in\VSHM{1/2}_T(\Gamma)$ with
	\begin{align*}
		\norm{\solv}{\VSHM{1/2}_T(\Gamma)}\leq C\left(\norm{\diverg_{\Gamma}\solv}{\SHM{-1/2}(\Gamma)}+\norm{\curl_{\Gamma}(\lambda\solv)}{\SHM{-1/2}(\Gamma)}\right),
	\end{align*}
	where the constant $C>0$ depends only on $\lambda$ and $\Gamma$. 

	Similarly, when $\solv\in\VSL_T(\Gamma)$ satisfies $\curl_{\Gamma}\solv\in\SHM{-1/2}(\Gamma)$ and $\diverg_{\Gamma}(\lambda\solv)\in\SHM{-1/2}(\Gamma)$, there holds 
$\solv\in\VSHM{1/2}_T(\Gamma)$ with
\begin{align*}
		\norm{\solv}{\VSHM{1/2}_T(\Gamma)}\leq C\left(\norm{\diverg_{\Gamma}(\lambda\solv)}{\SHM{-1/2}(\Gamma)}+\norm{\curl_{\Gamma}\solv}{\SHM{-1/2}(\Gamma)}\right),
\end{align*}
where, again, the constant $C>0$ depends only on $\lambda$ and $\Gamma$.
\end{proposition}

Let us mention that Proposition~\ref{finiteshiftprop} and Proposition~\ref{shiftproptangent} are versions of the more general results provided by \cite[Thm.~2.6]{MaxwellMyself}, \cite[Thm~2.7]{MaxwellMyself} and \cite[Lem.~3.6]{MaxwellMyself}.
We highlight that these two propositions are crucial for the subsequent proof of Theorem~\ref{Mainresult1}.

\subsection{Auxiliary estimates for interior analyticity}\label{secauxest}

The aim of this subsection is to provide some helpful inequalities which are necessary for the proof of Theorem~\ref{Mainresult1}. For a given right-hand side $\solf$, let $\solu$ be a smooth vector field satisfying
\begin{align}\label{compactweakform}
		\SCP{\mu^{-1}\curl\solu}{\curl\solv}{\VSL(\Omega)}-k^2\SCP{\varepsilon\solu}{\solv}{\VSL(\Omega)} = \SCP{\solf}{\solv}{\VSL(\Omega)}
\end{align}
for all $\solv\in\Hcurl$ with compact support in $\Omega$. The subsequent Lemma~\ref{curldivcontrol} provides important estimates for $\solu$.

To simplify notation, we assume that the origin of the coordinate system is contained in $\Omega$; this is justified by the translation invariance of \eqref{compactweakform}. 
Furthermore, we recall the sets $\anavec{\Omega}$ and $\anatens{\Omega}$ defined in Definition~\ref{analytictensors}, and that $\chi_{R,\delta}$ denotes a smooth cutoff function with the properties from Definition~\ref{cutoffdef}.

\begin{lemma}\label{curldivcontrol}
		Let $\Omega\subseteq\R^3$ be a bounded Lipschitz domain, let $\Rmax>0$ be such that $\PW_{\Rmax}\subseteq\Omega$ and let $0<R\leq\Rmax$. Furthermore, suppose that the tensor fields $\mu^{-1},\varepsilon\in\anatens{\Omega}$ satisfy \eqref{coercivedomain}. 

		Under these assumptions, consider a vector field $\solu$ which is smooth in $\Omega$ and which satisfies \eqref{compactweakform} for some $k\in\Co\setminus\{0\}$ and a (necessarily smooth) right-hand side $\solf$. 

		\medskip

		Then, for all $\ell\in\N_0$, any multiindex $\beta\in\N_0^3$ with $|\beta|=\ell$, all $\rho\in\left(0,\frac{R}{2(\ell+1+|k|)}\right)$ \revision{and $\tau:=\ell+|k|$ there holds}
	\begin{align}\label{divestint}
	\rho^{\ell+1}\norm{\curl\chi_{R-\tau\rho,\rho}\D^{\beta}(\solu)}{R-\tau\rho}\leq C\sum_{d=0}^{\ell}A^{\ell-d}\summorrey{\solu}{d,R}
	\end{align}
	as well as
	\begin{align}\label{curlestint}
	\begin{split}
	\rho^{\ell+1}\norm{\diverg\chi_{R-\tau\rho,\rho}\ \varepsilon\D^{\beta}(\solu)}{R-\tau\rho} \leq |k|^{-2}\morreyone{\diverg\solf}_{\ell,R}+C\sum_{d=0}^{\ell}A^{\ell-d}\morrey{\solu}_{d,R},
	\end{split}	
	\end{align}
	where the constants $C,A\geq 1$ depend only on $\mu^{-1},\varepsilon$ and $\Rmax$.
	
	In addition, we have 
	\begin{align}\label{curlest2}
	|k|^{-1}\rho^{\ell+1}\norm{\curl\chi_{R-\tau\rho,\rho}\D^{\beta}(\mu^{-1}\curl\solu)}{R-\tau\rho}\leq C |k|^{-2}\morrey{\solf}_{\ell,R} +C\sum_{d=0}^{\ell}A^{\ell-d}\summorrey{\solu}{d,R}
	\end{align}
	and
	\begin{align}\label{divest2}
	\begin{split}
	|k|^{-1}\rho^{\ell+1}\norm{\diverg\chi_{R-\tau\rho,\rho}\mu\D^{\beta}(\mu^{-1}\curl\solu)}{R-\tau\rho}\leq C|k|^{-1}\sum_{d=0}^{\ell}A^{\ell-d}\morrey{\mu^{-1}\curl\solu}_{d,R},
	\end{split}
	\end{align}
	where, again, the constants $C,A\geq 1$ depend only on $\mu^{-1},\varepsilon$ and $\Rmax$.
\end{lemma}

\begin{fatproof}
	We only show \eqref{divestint}-\eqref{curlestint}, the estimates \eqref{curlest2}-\eqref{divest2} can be proved similarly.

	 The product rule and the properties of $\chi_{R-\tau\rho,\rho}$ yield 
	\begin{align}\label{tmptech2}
	\begin{split}
	\rho^{\ell+1}\norm{\curl\chi_{R-\tau\rho,\rho}\D^{\beta}(\solu)}{R-\tau\rho} &= \rho^{\ell+1}\norm{\chi_{R-\tau\rho,\rho}\curl\D^{\beta}(\solu)+\nabla\chi_{R-\tau\rho,\rho}\times\D^{\beta}(\solu)}{R-\tau\rho}\\
	&\leq \rho^{\ell+1}\norm{\curl\D^{\beta}(\solu)}{R-\tau\rho}+C\rho^{\ell}\norm{\D^{\beta}(\solu)}{\PW_{R-\tau\rho}}.
	\end{split}
	\end{align}
	By employing Lemma~\ref{commutatorestimate} and $\rho\leq R|k|^{-1}$ we get 
	\begin{align*}
			\rho^{\ell+1}\norm{\curl\D^{\beta}(\solu)}{R-\tau\rho}&\leq C\rho^{\ell}R|k|^{-1}\norm{\D^{\beta}(\mu^{-1}\curl\solu)}{R-\tau\rho}+C\rho^{\ell}R|k|^{-1}\norm{\commute(\mu,\mu^{-1}\curl\solu,\beta)}{R-(\tau-1)\rho}\\
		&\leq C|k|^{-1}\sum_{d=0}^{\ell}A^{\ell-d}\morrey{\mu^{-1}\curl\solu}_{d,R}.
	\end{align*}
	Together with \eqref{tmptech2} this proves 
	\begin{align*}
		\rho^{\ell+1}\norm{\curl\chi_{R-\tau\rho,\rho}\D^{\beta}(\solu)}{R-\tau\rho}\leq  C\sum_{d=0}^{\ell}A^{\ell-d}\summorrey{\solu}{d,R}
	\end{align*}
	and thus concludes the proof of \eqref{divestint}.
	
	\medskip 
	
	For the proof of \eqref{curlestint} we note that the product rule and the properties of $\chi_{R-\tau\rho,\rho}$ imply 
	\begin{align*}
		\rho^{\ell+1}\norm{\diverg\chi_{R-\tau\rho,\rho}\ \varepsilon\D^{\beta}(\solu)}{R-\tau\rho}&\leq \rho^{\ell+1}\norm{\diverg \varepsilon\D^{\beta}(\solu)}{R-\tau\rho}+C\rho^{\ell}\norm{\D^{\beta}(\solu)}{R-\tau\rho} \\
		&\leq \rho^{\ell+1}\norm{\diverg \D^{\beta}(\varepsilon\solu)}{R-\tau\rho}+\rho^{\ell+1}\norm{\diverg\commute(\varepsilon,\solu,\beta)}{R-\tau\rho}\\
		&\hskip 5cm +C\rho^{\ell}\norm{\D^{\beta}(\solu)}{R-\tau\rho}.
	\end{align*}
	From \eqref{compactweakform} we infer $-k^2\diverg(\varepsilon\solu) = \diverg\solf$ in $\PW_R$, and Lemma~\ref{commutatorestimate} proves 
	\begin{align*}
		\rho^{\ell+1}\norm{\diverg\commute(\varepsilon,\solu,\beta)}{R-\tau\rho}\leq C\sum_{d=0}^{\ell}K^{\ell-d}\morrey{\solu}_{d,R}.
	\end{align*}
	In total we obtain 
	\begin{align*}
			\rho^{\ell+1}\norm{\diverg\chi_{R-\tau\rho,\rho}\ \varepsilon\D^{\beta}(\solu)}{R-\tau\rho}\leq |k|^{-2}\morreyone{\diverg\solf}_{\revision{\ell},R}+C\sum_{d=0}^{\ell}A^{\ell-d}\morrey{\solu}_{d,R},
	\end{align*}
	which finishes the proof of \eqref{curlestint}.
\end{fatproof}

\subsection{Proof of Theorem~\ref{Mainresult1}}\label{Mainresult1proof}

The auxiliary estimates derived in Lemma~\ref{curldivcontrol} allow for a rather elegant proof of Theorem~\ref{Mainresult1}. The key is to control the growth of the cumulative quantity $\summorrey{\solu}{\ell+1,R}$ from Definition~\ref{defmorrey}. We start by proving a recursive estimate for $\summorrey{\solu}{\ell+1,R}$ in terms of the right-hand side $\solf$ and lower order derivatives of both $\solu$ and $\curl\solu$.

\begin{lemma}\label{analyticinteriorprelemma}
		Let $\Omega\subseteq\R^3$ be a bounded Lipschitz domain, let $\Rmax>0$ be such that $\PW_{\Rmax}\subseteq\Omega$ and let $0<R\leq\Rmax$. 
		Furthermore, suppose that $\mu^{-1},\varepsilon\in\anatens{\Omega}$ satisfy \eqref{coercivedomain}. 

		Under these assumptions, consider a vector field $\solu$ which is smooth in $\Omega$ and which satisfies \eqref{compactweakform} for some $k\in\Co\setminus\{0\}$ and a (necessarily smooth) right-hand side $\solf$. 

Then, for all $\ell\in\N_0$ there holds
	\begin{align}\label{summorreyiteration}
		\summorrey{\solu}{\ell+1,R}\leq C|k|^{-2}\left[\morrey{\solf}_{\ell,R}+\morreyone{\diverg\solf}_{\ell,R}\right] +C\sum_{d=0}^{\ell}A^{\ell-d}\summorrey{\solu}{d,R},
	\end{align}
	where the constants $C,A\geq 1$ depend only on $\Rmax$, $\mu^{-1}$ and $\varepsilon$.
\end{lemma}

\begin{fatproof}
		The proof is divided \revision{into} two steps. In the first step we derive an estimate for $\morrey{\solu}_{\ell+1,R}$ and in the second step we prove an upper bound for $|k|^{-1}\morrey{\mu^{-1}\curl\solu}_{\ell+1,R}$.
	
	\medskip
	
	{\bf Step 1:} Let $\ell\in\N_0$ be fixed, choose a multiindex $\delta\in\N_0^3$ with $|\delta|=\ell+1$ and choose $\rho\in\left(0, \frac{R}{2(\ell+1+|k|)}\right)$. By definition, there exists a decomposition $\delta = \alpha+\beta$ with multiindices $\alpha,\beta\in\N_0^3$ with $|\alpha|=1$ and $|\beta| = \ell$. We define $\tau := \ell+|k|$ and notice that $\chi_{R-\tau\rho,\rho}\equiv 1$ in $\PW_{R-(\tau+1)\rho}$ implies 
	\begin{align}\label{fundestint1}
	\begin{split}
	\rho^{\ell+1}\norm{\D^{\delta}(\solu)}{R-(\tau+1)\rho}&\leq \rho^{\ell+1}\norm{\D^{\alpha}\chi_{R-\tau\rho,\rho}\D^{\beta}(\solu)}{R-\tau\rho} \\ &\leq \rho^{\ell+1}\norm{\chi_{R-\tau\rho,\rho}\D^{\beta}(\solu)}{\VSHM{1}(\PW_{R-\tau\rho})}.
	\end{split}	
	\end{align}
	
	According to Proposition~\ref{finiteshiftprop} we have
	\begin{align}\label{fundestint2}
		\norm{\chi_{R-\tau\rho,\rho}\D^{\beta}(\solu)}{\VSHM{1}(\PW_{R-\tau\rho})}\leq C\left(\norm{\curl\chi_{R-\tau\rho,\rho}\D^{\beta}(\solu)}{R-\tau\rho}+\norm{\diverg\chi_{R-\tau\rho,\rho}\varepsilon\D^{\beta}(\solu)}{R-\tau\rho}\right).
	\end{align}
	From Lemma \ref{curldivcontrol} we know
	\begin{align*}
		\rho^{\ell+1}\norm{\curl\chi_{R-\tau\rho,\rho}\D^{\beta}(\solu)}{R-\tau\rho}\leq C\sum_{d=0}^{\ell}A^{\ell-d}\summorrey{\solu}{d,R}
	\end{align*}
	and 
	\begin{align*}
		\begin{split}
		\rho^{\ell+1}\norm{\diverg\chi_{R-\tau\rho,\rho}\ \varepsilon\D^{\beta}(\solu)}{R-\tau\rho} \leq |k|^{-2}\morreyone{\diverg\solf}_{\ell,R}+C\sum_{d=0}^{\ell}A^{\ell-d}\morrey{\solu}_{d,R}.
		\end{split}
	\end{align*}

	Together with \eqref{fundestint1} and \eqref{fundestint2}, these two inequalities imply
	
	\begin{align}\label{estimateanalu}
		\morrey{\solu}_{\ell+1,R}\leq C|k|^{-2}\morreyone{\diverg\solf}_{\ell,R}+C\sum_{d=0}^{\ell}A^{\ell-d}\summorrey{\solu}{d,R}
	\end{align}
	for appropriate constants $C,A\geq 1$.
	
	\medskip 
	
	{\bf Step 2:} Let $\delta$ again be a multiindex of order $\ell+1$ and decompose $\delta = \alpha+\beta$ with $|\alpha|=1$ and $|\beta|=\ell$. We observe that for all $\rho\in\left(0,\frac{R}{2(\ell+1+|k|)}\right)$ there holds 
	\begin{align}\label{fundanal3}
		|k|^{-1}\rho^{\ell+1}\norm{\D^{\delta}(\mu^{-1}\curl\solu)}{R-(\tau+1)\rho}\leq |k|^{-1}\rho^{\ell+1}\norm{\chi_{R-\tau\rho,\rho}\D^{\beta}(\mu^{-1}\curl\solu)}{\VSHM{1}(\PW_{R-\tau\rho})}.
	\end{align}
	By applying Proposition~\ref{finiteshiftprop} to the function $\chi_{R-\tau\rho,\rho}\D^{\beta}(\mu^{-1}\curl\solu)$  we get
	
	\begin{align}\label{fundanal4}
	\begin{split}
	\norm{\chi_{R-\tau\rho,\rho}\D^{\beta}(\mu^{-1}\curl\solu)}{\VSHM{1}(\PW_{R-\tau\rho})}\leq C&\norm{\curl\chi_{R-\tau\rho,\rho}\D^{\beta}(\mu^{-1}\curl\solu)}{R-\tau\rho} \\
	&\hskip 1cm + C\norm{\diverg\chi_{R-\tau\rho,\rho}\mu\D^{\beta}(\mu^{-1}\curl\solu)}{R-\tau\rho}.
	\end{split}
	\end{align}
	According to Lemma \ref{curldivcontrol} we have 
	\begin{align*}
		|k|^{-1}\rho^{\ell+1}\norm{\curl\chi_{R-\tau\rho,\rho}\D^{\beta}(\mu^{-1}\curl\solu)}{R-\tau\rho}\leq C |k|^{-2}\morrey{\solf}_{\ell,R} +C\sum_{d=0}^{\ell}A^{\ell-d}\summorrey{\solu}{d,R}
	\end{align*}
	and
	\begin{align*}
		\begin{split}
			|k|^{-1}\rho^{\ell+1}\norm{\diverg\chi_{R-\tau\rho,\rho}\mu\D^{\beta}(\mu^{-1}\curl\solu)}{R-\tau\rho}\leq C|k|^{-1}\sum_{d=0}^{\ell}A^{\ell-d}\morrey{\mu^{-1}\curl\solu}_{d,R}.
		\end{split}
	\end{align*}
	Together with \eqref{fundanal3} and \eqref{fundanal4}, these two inequalities prove 
	\begin{align*}
		|k|^{-1}\morrey{\mu^{-1}\curl\solu}_{\ell+1,R}\leq C|k|^{-2}\morrey{\solf}_{\ell,R}+ C\sum_{d=0}^{\ell}A^{\ell-d}\summorrey{\solu}{d,R},
	\end{align*} 
	and in combination with \eqref{estimateanalu} this proves \eqref{summorreyiteration} and thus concludes the proof.
	
\end{fatproof}

Lemma~\ref{analyticinteriorprelemma} provides a recursion formula for the cumulative quantity $\summorrey{\solu}{\ell,R}$. The next step is to resolve the obtained recursive relation and to derive an estimate for $\summorrey{\solu}{\ell,R}$ solely in terms of $\summorrey{\solu}{0,R}$ and derivatives of the right-hand side $\solf$. Doing this is the purpose of the following lemma. The proof exploits Lemma~\ref{analyticinteriorprelemma} and mainly consists of an induction argument and some computations; it is postponed to Appendix \ref{appendix:commest}.

\begin{lemma}\label{analyticinterior}
		Let $\Omega\subseteq\R^3$ be a bounded Lipschitz domain, let $\Rmax>0$ such that $\PW_{\Rmax}\subseteq\Omega$. 
		Furthermore, suppose that $\mu^{-1},\varepsilon\in\anatens{\Omega}$ satisfy \eqref{coercivedomain}.

		Under these assumptions, consider a vector field $\solu$ which is smooth in $\Omega$ and which satisfies \eqref{compactweakform} for some $k\in\Co\setminus\{0\}$ and a (necessarily smooth) right-hand side $\solf$. 

		\medskip 

Then, there exists a constant $A\geq 1$ depending only on $\Rmax$, $\mu^{-1}$ and $\varepsilon$, such that for all $\ell\in\N_0$ there holds
	\begin{align}\label{estimateomega}
		\summorrey{\solu}{\ell+1,R}\leq A^{\ell+1}\summorrey{\solu}{0,R}+|k|^{-2}\sum_{d=0}^{\ell}A^{\ell+1-d}\left[\morrey{\solf}_{d,R}+\morreyone{\diverg\solf}_{d,R}\right].
	\end{align}
	As a consequence, if $0<\Rmin\leq R\leq \Rmax$ then for all $\ell\in\N_0$ we have
	\begin{align}\label{analyticityfinal}
	\begin{split}
			\frac{P^{-\ell-1}}{(\ell+1+|k|)^{\ell+1}}\ \max_{|\delta|=\ell+1}\norm{\D^{\delta}(\solu)}{R/2}&\leq \summorrey{\solu}{0,R}+|k|^{-2}\sum_{d=0}^{\ell}P^{-d}\frac{1}{(d+|k|)^d}\ \max_{|\beta|=d}\norm{\D^{\beta}(\solf)}{R} \\
	& \hskip 2.2cm +|k|^{-3}\sum_{d=0}^{\ell}P^{-d}\frac{1}{(d+|k|)^d}\ \max_{|\beta|=d}\norm{\D^{\beta}(\diverg\solf)}{R},
	\end{split}
	\end{align}
	where the constant $P\geq 1$ depends only on $\Rmin$, $\Rmax$, $\mu^{-1}$ and $\varepsilon$.
\end{lemma}

We conclude this section by proving Theorem~\ref{Mainresult1}. We recall that for given $0<R\leq R^*$ and $\delta\in (0,R)$ we assume that $\chi_{R,\delta}:\R^3\rightarrow[0,1]$ is a smooth cutoff function satisfying the properties from Definition~\ref{cutoffdef}.

\bigskip

\begin{fatproofmod}{Theorem~\ref{Mainresult1}}
		Under the hypotheses of Theorem~\ref{Mainresult1}, let $\solu\in\Hcurl$ satisfy \eqref{compactweakform} for a right-hand side $\solf\in\anavec{\Omega}$. 
		First, we need to show that $\solu$ is smooth in $\Omega$. To that end, let $x_0\in\Omega$ and $R>0$ be such that $\PW_{2R}(x_0)\subseteq\Omega$. Due to the translation invariance of \eqref{compactweakform} we may assume that $x_0$ coincides with the origin, hence in order to prove that $\solu$ is smooth in $\Omega$ it suffices to show that $\chi_{2R,R}\solu$ is smooth on the ball $\PW_{2R}$.

		Indeed, from \eqref{compactweakform} we infer the relations $-k^2\diverg\varepsilon\solu = \diverg\solf$ and $\curl\mu^{-1}\curl\solu = \solf+k^2\varepsilon\solu$. We can exploit these relations and use Proposition~\ref{finiteshiftprop}, Proposition~\ref{shiftproptangent} and similar techniques as in \cite[Proof of Thm.~2.10]{MaxwellMyself} to conclude $\chi_{2R,R}\solu\in\VSHM{\ell}(\PW_{2R})$ for all $\ell\in\N_0$. The Sobolev embedding theorem shows that $\solu$ is smooth in $\Omega$.

	\bigskip	

	It remains to prove \eqref{Mainestimate1}. We notice that it suffices to consider the case $\Omega_2 = \PW_R(x_0)$ for some $R>0$ and $x_0\in \Omega$ which satisfy $\PW_{2R}(x_0)\subseteq \Omega$; the general case follows by covering $\overline{\Omega_2}$ with finitely many such balls. Moreover, by translation invariance of \eqref{compactweakform} we may without loss of generality assume $x_0=0$, thus we restrict ourselves to the case $\Omega_2 = \PW_R$ for some $R>0$ satisfying $\PW_{2R}\subseteq\Omega$.

	Under these assumptions and since we showed that $\solu$ is smooth in $\Omega$, Lemma \ref{analyticinterior} is applicable and yields the estimate \eqref{analyticityfinal}. Since by assumption $\solf\in\anaomegavec{\Omega}{\omega_{\solf}}{M}$ and $\diverg\solf\in\anaomega{\Omega}{\omega_{\diverg\solf}}{M}$, there holds
	\begin{align*}
			\norm{\D^{\gamma}(\solf)}{2R}\leq \omega_{\solf}M^{|\gamma|}(|k|+|\gamma|)^{|\gamma|}\quad {\rm and}\quad \norm{\D^{\gamma}(\diverg \solf)}{2R}\leq \omega_{\diverg\solf}M^{|\gamma|}(|k|+|\gamma|)^{|\gamma|}
	\end{align*}
	for all mutliindices $\gamma\in\N_0^3$.
	With this, \eqref{analyticityfinal} implies
	\begin{align}\label{tmpana}
			\frac{1}{(\ell+|k|)^{\ell}}\norm{\D^{\gamma}(\solu)}{R} \leq L^{\ell}\left(|k|^{-2}\omega_{\solf}+|k|^{-3}\omega_{\diverg\solf}+\summorrey{\solu}{0,2R}\right)
	\end{align}
	for all multiindices $\gamma$ with $|\gamma| = \ell$, where $L>0$ depends only on $M$ and the constant $P$ from \eqref{analyticityfinal}. This concludes the proof.
\end{fatproofmod}

%% file: analytic/analytic_bdr.tex
\section{Local flattenings of boundaries and subdomain interfaces}\label{flatteningsec}

It remains to prove the second main result of this work, namely Theorem~\ref{Mainresult2}. The main idea is to locally flatten the boundary $\Gamma$ and the subdomain interfaces $\interf_1,\ldots,\interf_r$ to the disc $\PZ_R:=\PW_R\cap\{x_3=0\}$ and to subsequently derive estimates for tangential and normal derivatives of $\solu$. The aim of this intermediate section is to describe the flattening process in more detail.

\subsection{Definition and properties of local flattenings}

In order to properly introduce the type of transformations that we rely upon in the following, we adopt the subsequent notation: As previously, for $R>0$ and $x_0\in\R^3$, the symbol $\PW_R(x_0)$ shall denote the open ball with radius $R$ centered at $x_0$. If $x_0$ coincides with the origin, we write $\PW_R$ for short. Moreover, for $R>0$ we set 
\begin{align}\label{halfballdef}
	\PP_R:=\PW_R\cap\{x_3>0\},\quad \PM_R:=\PW_R\cap\{x_3<0\},\quad{\rm and}\quad   \PZ_R:=\PW_R\cap\{x_3=0\}.
\end{align}

\revision{Note that this definition of $\PZ_R$ agrees with \eqref{PZMdef}.} Let $\Sigma$ be a closed, orientable and analytic\footnote{Most of the definitions and results of this section can analogously be formulated for smooth surfaces. However, throughout this work we will apply local flattenings only to analytic surfaces.} surface with normal unit vector $\soln$, and let 
\begin{align*}
\sigma:\{u,v\in\R\ |\ u^2+v^2 < 1\}\rightarrow\Sigma
\end{align*}
be an analytic chart. Then, for $R>0$ sufficiently small, the map $\Upsilon:\PW_R\rightarrow\R^3$ defined by 
\begin{align}\label{framebased}
\Upsilon(u,v,s):=\sigma(u,v)+s\soln(\sigma(u,v))
\end{align}
is an analytic diffeomorphism which maps $\PW_R$ to an open neighbourhood of $\sigma(0,0)\in\Sigma$. 
Henceforth, we write $F$ for the Jacobian of $\Upsilon$ and set $J:=\operatorname{det}F$.

\begin{definition}\label{framebaseddef}
A local flattening of $\Sigma$ is a coordinate transformation $\Upsilon$ of the form \eqref{framebased} which satisfies the condition $J>0$.
\end{definition}

Note that by interchanging $u$ and $v$ in \eqref{framebased} we change the sign of $J$, i.e., the requirement $J>0$ in Definition~\ref{framebaseddef} is no real restriction. 
The following lemma captures an important property of local flattenings:

\begin{lemma}\label{FTF}
		Let $\Sigma$ be a closed, orientable and analytic surface with unit normal vector $\soln$ and let $\Upsilon:\PW_R\rightarrow\R^3$ be a local flattening of $\Sigma$. Then, the Jacobian $F$ of $\Upsilon$ satisfies
	\begin{align*}
		F^TF = \begin{pmatrix}
		E & H & 0 \\
		H & G & 0 \\
		0 & 0& 1
\end{pmatrix} \quad \text{on}\ \PZ_R,
	\end{align*}
where $E,H$ and $G$ are the coefficients of the first fundamental form of $\Sigma$. Moreover, we have $F^T(\soln\circ\Upsilon) = \sole_3$, where $\sole_3 := (0,0,1)^T$. 
\end{lemma}
\begin{fatproof}
	Follows from elementary manipulations and the observation that on $\Sigma$ (i.e., when $s=0)$ the derivatives $\Upsilon_u$ and $\Upsilon_v$ are tangent to $\Sigma$ and thus normal to $\soln$.
\end{fatproof}
\begin{remark}\label{lambdainvariance}
		A vector field $\solt:\PZ_R\rightarrow\Co^3$ is called a tangent field on $\PZ_R$ if its third component vanishes.
	Note that Lemma \ref{FTF} implies that for any tangent field $\solt$ on $\PZ_R$ the vector fields $F^TF\solt$ and $F^{-1}F^{-T}\solt$ are again tangential. It also shows that $\sole_3$ is invariant under $F^TF$ and $F^{-1}F^{-T}$.
\end{remark}

Subsequently, we will assume that $\Gp=\geom$ is an $\ana$-partition in the sense of Definition~\ref{partitiondef}, and the surface $\Sigma$ will either be $\Gamma$ or a subdomain interface $\interf_j$. By construction, if $\Upsilon:\PW_R\rightarrow\R^3$ is a local flattening of $\Gamma$, we have that $\Upsilon(\PM_R)\subseteq\Omega$. Similarly, if $\Upsilon:\PW_R\rightarrow\R^3$ is a local flattening of $\interf_j=\Gp_{j_1}\cap\Gp_{j_2}$ we have that $\Upsilon(\PM_R)\subseteq \Gp_{j_1}$ and $\Upsilon(\PP_R)\subseteq \Gp_{j_2}$ or vice versa.

\medskip

So far, we discussed how we transform geometries, the follow-up question is how differential operators behave under these coordinate transformations. The next result (\cite[Lem.~14]{MaxwellSchoeberl} or \cite[Sec.~3.9]{BookMonk}) provides the answer.
\begin{lemma}\label{schöberltrafo}
	Let $\Omega_1$ and $\Omega_2$ be domains in $\R^3$ and let $\phi:\Omega_1\rightarrow \Omega_2$ be a smooth diffeomorphism, let $F$ denote the Jacobian of $\phi$ and let $J:=\det F$. Then, for all smooth vector fields $\solv$ there holds 
	\begin{align*}
		\curl\left[F^T(\solv\circ\phi)\right] = JF^{-1}\left[(\curl\solv)\circ\phi\right], \quad{\rm as\ well\ as}\
		\diverg\left[JF^{-1}(\solv\circ\phi)\right] = J\left[(\diverg\solv)\circ\phi\right].
	\end{align*}
	The mapping $\solv\mapsto F^T(\solv\circ\phi)$ is called covariant transformation, whereas $\solv\mapsto JF^{-1}(\solv\circ\phi)$ is called Piola transformation.
\end{lemma}

With this result we finally have the tools to discuss the behaviour of Maxwell's equations \eqref{Maxwellorig} under local flattenings. Subsequently, we prove invariance of Maxwell's equations under local flattenings; this invariance will be the starting point of deriving analytic regularity of a solution $\solu$ up to the boundary and up to subdomain interfaces.

\subsection{Transformation invariance of Maxwell's equations}

To simplify further calculations, we transform the coefficients $\mu^{-1}$, $\varepsilon$ and $\zeta$ according to the subsequent definition: 
%We notice that in the case of $\Upsilon:\PW_R\rightarrow\R^3$ being a local flattening of $\Gamma$ and $\mu^{-1},\varepsilon\in\anatenspw{\Gp}$, the vector fields $\mu^{-1}\circ\Upsilon$ and $\varepsilon\circ\Upsilon$ are formally only defined on $\PM_R$. We circumvent this issue by noticing that $\mu^{-1},\varepsilon\in\anatenspw{\Gp}$ can be extended analytically to a small region outside of $\Omega$, hence by implicitly assuming $R>0$ sufficiently small we may consider $\mu^{-1}\circ\Upsilon$ and $\varepsilon\circ\Upsilon$ to be well-defined on $\PW_R$.

\begin{definition}\label{trafoquantitiesdef}
		Let $\Gp=\geom$ be an $\ana$-partition, suppose that $\mu^{-1},\varepsilon\in\anatenspw{\Gp}$ satisfy \eqref{coercivedomain} and assume that $\zeta\in\anagamma$ satisfies \eqref{coercivebdr} and \eqref{zetaproperty}. 
Furthermore, assume that $\Upsilon:\PW_R\rightarrow\R^3$ is a local flattening of the boundary $\Gamma$. With $F$ denoting the Jacobian of $\Upsilon$ and $J:=\operatorname{det}F$ we define the tensor fields $\wmu^{-1}$, $\weps$ and $\wzeta$ as
		\begin{align*}
				\wmu^{-1}:=J^{-1}F^T(\mu^{-1}\circ\Upsilon\vert_{\PM_R})F, \quad \weps:=JF^{-1}(\varepsilon\circ\Upsilon\vert_{\PM_R})F^{-T} \quad{\rm and}\quad \wzeta:=JF^{-1}(\zeta\circ\Upsilon\vert_{\PZ_R})F^{-T}.
	\end{align*}

	In addition, for any vector field $\solv$ on $\Omega$  we set $\tildf{\solv}:= F^T(\solv\circ\Upsilon\vert_{\PM_R})$ and $\hatf{\solv}:=JF^{-1}(\solv\circ\Upsilon\vert_{\PM_R})$. Analogously, for any tangent field $\solg$ on $\Gamma$ we define $\tildf{\solg}:=F^T(\solg\circ\Upsilon\vert_{\PZ_R})$ and $\hatf{\solg}:=JF^{-1}(\solg\circ\Upsilon\vert_{\PZ_R})$, respectively.

	\medskip

If, instead of being a local flattening of $\Gamma$, the transformation $\Upsilon:\PW_R\rightarrow\R^3$ is a local flattening of a subdomain interface $\interf_j$, we define
\begin{align*}
		\wmu^{-1}:=J^{-1}F^T(\mu^{-1}\circ\Upsilon)F\quad {\rm and}\ \quad \weps:=JF^{-1}(\varepsilon\circ\Upsilon)F^{-T},
	\end{align*}
	as well as $\tildf{\solv}:= F^T(\solv\circ\Upsilon)$ and $\hatf{\solv}:=JF^{-1}(\solv\circ\Upsilon)$ for all tensor fields $\mu^{-1},\varepsilon\in\anatenspw{\Gp}$ and all vector fields $\solv$ on $\Omega$.

\end{definition}

From Lemma~\ref{FTF} we notice that for a local flattening $\Upsilon:\PW_R\rightarrow\R^3$ of $\Gamma$, the outer unit normal $\soln$ to $\Gamma$ satisfies $\tildf{\soln}=\sole_3$, where $\sole_3=(0,0,1)^T$. This fact will be useful later.
At this point it is natural to ask which properties the transformed coefficients $\wmu^{-1}$, $\weps$ and $\wzeta$ inherit from $\wmu^{-1}$, $\weps$ and $\wzeta$. This question will be answered in the subsequent Section~\ref{adaptnotation}, since for now we are more interested in the following result, which shows that Maxwell's equations \eqref{Maxwellorig} are invariant under transformations by local flattenings.

\begin{lemma}\label{trafolemma}
		Let $\Gp=\geom$ be an $\ana$-partition, suppose that $\mu^{-1},\varepsilon\in\anatenspw{\Gp}$ satisfy \eqref{coercivedomain} and let $\zeta\in\anagamma$ satisfy \eqref{coercivebdr} and \eqref{zetaproperty}. 

		Furthermore, suppose that $\solu$ is a weak solution of $\eqref{Maxwellorig}$ corresponding to a right-hand side $\solf$ and boundary data $\solgi$. 

		Finally, suppose that $\Upsilon:\PW_R\rightarrow\R^3$ is a local flattening of the boundary $\Gamma$. Then, with the transformed quantities $\wmu^{-1}$, $\weps$ and $\wzeta$ from Definition~\ref{trafoquantitiesdef}, the transformed solution $\wu$ satisfies
	\begin{align}\label{trafosystem}
			\curl\wmu^{-1} \curl\wu-k^2\weps\wu = \wf\quad\text{in}\ \PM_R,
	\end{align}
	as well as the boundary condition 
	\begin{align}\label{transformedbdry}
			\wmu^{-1}\curl\wu\times\sole_3-ik\wzeta\tildf{\solu}_T = \solgicheck\quad\text{on}\ \PZ_R.
	\end{align}
	
	If, instead of being a local flattening of $\Gamma$, the transformation $\Upsilon:\PW_R\rightarrow\R^3$ is a local flattening of a subdomain interface $\I_j$, the transformed solution $\tildf{\solu}$ satisfies 
	\begin{align}\label{trafosystem2}
			\curl\wmu^{-1} \curl\tildf{\solu}-k^2\weps\tildf{\solu} &= \hatf{\solf}\quad{\rm in}\ \PP_R\cup\PM_R.
	\end{align}
	This shows that the Maxwell equations \eqref{Maxwellorig} are invariant under transformations by local flattenings.
\end{lemma}

\begin{fatproof}
		\revision{Let $\solv:=\mu^{-1}\curl\solu$. We observe that }any weak solution $\solu$ of \eqref{Maxwellorig} satisfies
	\begin{align*}
	\curl\solv - k^2\varepsilon\solu = \solf \quad{\rm in}\ \Omega
	\end{align*}
	with the boundary condition 
			$\solv\times\soln-ik\zeta\solu_T = \solgi$ on $\Gamma$. According to Lemma \ref{schöberltrafo}, the transformed vector fields $\tildf{\solv}$ and $\tildf{\solu}$ satisfy 
	\begin{align*}
		\curl\tildf{\solv} = JF^{-1}\left[(\curl\solv)\circ\Upsilon\right] = JF^{-1}\left[(k^2\varepsilon\solu+\solf)\circ\Upsilon\right] = JF^{-1}\left[k^2(\varepsilon\circ\Upsilon)F^{-T}\tildf{\solu}+\solf\circ\Upsilon\right],
	\end{align*}
	that is,
		$\curl\tildf{\solv} -k^2\weps\tildf{\solu}=\hatf{\solf}$.
	Moreover, Lemma \ref{schöberltrafo} shows 
	\begin{align*}
		\tildf{\solv} = F^T\left[(\mu^{-1}\curl\solu)\circ\Upsilon\right] = F^T\left[(\mu^{-1}\circ\Upsilon)J^{-1}F\curl\tildf{\solu}\right] = \wmu^{-1}\curl\tildf{\solu},
	\end{align*}
	thus proving \eqref{trafosystem} and \eqref{trafosystem2}. 

	If $\Upsilon$ is a local flattening of the boundary $\Gamma$, \revision{the formula $Ma\times Mb = \operatorname{det}(M)M^{-T}(a\times b)$ yields}
	\begin{align*}
		\wmu^{-1}\curl\tildf{\solu}\times\sole_3 &= JF^{-1}\left[\mu^{-1}(\curl\solu)\times\soln\right]\circ\Upsilon \\
		&=ikJF^{-1}((\zeta\solu)_T\circ\Upsilon)+\solgicheck,
	\end{align*}
	and in combination with the identity $\sola\times(\solb\times \solc) = (\solc \cdot \sola)\solb-(\sola\cdot \solb)\solc$ and Lemma~\ref{FTF} we get
	\begin{align*}
		JF^{-1}((\zeta\solu)_T\circ\Upsilon) &= J\left[(F^TF)^{-1}\tildf{\zeta\solu}-\left[(F^TF)^{-1}\tildf{\zeta\solu}\cdot\tildf{\soln}\right] (F^TF)^{-1}\tildf{\soln}\right] \\
		&= J\left[(F^TF)^{-1}\tildf{\zeta\solu}-\left[(F^TF)^{-1}\tildf{\zeta\solu}\cdot\sole_3\right] \sole_3\right] \\
		&= J(F^TF)^{-1}\tildf{\zeta\solu}_T.
	\end{align*}
	Finally, some calculations show $J(F^TF)^{-1}\tildf{\zeta\solu}_T = (\wzeta\tildf{\solu})_T=\wzeta\tildf{\solu}_T$,
	which proves \eqref{transformedbdry}.

\end{fatproof}

\section{Vector- and tensor fields on half-balls and discs}\label{adaptnotation}

In the previous section we introduced the notion of local flattenings, which we will use to change the underlying geometry given by an $\ana$-partition $\Gp$ to a geometry governed by balls, half-balls and discs. As before, for $R>0$ we define half-balls $\PM_R$, $\PP_R$ and discs $\PZ_R$ by
\begin{align*}
		\PM_R := \PW_R\cap\{x_3<0\},\quad \PP_R:=\PW_R\cap\{x_3>0\}\quad{\rm and} \quad \PZ_R:=\PW_R\cap\{x_3=0\},
\end{align*}
where $\PW_R\subseteq\R^3$ denotes the open ball with radius $R$ centered at the origin.

Changing the underlying geometry requires certain adaptions concerning our notation. For example, the definition of the sets of piecewise analytic vector- and tensor fields, Definition~\ref{pwanalyticdef}, strongly relies on an underlying $\ana$-partition. Also the shift results given by Lemma~\ref{finiteshiftprop} and Lemma~\ref{shiftproptangent} rely on underlying $\ana$-partitions. Furthermore, on the half-balls $\PM_R$ and $\PP_R$ there are two types of derivatives, namely, derivatives in tangential direction and normal derivatives. Hence we also want to adapt our notion of Sobolev-Morrey seminorms such that they include information about which type of derivatives we are currently considering.

\medskip

All in all, there is quite some notation that has to be adapted to the new geometric situation, and doing these adaptions is the aim of this section. The outline is as following: In Section~\ref{adaptnotation1} we define classes of analytic vector- and tensor fields on half-balls $\PM_R$, $\PP_R$ and balls $\PW_R$, and we discuss how piecewise analytic vector- and tensor fields on an $\ana$-partition behave, when they are transformed by a local flattening.

Subsequently, in Section~\ref{fracsob} we discuss (fractional) Sobolev spaces on a disc $\PZ_R$ and in Section~\ref{adaptnotation2} we formulate versions of the shift results given by Lemma~\ref{finiteshiftprop} and Lemma~\ref{shiftproptangent} which are valid in the new geometric setting. 

In Section~\ref{interfacesubsec} we discuss the regularity of interface problems on $\PM_R$ and $\PP_R$ and in Section~\ref{adaptnotation4} and Section~\ref{adaptnotation5} we define the concepts of Sobolev-Morrey seminorms on half-balls $\PM_R$, and $\PP_R$, and on discs $\PZ_R$.

\medskip

We begin by adapting the spaces $\Hcurl$ and $\Hdiv$. For the remainder of this work we set
\begin{align*}
		\Hcurldomain{\PM_R}:=\left\{\solv\in\VSL(\PM_R)\ |\ \curl\solv\in\VSL(\PM_R)\right\} \quad{\rm and}\quad 
		\Hdivdomain{\PM_R}:=\left\{\solv\in\VSL(\PM_R)\ |\ \diverg\solv\in\SL(\PM_R)\right\}, 
\end{align*}
as well as 
\begin{align*}
		\Hcurldomain{\PW_R}:=\left\{\solv\in\VSL(\PW_R)\ |\ \curl\solv\in\VSL(\PW_R)\right\} \quad{\rm and}\quad 
		\Hdivdomain{\PW_R}:=\left\{\solv\in\VSL(\PW_R)\ |\ \diverg\solv\in\SL(\PW_R)\right\}. 
\end{align*}
In accordance with the previously used notation, we call a function $u$ on $\PM_R$ (or on $\PP_R$) smooth, if $u\in\SHM{\ell}(\PM_R)$ (or $u\in\SHM{\ell}(\PP_R)$) for all $\ell\in\N_0$. Naturally, vector- and tensor fields on $\PM_R$ or $\PP_R$ are called smooth if their respective component functions are smooth on $\PM_R$ or $\PP_R$, respectively. 
In addition, we extend the notion of piecewise smooth functions or vector fields from Definition~\ref{pwsmooth} to the new geometric setting:
\begin{definition}
		A function $u$ on $\PW_R$ is called piecewise smooth if the restrictions of $u$ to $\PM_R$ and $\PP_R$ satisfy $u\vert_{\PM_R}\in\SHM{\ell}(\PM_R)$ and $u\vert_{\PP_R}\in\SHM{\ell}(\PP_R)$ for all $\ell\in\N_0$. Consequently, vector- or tensor fields on $\PW_R$ are called piecewise smooth if their respective component functions are piecewise smooth.
\end{definition}
Again, we notice that the terminology {\it piecewise smooth} is justified by the Sobolev embedding theorem.

%When working with functions and vector fields on half-balls, we sometimes require a certain support property which is explained in the next definition.

%\begin{definition}
%		A function $u$ on a half-ball $\PM_R$ or $\PP_R$ is said to have support near the origin if there exists $0<R^*<R$ such that $\operatorname{supp}(u)$ is contained in $\PM_{R^*}$ or $\PP_{R^*}$, respectively. Consequently, a vector field $\solu$ on $\PM_R$ or $\PP_R$ is said to have support near the origin if all its component functions have support near the origin. 
%\end{definition}
%

\subsection{Analytic vector- and tensor fields on half-balls}\label{adaptnotation1}

We highlight that since $\PM_R$ is a bounded Lipschitz domain, Definition~\ref{analytictensors} straightforwardly extends to the new geometric setting. Hence we may define the sets $\anaf{\PM_R}$, $\anavec{\PM_R}$ and $\anatens{\PM_R}$ as in Definition~\ref{analytictensors} simply by replacing $\Omega$ by $\PM_R$. In the same way we also define $\anaf{\PP_R}$, $\anavec{\PP_R}$ and $\anatens{\PP_R}$, i.e., we replace $\Omega$ by $\PP_R$ in Definition~\ref{analytictensors}.

While Definition~\ref{analytictensors} immediately extends to the new geometric setting, Definition~\ref{pwanalyticdef} does not. Hence, we adapt our notation in the following way:

\begin{definition}
		For $R,M>0, \omega\geq 0$ and $k\in\Co\setminus\{0\}$, the sets $\anapw{\PW_R}$ and $\anavecpw{\PW_R}$ consist of all functions $v:\PW_R\rightarrow\Co$ and vector fields $\solv:\PW_R\rightarrow\Co^3$ such that
		\begin{alignat*}{2}
				\forall\ell\in\N_0:\ &\sum_{|\alpha|=\ell}\norm{\D^{\alpha}v}{\SL(\PM_R)}+\sum_{|\alpha|=\ell}\norm{\D^{\alpha}v}{\SL(\PP_R)}&&\leq \omega M^{\ell}(|k|+\ell)^{\ell}, \\
				\forall\ell\in\N_0:\ &\sum_{|\alpha|=\ell}\norm{\D^{\alpha}\solv}{\VSL(\PM_R)}+\sum_{|\alpha|=\ell}\norm{\D^{\alpha}\solv}{\VSL(\PP_R)}&&\leq \omega M^{\ell}(|k|+\ell)^{\ell},
		\end{alignat*}
		respectively.

		Furthermore, $\anatenspw{\PW_R}$ denotes the space of all tensor fields $\nu:\PW_R\rightarrow\Co^{3\times 3}$ for which there exist constants $\omega\geq 0,M>0$ such that 
\begin{align*}
		\forall\ell\in\N_0:\ \sum_{i,j=1}^3\sum_{|\alpha|=\ell}\norm{\D^{\alpha}(\nu_{i,j})}{\SL(\PM_R)}+\sum_{|\alpha|=\ell}\norm{\D^{\alpha}(\nu_{i,j)}}{\SL(\PP_R)}\leq \omega M^{\ell}\ell^{\ell}.
\end{align*}

\end{definition}

In a similar way we also extend the definition of the boundary spaces $\anagammaf$, $\anagammavec{\Gamma}$ and $\anagamma$ from Definition~\ref{assumptionimpedance} to the new geometric setting.
We recall that a vector field $\solg:\PZ_R\rightarrow\Co^3$ is called a tangent field on $\PZ_R$ if and only if its third component is identically zero.
\begin{definition}
		For $R,M>0, \omega\geq 0$ and $k\in\Co\setminus\{0\}$, the sets $\anagammafdisc{\PZ_R}$ and $\anagammavec{\PZ_R}$ consist of all functions $g:\PZ_R\rightarrow\Co$ and tangent fields $\solg_T:\PZ\rightarrow\Co^3$ for which there exist continuations $g^*$ and $\solg^*$ to a small three-dimensional neighbourhood $\mathcal{D}$ of the disk $\PZ_R$ such that
		\begin{align*}
				\revision{\forall\ell\in\N_0:\ }\sum_{|\alpha|=\ell}\norm{\D^{\alpha}g^*}{\SL(\mathcal{D})}\leq \omega M^{\ell}(|k|+\ell)^{\ell} \quad {\rm and}\quad
				\revision{\forall\ell\in\N_0:\ }\sum_{|\alpha|=\ell}\norm{\D^{\alpha}\solg^*}{\VSL(\mathcal{D})}\leq \omega M^{\ell}(|k|+\ell)^{\ell},
		\end{align*}
		respectively.
		Furthermore, the space $\anagammadisc{\PZ_R}$ consists of all tensor fields $\nu:\PZ_R\rightarrow\Co^{3\times 3}$ for which there exist constants $\omega\geq 0, M>0$ and a continuation $\nu^*$ of $\nu$ onto a three-dimensional neighbourhood $\mathcal{D}$ of $\PZ_R$ such that
\begin{align*}
		\revision{\forall\ell\in\N_0:\ }\sum_{i,j=1}^3\sum_{|\alpha|=\ell}\norm{\D^{\alpha}(\nu^*_{i,j})}{\VSL(\mathcal{D})}\leq \omega M^{\ell}\ell^{\ell}.
\end{align*}
\end{definition}

The subsequent lemma shows that the transformed coefficients $\wmu^{-1}$, $\weps$ and $\wzeta$ from Definition~\ref{trafoquantitiesdef} inherit many properties from $\mu^{-1}$, $\varepsilon$ and $\zeta$. 
\begin{lemma}\label{trafotensors}
		Let $\Gp=\geom$ be an $\ana$-partition, suppose that $\mu^{-1},\varepsilon\in\anatenspw{\Gp}$ satisfy \eqref{coercivedomain} and assume that $\zeta\in\anagamma$ satisfies \eqref{coercivebdr} and \eqref{zetaproperty}. 
		Furthermore, for a local flattening $\Upsilon:\PW_R\rightarrow\R^3$ of $\Gamma$ consider the quantities $\wmu^{-1}$, $\weps$ and $\wzeta$ from Definition~\ref{trafoquantitiesdef}. 

		Then, there holds $\wmu^{-1}, \weps\in\anatens{\PM_R}$ and $\wzeta\in\anagammadisc{\PZ_R}$, and there exist a constant $c>0$ as well as complex numbers $\alpha_{\wmu^{-1}}$, $\alpha_{\weps}$ and $\alpha_{\wzeta}$ with $|\alpha_{\wmu^{-1}}| = |\alpha_{\weps}| = |\alpha_{\wzeta}|=1$ such that
\begin{itemize}
		\item there holds
				\begin{align}\label{coercivehalfball}
								\realpart\SCP{\alpha_{\wmu^{-1}}\wmu^{-1}\solw}{\solw}{}+\realpart\SCP{\alpha_{\weps}\weps\solz}{\solz}{}\geq c\norm{\solw}{}^2+c\norm{\solz}{}^2
				\end{align}
				for all $\solw,\solz\in\Co^3$ uniformly in $\PM_R$,
		\item there holds
				\begin{align}\label{coercivedisc}
								\realpart\SCP{\alpha_{\wzeta}\wzeta\solz}{\solz}{}\geq c\norm{\solz}{}^2
				\end{align}
				for all $\solz\in\Co^3$ uniformly on $\PZ_R$. Furthermore, for all smooth vector fields $\solv:\PZ_R\rightarrow\Co^3$ there holds $(\wzeta\solv)_T = \wzeta\solv_T$, where $\solz_T:=\sole_3\times(\solz\times\sole_3)$ denotes the tangential component of a vector field $\solz$ on $\PZ_R$. 
\end{itemize}

		\smallskip

		Similarly, when $\Upsilon:\PW_R\rightarrow\R^3$ is a local flattening of a subdomain interface $\interf_j$, then $\wmu^{-1},\weps\in\anatenspw{\PW_R}$ and there exist a constant $c>0$ as well as complex numbers $\alpha_{\wmu^{-1}}$ and $\alpha_{\weps}$ with $|\alpha_{\wmu^{-1}}| = |\alpha_{\weps}| = 1$ such that 
\begin{align}
								\realpart\SCP{\alpha_{\wmu^{-1}}\wmu^{-1}\solw}{\solw}{}+\realpart\SCP{\alpha_{\weps}\weps\solz}{\solz}{}\geq c\norm{\solw}{}^2+c\norm{\solz}{}^2
						\end{align}\label{coercivefullball}
				for all $\solw,\solz\in\Co^3$ uniformly in $\PW_R$.

\end{lemma}

\begin{fatproof}
		The claimed analytic regularity properties of $\wmu^{-1},\weps$ and $\wzeta$ follow from the (piecewise) analyticity of $\mu^{-1}, \varepsilon, \zeta$ and the fact that $\Upsilon$ is an analytic diffeomorphism. Furthermore, the estimates \eqref{coercivehalfball}, \eqref{coercivedisc} and \eqref{coercivefullball} follow immediately from \eqref{coercivedomain} and \eqref{coercivebdr}. Finally, the equation $(\wzeta\solv)_T = \wzeta\solv_T$ can be checked by a computation. 
\end{fatproof}

Actually, the equality $(\wzeta\solv)_T = \wzeta\solv_T$ from Lemma~\ref{trafotensors} was already used in the proof of Lemma~\ref{trafolemma}. For the sake of completeness we stated it again in Lemma~\ref{trafotensors}.

Our next aim is to discuss the interplay between vector fields $\solv$ on $\Omega$ and their transforms $\hatf{\solv}$ from Definition~\ref{trafoquantitiesdef}. The proof of the following result is based on the Fa\`{a} di Bruno formula from \cite[Lem.~1.1.1]{GrandLivre}, the Leibniz rule and the estimate \eqref{elementaryestimate} from Appendix~\ref{appendix:commest}. For the sake of brevity we omit it.

%While it is clear that for a local flattening $\Upsilon:\PW_R\rightarrow\R^3$ of $\Gamma$ and $\solv\in\anavecpwCM{\Gp}{C_{\solv}}{M_{\solv}}$ there holds $\hatf{\solv}\in\anavecpwCM{\PM_R}{C'}{M'}$ for some $C',M'>0$, it is not immediately clear how the constants $C',M'$ are related to $C_{\solv}$ and $M_{\solv}$. The following lemma at least clarifies the dependence of $C'$ on $C_{\solv}$. Its proof is based on the Fa\`{a} di Bruno formula from \cite[Lemma~1.1.1]{GrandLivre} and the Leibniz rule. For the sake of brevity we omit it.

\begin{lemma}\label{interplay1}
Let $\Gp$ be an $\ana$-partition and suppose $\solv\in\anaomegavecpw{\Gp}{\omega_{\solv}}{M}$ and $\solg_T\in\anaomegavec{\Gamma}{\omega_{\solg_T}}{M}$. Furthermore, let $\Upsilon:\PW_R\rightarrow \R^3$ be a local flattening of the boundary~$\Gamma$. 
Then, the transformed quantities $\hatf{\solv}$ and $\hatf{\solg}_T$ from Definition~\ref{trafoquantitiesdef} satisfy $\hatf{\solv}\in\anaomegavecpw{\PM_R}{\rho\omega_{\solv}}{M'}$ and $\hatf{\solg}_T\in\anaomegavec{\PZ_R}{\rho\omega_{\solg_T}}{M'}$, where $\rho>0$ depends only on $\Upsilon$ and $M'>0$ depends only on $\Upsilon$ and $M$.

		\smallskip

		Similarly, when $\Upsilon:\PW_R\rightarrow\R^3$ is a local flattening of a subdomain interface $\interf_j$, the transformed quantity $\hatf{\solv}$ satisfies $\hatf{\solv}\in\anaomegavecpw{\PW_R}{\rho\omega_{\solv}}{M'}$, where again, $\rho>0$ depends only on $\Upsilon$ and $M'$ depends only on $\Upsilon$ and $M$. 
\end{lemma}

In addition to the interplay between $\solv$ and its transform $\hatf{\solv}$, we will also need to discuss the interplay between $\diverg\solv$ and $\diverg\hatf{\solv}$. The following result provides the needed information, for the sake of brevity we again omit the proof.

\begin{lemma}\label{interplay2}
		Let $\Gp$ be an $\ana$-partition and suppose that $\solv\in\Hdiv$ satisfies $\diverg\solv\in\anaomegapw{\Gp}{\omega_{\diverg\solf}}{M}$. 

		Then, with $\Upsilon:\PW_R\rightarrow\R^3$ being a local flattening of $\Gamma$ and $\hatf{\solv}$ being the transformed quantity from Definition~\ref{trafoquantitiesdef}, there holds that $\diverg\hatf{\solv}\in\anaomega{\PM_R}{\rho\omega_{\diverg\solv}}{M'}$, where $\rho>0$ depends only on $\Upsilon$ and $M'>0$ depends only on $\Upsilon$ and $M$.

		Furthermore, if $\Upsilon:\PW_R\rightarrow\R^3$ is a local flattening of a subdomain interface $\interf_j$, the transformed quantity $\hatf{\solv}$ satisfies $\diverg\hatf{\solv}\in\anaomegavecpw{\PW_R}{\rho\omega_{\diverg\solv}}{M'}$, where $\rho>0$ depends only on $\Upsilon$ and $M'>0$ depends only on $\Upsilon$ and $M$.
\end{lemma}

\subsection{Fractional Sobolev spaces on two-dimensional domains}\label{fracsob}

To properly discuss the influence of the boundary data on the solution $\solu$ we need to discuss (fractional) Sobolev spaces on $\PZ_R$. For a domain\footnote{Here, $\Omega$ denotes a two-dimensional domain and must not be confused with the three-dimensional domain $\Omega$ on which we pose Maxwell's equations.} $\Omega\subseteq\R^2$ and $\ell\in\N_0$ the spaces $\SHM{\ell}(\Omega)$ denote the usual Sobolev spaces of integer order on $\Omega$. 
Following \cite[Ch.~3]{BookMcLean} we define for $r\in (0,1)$ the spaces $\SHM{s}(\Omega)$ and $\SHMZ{s}(\Omega)$ in terms of the Aronstein-Slobodeckij seminorm
\begin{alignat*}{2}
		&\ASnorm{u}{r}{\Omega}^2:= \int_{\Omega}\int_{\Omega}\frac{|u(x)-u(y)|^2}{|x-y|^{2+2r}} \diff x\diff y, \quad &&\norm{u}{\SHM{r}(\Omega)}^2:=\norm{u}{\SL(\Omega)}^2+\ASnorm{u}{r}{\Omega}^2, \\
		&\SHMZ{r}(\Omega) := \left\{u\in\SHM{r}(\R^2)\ |\ \operatorname{supp} u\subseteq\Omega\right\}, \quad && \norm{u}{\SHMZ{r}(\Omega')}:=\norm{u}{\SHM{r}(\R^2)}.
\end{alignat*}

For $s = \ell+r$ with $\ell\in\N_0$ and $r\in (0,1)$ we set 
\begin{align*}
		&\norm{u}{\SHM{s}(\Omega)}^2:=\norm{u}{\SHM{\ell}(\Omega)}^2+\sum_{|\alpha|=\ell}\norm{\D^{\alpha}u}{\SHM{r}(\Omega)}^2, \\
		&\SHMZ{s}(\Omega):= \left\{u\in\SHM{s}(\R^2)\ |\ \operatorname{supp} u\subseteq\overline{\Omega}\right\}, \quad \norm{u}{\SHMZ{s}(\Omega)}:=\norm{u}{\SHM{s}(\R^2)}.
\end{align*}
\begin{remark}
		We highlight the fundamental difference between the spaces $\SHM{s}(\Omega)$ and $\SHMZ{s}(\Omega)$: While $\SHM{s}(\Omega)$ consists of distributions on $\Omega$, the space $\SHMZ{s}(\Omega)$ consists of distributions \textit{on the whole $\R^2$} whose support is contained in $\overline{\Omega}$.
\end{remark}

For $s\geq 0$ we define the dual spaces
\begin{align*}
		\SHM{-s}(\Omega):=\left(\SHMZ{s}(\Omega)\right)' \quad {\rm and}\quad \SHMZ{-s}(\Omega):=\left(\SHM{s}(\Omega)\right)'.
\end{align*}

For our purposes, the most important case will be $s=1/2$ and $\Omega=\PZ_R$ for some $R>0$. For $u\in\SHM{1/2}(\PZ_1)$ and $u_R(x):=u(R^{-1}x)$ there holds $u_R\in\SHM{1/2}(\PZ_R)$ and
\begin{align*}
		\norm{u_R}{\SHM{1/2}(\PZ_R)}^2 = R\norm{u}{\SL(\PZ_1)}^2+R^{-1}\ASnorm{u}{\frac{1}{2}}{\PZ_1}^2.
\end{align*}
Analogously, if $u\in\SHMZ{1/2}(\PZ_1)$ then $u_R\in\SHMZ{1/2}(\PZ_R)$ and
\begin{align*}
\norm{u_R}{\SHMZ{1/2}(\PZ_R)}^2 = R\norm{u}{\SL(\R^2)}^2+R^{-1}\ASnorm{u}{\frac{1}{2}}{\R^2}^2.
\end{align*}
In total we get that for $0<\Rmin\leq R\leq \Rmax$ there exists a constant $C>0$ depending only on $\Rmin$ and $\Rmax$ such that
\begin{align}\label{eq:scaling}
		\begin{split}
			\forall u\in\SHM{1/2}(\PZ_1):\ C^{-1}\norm{u_R}{\SHM{1/2}(\PZ_R)}\leq \norm{u}{\SHM{1/2}(\PZ_1)}\leq C\norm{u_R}{\SHM{1/2}(\PZ_R)}, \\
\forall u\in\SHMZ{1/2}(\PZ_1):\ C^{-1}\norm{u_R}{\SHMZ{1/2}(\PZ_R)}\leq \norm{u}{\SHMZ{1/2}(\PZ_1)}\leq C\norm{u_R}{\SHMZ{1/2}(\PZ_R)}.
		\end{split}
	\end{align}

	The estimates provided by \eqref{eq:scaling} allow us to apply scaling arguments to the norms $\norm{\cdot}{\SHM{1/2}(\PZ_R)}$ and $\norm{\cdot}{\SHMZ{1/2}(\PZ_R)}$, and by duality also to $\norm{\cdot}{\SHM{-1/2}(\PZ_R)}$ and $\norm{\cdot}{\SHMZ{-1/2}(\PZ_R)}$.

It remains to define vector-valued versions of $\SHM{1/2}(\Omega)$ and $\SHMZ{1/2}(\Omega)$. For a domain $\Omega\subseteq\R^2$ and $s\geq 0$ we define the spaces
\begin{align*}
		\VSHM{s}_T(\Omega):=\{(u_1,u_2,0)^T\ |\ u_1,u_2\in\SHM{s}(\Omega)\} \quad {\rm and} \quad \VSHMZ{s}_T(\Omega):=\{(u_1,u_2,0)^T\ |\ u_1,u_2\in\SHMZ{s}(\Omega)\} 
\end{align*}
which are equipped with their respective graph norms. Furthermore, for $s\geq 0$ we set
\begin{align*}
		\VSHM{-s}_T(\Omega) := \left(\VSHMZ{s}_T(\Omega)\right)'\quad {\rm and}\quad \VSHMZ{-s}_T(\Omega) := \left(\VSHM{s}_T(\Omega)\right)'.
\end{align*}
It is clear that the estimates \eqref{eq:scaling} can be extended to $\solu\in\VSHM{1/2}_T(\PZ_1)$ and $\solu\in\VSHMZ{1/2}_T(\PZ_1)$. That is, we may apply scaling arguments to the norms $\norm{\cdot}{\VSHM{1/2}_T(\PZ_R)}$ and $\norm{\cdot}{\VSHMZ{1/2}_T(\PZ_R)}$, and by duality also to the dual norms $\norm{\cdot}{\VSHM{-1/2}_T(\PZ_R)}$ and $\norm{\cdot}{\VSHMZ{-1/2}_T(\PZ_R)}$.

The following result shows that multiplication by a smooth cutoff function maps $\VSHM{1/2}_T(\PZ_R)$ to $\VSHMZ{1/2}_T(\PZ_R)$.
Moreover, it shows that if $\solw\in\VSHM{1/2}_T(\PZ_R)$ has compact support inside of $\PZ_R$ then the following lemma shows that we can relate $\norm{\solw}{\VSHM{1/2}_T(\PZ_R)}$ and $\norm{\solw}{\VSHMZ{1/2}(\PZ_R)}$. 

\begin{lemma}\label{interpolation}
		Let $0<\Rmin\leq R\leq \Rmax$ and $\solw\in\VSHM{1/2}_T(\PZ_R)$. Then, for all $\delta\in (0,R)$ there holds
		\begin{align*}
				\norm{\chi_{R,\delta} \solw}{\VSHMZ{1/2}_T(\PZ_R)}\leq C' \frac{1}{\sqrt{\delta}}\norm{\solw}{\VSHM{1/2}_T(\PZ_R)},
		\end{align*}
		where $\chi_{R,\delta}$ is a smooth cutoff function with the properties described in Definition~\ref{cutoffdef} and $C'>0$ depends only on \revisionc{$\Rmin, \Rmax$} and the constant $C>0$ from \eqref{cutoffestimate}. 
\end{lemma}

\begin{fatproof}
	%	\revisionc{The proof uses standard interpolation techniques that are common in the literature and is very similar to e.g. \cite[Theorem~3.20]{BookMcLean}; for the sake of exposition we state it.}
		We consider only the case $R=1$, the general case then follows from a scaling argument.
		According to \cite[Ch.33, Ch.36]{BookTartar} the spaces $\SHM{1/2}(\PZ_1)$ and $\SHMZ{1/2}(\PZ_1)$ are interpolation spaces, namely 
		\begin{align*}
				\SHM{1/2}(\PZ_1) = \left(\SL(\PZ_1),\SHM{1}(\PZ_1)\right)_{\frac{1}{2},2}\quad{\rm and}\quad 
				\SHMZ{1/2}(\PZ_1) = \left(\SL(\PZ_1),\SHM{1}_0(\PZ_1)\right)_{\frac{1}{2},2}
		\end{align*}
		where $\SHM{1}_0(\PZ_1)$ is the closure of smooth functions with compact support in $\PZ_1$ in $\SHM{1}(\PZ_1)$. \revisionc{The statement then follows from standard interpolation arguments, cf. \cite[Lem.~22.3]{BookTartar}.}

\end{fatproof}

In the following we will often work with function $u$ and vector fields $\solu$ that have compact support in $\PM_R\cup\PZ_R$. We highlight that the condition of $u$ and $\solu$ having compact support in $\PM_R\cup\PZ_R$ is equivalent to the existence of $\widetilde{R}\in (0,R)$ such that $\operatorname{supp}u$ and $\operatorname{supp}\solu$ are contained in the set
		$\{x\in\R^3\ |\ x_3\geq 0, |x|\leq \widetilde{R}\}$.

		\begin{lemma}\label{lem:trbdr}
				Let $0<\Rmin\leq R\leq \Rmax$ and suppose that $\solu\in\Hcurldomain{\PM_R}$ has compact support in $\PM_R\cup\PZ_R$. Then, $\solu_t:={\solu|}_{\PZ_R}\times\sole_3\in\VSHMZ{-1/2}(\PZ_R)$ and $\solu_T:=\sole_3\times({\solu|}_{\PZ_R}\times\sole_3)\in\VSHMZ{-1/2}(\PZ_R)$, and there holds the trace estimate
\begin{align*}
		\norm{\solu_t}{\VSHMZ{-1/2}_T(\PZ_R)}+\norm{\solu_T}{\VSHMZ{-1/2}_T(\PZ_R)}\leq C\left(\norm{\solu}{\VSL(\PM_R)}+\norm{\curl\solu}{\VSL(\PM_R)}\right)
\end{align*}
for a constant $C>0$ that depends only on $\Rmin$ and $\Rmax$.

Similarly, if $\solu\in\Hdivdomain{\PM_R}$ has compact support in $\PM_R\cup\PZ_R$, then ${\solu|}_{\PZ_R}\cdot\sole_3\in\SHMZ{-1/2}(\PZ_R)$ and
\begin{align*}
				\norm{\solu\cdot\sole_3}{\SHMZ{-1/2}(\PZ_R)}\leq C\left(\norm{\solu}{\VSL(\PM_R)}+\norm{\diverg\solu}{\SL(\PM_R)}\right) 
\end{align*}
for a constant $C>0$ that again depends only on $\Rmin$ and $\Rmax$. 
\end{lemma}
\begin{fatproof}
		We consider only the case $\solu\in\Hcurldomain{\PM_R}$, the case of $\solu\in\Hdivdomain{\PM_R}$ follows from similar arguments. Moreover, we restrict ourselves to the case $R=1$ since the more general case then follows from scaling arguments. By using Stein's extension operator \cite{Stein} and interpolation theory we notice that for every $\solv\in\VSHM{1/2}_T(\PZ_1)$ there exists an extension $\widetilde{\solv}\in\VSHM{1/2}(\R^2)$ be an extension of $\solv$ to $\R^2$ which satisfies
		$$\norm{\widetilde{\solv}}{\VSHM{1/2}_T(\R^2)}\leq C\norm{\solv}{\VSHM{1/2}_T(\PZ_1)}$$
		for a universal constant $C>0$. Moreover, according to \cite[Thm.~7.36]{BookRenardyRogers} we can further extend $\widetilde{\solv}$ to a vector field $\mathbf{V}\in\VSHM{1}(\R^3_{-})$, where $\R^3_{-}$ denotes the lower half-space, such that $\mathbf{V}$ satisfies
		$$\norm{\mathbf{V}}{\VSHM{1}(\R^3_{-})}\leq C\norm{\solv}{\VSHM{1/2}_T(\PZ_1)}$$
for a universal constant $C>0$.

Finally, we use partial integration \cite[Thm.~3.29]{BookMonk} and the fact that $\solu$ has compact support in $\PM_R\cup\PZ_R$ to get for all $\solv\in\VSHM{1/2}_T(\PZ_1)$
\begin{align}\label{eq:tracesbdr}
		\begin{split}
			|\langle \solu_t, \solv\rangle | & = |\langle \SCP{\curl\solu}{\mathbf{V}}{\VSL(\R^3_{-})}-\SCP{\solu}{\curl\mathbf{V}}{\VSL(\R^3_{-})}| \\
										 &\leq C\left(\norm{\solu}{\VSL(\PM_R)}+\norm{\curl\solu}{\VSL(\PM_1)}\right)\norm{\mathbf{V}}{\VSHM{1}(\R^3_{-})} \\
										 &\leq C\left(\norm{\solu}{\VSL(\PM_R)}+\norm{\curl\solu}{\VSL(\PM_1)}\right)\norm{\solv}{\VSHM{1/2}_T(\PZ_1)}
		\end{split}
	\end{align}
for a universal constant $C>1$.
Therefore, 
\begin{align*}
		\norm{\solu_t}{\VSHMZ{-1/2}_T(\PZ_1)} = \sup_{\solv\in\VSHM{1/2}_T(\PZ_1)\setminus\{0\}}\frac{|\langle \solu_t, \solv\rangle|}{\norm{\solv}{\VSHM{1/2}_T(\PZ_1)}} \leq C\left(\norm{\solu}{\VSL(\PM_1)}+\norm{\curl\solu}{\VSL(\PM_1)}\right).
\end{align*}
This shows the trace estimate for $\solu_t$. To prove the trace estimate for $\solu_T$ we note that 
$$\langle \solu_T, \solv\rangle = -\langle\solu_t, r(\solv)\rangle,$$
where $r(\solv):=\sole_3\times\solv$. Plugging this identity into \eqref{eq:tracesbdr} then proves
$$\norm{\solu_T}{\VSHMZ{-1/2}_T(\PZ_1)}\leq C\left(\norm{\solu}{\VSL(\PM_R)}+\norm{\curl\solu}{\VSL(\PM_1)}\right)\norm{r(\solv)}{\VSHM{1/2}_T(\PZ_1)},$$
and the estimate $\norm{r(\solv)}{\VSHM{1/2}_T(\PZ_1)}\leq C\norm{\solv}{\VSHM{1/2}_T(\PZ_1)}$ then concludes the proof.

\end{fatproof}

\begin{remark}
		Let us mention that Lemma~\ref{lem:trbdr} is a special case of the more general statement \cite[Prop.~3.3]{PaperFernandesGilardi}. Nevertheless, we gave an explicit proof of Lemma~\ref{lem:trbdr} because we think that it gives a good impression on how to work with the spaces $\VSHM{1/2}_T(\Omega)$ and $\VSHMZ{-1/2}_T(\Omega)$.
\end{remark}

%To see this we note that we can extend $\solu$ by zero to $\widetilde{\solu}\in\Hcurldomain{\PM_{R+1}}$ and due to $\widetilde{\solu}_T$ having compact support in $\PZ_R$ and Proposition~\ref{traceprop} we have
%		\begin{align*}
%				\norm{\widetilde{\solu}_T}{\VSHM{-1/2}(\PZ_{R+1})} \leq C\norm{\solu_T}{\VSHM{-1/2}{\PZ_R}}\leq C\left(\norm{\solu}{\VSL{\PM_R}}+\norm{\curl\solu}{\VSL{\PM_R}}\right),
%		\end{align*}
%and a scaling argument shows that the constant $C>0$ depends only on $R'$ and $R^*$. 
%

We conclude this intermediate section by discussing differential operators acting on $\Omega\subseteq\R^2$. 
With $\partial_x$ and $\partial_y$ denoting the partial derivatives in $x$ and $y$-direction, we define for every smooth function $v$ and vector field $\solv = (v_1, v_2, 0)^T$ on $\PZ_R$ the operators
		\begin{align*}
		\nablat v := (\partial_x v, \partial_y v, 0)^T, \quad  \divergt\solv:=\partial_x v_1+\partial_y v_2, \quad
		\rcurlt v := (\partial_y v, -\partial_x v, 0)^T ,\quad  \curlt\solv := \partial_x v_2-\partial_y v_1.
		\end{align*}

		\begin{lemma}\label{lem:divtmappings}
		For $s\in\R$, the operators $\nablat$, $\divergt$, $\rcurlt$ and $\curlt$ have the mapping properties
\begin{alignat}{2}
	&\nablat ,\rcurlt:\SHM{s}(\PZ_R)\rightarrow\VSHM{s-1}_T(\PZ_R), \quad &&\divergt, \curlt: \VSHM{s}_T(\PZ_R)\rightarrow\SHM{s-1}(\PZ_R),\label{divtmapping1} \\
	&\nablat ,\rcurlt:\SHMZ{s}(\PZ_R)\rightarrow\VSHMZ{s-1}_T(\PZ_R), \quad &&\divergt, \curlt: \VSHMZ{s}_T(\PZ_R)\rightarrow\SHMZ{s-1}(\PZ_R).\label{divtmapping2}
\end{alignat}
		\end{lemma}
		\begin{fatproof}
		We consider only the operator $\nablat$, the mapping properties of $\rcurlt$, $\divergt$ and $\curlt$ are proved using similar arguments. The claimed mapping properties of $\nablat$ are clear for $s\geq 1$ and $s\leq 0$. It remains to consider $0<s<1$. 
		For $s\in (0,1)$, the mapping properties \eqref{divtmapping1} follow from \cite[Rem.~1.4.4.7]{BookGrisvard}. For \eqref{divtmapping2} we proceed using the same arguments as in \cite[Proof of Thm.~1.4.4.6]{BookGrisvard}. For further details we also refer to \cite[Lem.~7.34]{MyDiss}.
		\end{fatproof}

For any smooth vector field $\solv$ on $\PM_R$ (or $\PP_R$), a direct calculation shows the analogue to \eqref{surfacecurlconnect}, namely 
\begin{align}\label{surfacetraceprop}
		\curlt\solv_T = \divergt\solv_t = \sole\cdot(\curl\solv)\vert_{\PZ_R},
\end{align}
where $\sole:=(0,0,1)^T$ if $\solv$ is defined on $\PM_R$ and $\sole:=(0,0,-1)^T$ if $\solv$ is defined on $\PP_R$.
%where $\sole_3:=(0,0,1)^T$, and $\solv_T:=\sole_3\times(\solv\times\sole_3)$ and $\solv_t:=\solv\times\sole_3$ denote the tangential component and the tangential trace of $\solv$ on $\PZ_R$, respectively.

\subsection{Regularity shift results of vector fields on half-balls and discs}\label{adaptnotation2}

We proceed by adapting Proposition~\ref{finiteshiftprop} and Proposition~\ref{shiftproptangent} to the geometric settings of half-balls and discs. These results will be crucial in the discussion of analytic regularity of a solution $\solu$ of Maxwell's equations \eqref{Maxwellorig} near the boundary $\Gamma$. 

The equivalent of Proposition~\ref{finiteshiftprop} reads as follows:

\begin{proposition}\label{finiteshifttrafo}
		Let $0<\Rmin\leq R\leq\Rmax$ and suppose that $\nu$ is a smooth tensor field on $\PM_R$ that satisfies the coercivity condition \eqref{coercivenudomain} for $\Omega=\PM_R$. In addition, let $\solv\in\Hcurldomain{\PM_R}$ satisfy $\operatorname{supp}(\solv)\subseteq\PM_R\cup\PZ_R$ as well as $\nu\solv\in\Hdivdomain{\PM_R}$. 
		Under these assumptions, suppose that either 
		\begin{enumerate}
				\item the tangential component $\solv_T$ of $\solv$ satisfies $\solv_T= \solg_T$ for a tangent field $\solg_T\in\VSHMZ{1/2}_T(\PZ_R)$ on $\PZ_R$,
				\item or the tangential trace $\solv_t$ of $\solv$ satisfies $\solv_t = \solg_T$ for a tangent field $\solg_T\in\VSHMZ{1/2}_T(\PZ_R)$ on $\PZ_R$,
				\item or the normal trace of $\nu\solv$ satisfy $\nu\solv\cdot\sole_3 = h$ for a function $h\in\SHMZ{1/2}(\PZ_R)$ on $\PZ_R$, where $\sole_3 := (0,0,1)^T$.
		\end{enumerate}
		Then, $\solv\in\VSHM{1}(\PM_R)$ and 
		\begin{align*}
				\norm{\solv}{\VSHM{1}(\PM_R)}\leq C\left(\norm{\curl\solv}{\VSL(\PM_R)}+\norm{\diverg\nu\solv}{\SL(\PM_R)}+\norm{\solg_T}{\VSHMZ{1/2}_T(\PZ_R)}\right)
		\end{align*}
		\revision{in the cases 1. and 2., and}
		\begin{align*}
				\norm{\solv}{\VSHM{1}(\PM_R)}\leq C\left(\norm{\curl\solv}{\VSL(\PM_R)}+\norm{\diverg\nu\solv}{\SL(\PM_R)}+\norm{h}{\SHMZ{1/2}(\PZ_R)}\right)
		\end{align*}
		\revision{in case 3. In all three cases, the constant} $C>0$ depends only on $\Rmin,\Rmax$ and $\nu$.
\end{proposition}
\begin{fatproof}
		We only consider the case of $\solv_T=\solg_T$ for a tangent field $\solg_T\in\VSHMZ{1/2}_T(\PZ_R)$, the other cases follow analogously.
		First, let us assume that $R=1$. We choose a smooth domain $\Ud\subseteq\R^3$ such that $\PM_{3/2}\subseteq\Ud\subseteq\PM_{2}$. We extend $\solv$ by zero to a vector field $\wt{\solv}$ on $\Ud$, and after extending $\nu$ to a smooth tensor field on $\Ud$ we have $\wt{\solv}\in\Hcurldomain{\Ud}$ together with $\nu\wt{\solv}\in\Hdivdomain{\Ud}$ and $\wt{\solv}_T = \wt{\solg}_T$ on $\partial\Ud$, where $\wt{\solg}_T$ denotes the extension by zero of $\solg_T$ to $\partial\Ud$. 

		Due to $\operatorname{supp}\solg_T\subseteq\PZ_1$ we have $\wt{\solg}_T\in\VSHM{1/2}_T(\partial\Ud)$ and
		$$ \norm{\wt{\solg}_T}{\VSHM{1/2}_T(\partial\Ud)}\leq C\norm{\solg_T}{\VSHMZ{1/2}_T(\PZ_1)}$$
		for a constant $C>0$, for more details see e.g. \cite[Prop.~3.3]{PaperFernandesGilardi}.

				By applying Proposition~\ref{finiteshiftprop} we obtain $\wt{\solv}\in\VSHM{1}(\Ud)$ together with
				\begin{align*}
						\norm{\wt{\solv}}{\VSHM{1}(\Ud)}\leq C\left(\norm{\curl\wt{\solv}}{\VSL(\Ud)}+\norm{\diverg\nu\wt{\solv}}{\SL(\Ud)}+\norm{\wt{\solg}_T}{\VSHM{1/2}_T(\partial\Ud)}\right)
				\end{align*}
				which implies $\solv\in\VSHM{1}(\PM_1)$ and
				\begin{align*}
						\norm{\solv}{\VSHM{1}(\PM_1)}\leq C\left(\norm{\curl\solv}{\VSL(\PM_1)}+\norm{\diverg\nu\solv}{\SL(\PM_1)}+\norm{\solg_T}{\VSHMZ{1/2}_T(\PZ_1)}\right).
				\end{align*}
				This proves the claim for $R=1$. For general $0<\Rmin\leq R\leq\Rmax$ the claim follows from \eqref{eq:scaling} and a scaling argument.
\end{fatproof}

Another important tool for the discussion of the influence of boundary conditions on the transformed solution $\wu$ are the surface differential operators $\curlt$ and $\divergt$ from Section~\ref{fracsob}.
These operators allow for an adaption of Proposition~\ref{shiftproptangent} to the new geometric setting of half-balls and discs, which reads as follows: 

\begin{proposition}\label{shiftprop}
		Let $0<\Rmin\leq R\leq \Rmax$, let $\solv\in\VSL_T(\PZ_R)$ have compact support in $\PZ_R$ and let $\lambda:\PZ_R\rightarrow\Co^{3\times 3}$ be smooth and satisfy \eqref{coercivenubdr} for $\Gamma = \PZ_R$, as well as $\lambda\solw_T = (\lambda\solw)_T$ for all smooth vector fields $\solw$ on $\PZ_R$. Furthermore, suppose that $\divergt\solv\in\SHMZ{-1/2}(\PZ_R)$ and $\curlt(\lambda\solv)\in\SHMZ{-1/2}(\PZ_R)$. Then, $\solv\in\VSHMZ{1/2}_T(\PZ_R)$ with 
		\begin{align*}
				\norm{\solv}{\VSHMZ{1/2}_T(\PZ_R)}\leq C\left(\norm{\divergt\solv}{\SHMZ{-1/2}(\PZ_R)}+\norm{\curlt(\lambda\solv)}{\SHMZ{-1/2}(\PZ_R)}\right),
		\end{align*}
		where the constant $C>0$ depends only on $\lambda, \Rmin$ and $\Rmax$.

		Similarly, when a tangent field $\solv\in\VSL_T(\PZ_R)$ has compact support in $\PZ_R$ and satisfies $\divergt(\lambda\solv)\in\SHMZ{-1/2}(\PZ_R)$ as well as $\curlt\solv\in\SHMZ{-1/2}(\PZ_R)$, then $\solv\in\VSHMZ{1/2}_T(\PZ_R)$ with
		\begin{align*}
				\norm{\solv}{\VSHMZ{1/2}_T(\PZ_R)}\leq C\left(\norm{\divergt(\lambda\solv)}{\SHMZ{-1/2}(\PZ_R)}+\norm{\curlt\solv}{\SHMZ{-1/2}(\PZ_R)}\right),
		\end{align*}
		where $C>0$ again depends only on $\lambda$, $\Rmin$ and $\Rmax$.
\end{proposition}

For the proof we need the following result:
\begin{proposition}\label{CostabelMcIntosh}
		Let $s,t\in\R$ and let $u\in\SHMZ{s}(\PZ_1)$ satisfy $u = \curlt\solv$ for some $\solv\in\VSHMZ{t}(\PZ_1)$. Then, there exists $\solw\in\VSHMZ{s+1}(\PZ_1)$ such that $u=\curlt\solw$ and
		\begin{align*}
				\norm{\solw}{\VSHMZ{s+1}(\PZ_1)}\leq C\norm{u}{\VSHMZ{s}(\PZ_1)}
		\end{align*}
		for a constant $C>0$ that depends only on $s$. 
\end{proposition}
\begin{fatproof}
		According to \cite[Cor.~3.4, (3.13)]{CurlInverse} there exist operators $T_1$ and $T_2$ with the mapping properties $T_1:\VSHMZ{r}_T(\PZ_1)\rightarrow\SHMZ{r+1}(\PZ_1)$ and $T_2:\SHMZ{r}(\PZ_1)\rightarrow\VSHMZ{r+1}_T(\PZ_1)$ for all $r\in\R$ such that
		$$\solv = \nablat T_1\solv+T_2\curlt\solv.$$
		Hence, $\solw:=T_2\curlt\solv = T_2u$ satisfies $u = \curlt\solw$ together with the desired estimate.
\end{fatproof}

\begin{fatproofmod}{Proposition~\ref{shiftprop}}
		We consider only the case $R=1$, the general case $0<\Rmin\leq R\leq \Rmax$ follows from a scaling argument. By implicitly extending $\solv$ by zero we may consider $\solv$ as an element of $\VSL_T(\R^2)$. Due to $\solv$ having compact support in the ball $\PZ_1$, it satisfies\footnote{We implicitly extend $\lambda$ to a smooth tensor field on $\R^2$.} $\curlt(\lambda\solv)\in\SHM{-1/2}(\R^2)$. Hence, by Proposition~\ref{CostabelMcIntosh} there exists $\solw\in\VSHMZ{1/2}_T(\PZ_1)$ such that $\curlt\solw =\curlt(\lambda\solv)$. By \cite[Prop.~4.1]{CurlInverse} this implies $\lambda\solv-\solw = \nabla_{\Gamma}\xi$ for some $\xi\in\SHM{1}(\R^2)$ satisfying $\operatorname{supp}\xi\subseteq\PZ_1$. Therefore, $\divergt(\lambda^{-1}\nablat\xi) =\divergt\solv-\diverg(\lambda^{-1}\solw)\in\SHMZ{-1/2}(\PZ_1)$, and we stress that the space $\SHMZ{-1/2}(\PZ_1)$ consists of distributions on the whole $\R^2$, that is, the equation $\divergt(\lambda^{-1}\nablat\xi) =\divergt\solv-\diverg(\lambda^{-1}\solw)$ holds on the whole $\R^2$. Hence, elliptic regularity theory proves 
		\begin{align*}
				\norm{\xi}{\SHMZ{3/2}(\PZ_1)}&\leq C\left(\norm{\divergt\solv}{\SHMZ{-1/2}(\PZ_1)}+\norm{\solw}{\SHMZ{1/2}(\PZ_1)}\right) \\
				&\leq C\left(\norm{\divergt\solv}{\SHMZ{-1/2}(\PZ_1)}+\norm{\curlt\solv}{\VSHMZ{-1/2}_T(\PZ_1)}\right)
		\end{align*}
		for some constant $C>0$ depending only on $\lambda$. 
		
%		Follows from Proposition~\ref{shiftproptangent} and similar arguments as in the proof of Proposition~\ref{finiteshifttrafo}.
\end{fatproofmod}

\subsection{Regularity of interface problems on half-balls and discs}\label{interfacesubsec}	

When transforming a piecewise smooth solution $\solu$ of Maxwell's equations \eqref{Maxwellorig} by a local flattening $\Upsilon$ of a subdomain interface $\interf_j$, we obtain a vector field $\wu$ which is piecewise smooth on $\PW_R$ and satisfies certain transmission conditions across the disc $\PZ_R$. Subsequently, we provide regularity shift results for such transformed vector fields. As a first step we follow the notation in \cite{MaxwellMyself} and define a normal jump operator: For any vector field $\solv$ that is piecewise smooth on $\PW_R$ we set
\begin{align*}
		\njump{\solv}:=\solv_{1}\vert_{\PZ_R}\cdot\sole_3-\solv_{2}\vert_{\PZ_R}\cdot\sole_3%\quad{\rm and}\quad \tjump{\solv}:=\sole_3\times\left[\left(\solv_1\vert_{\PZ_R}-\solv_2\vert_{\PZ_R}\right)\times\sole_3\right],
\end{align*}
where $\solv_1:=\solv\vert_{\PM_R}$ and $\solv_2:=\solv\vert_{\PP_R}$, and $\sole_3:=(0,0,1)^T$.

{\revisioncolorred
The proof of the following result is based on \cite[Lem.~4.2, Lem.~4.5]{MaxwellMyself}.

\begin{proposition}\label{helmholtzmod}
		Let $0<\Rmin\leq R\leq \Rmax$ and let $\nu:\PW_{R+1}\rightarrow\Co^{3\times 3}$ be piecewise smooth and satisfy \eqref{coercivenudomain} for $\Omega=\PW_{R}$. Furthermore, let $\solv\in\Hcurldomain{\PW_{R}}$ be piecewise smooth on $\PW_{R+1}$ with compact support in $\PW_R$ as well as $\njump{\nu\solv}=h\in\SHMZ{1/2}(\PZ_{R})$. Then, there holds 
		\begin{align*}
				\norm{\solv}{\VSHM{1}(\PP_R)}+\norm{\solv}{\VSHM{1}(\PM_R)}\leq C\left(\norm{\curl\solv}{\VSL(\PW_R)}+\norm{\diverg\nu\solv}{\SL(\PP_R)}+\norm{\diverg\nu\solv}{\SL(\PM_R)}+\norm{h}{\SHMZ{1/2}(\PZ_{R})}\right),	
		\end{align*}
		where the constant $C>0$ depends only on $\nu, \Rmin$ and $\Rmax$.
\end{proposition}

\begin{fatproof}
		We consider only the case $R=1$, the general case follows from a scaling argument. Due to $\solv$ having compact support in $\PW_1$ we may extend $\solv$ by zero to $\PW_{2}$. Let $\widetilde{\nu}$ be a smooth continuation of $\nu$ to $\PW_2$ which satisfies \eqref{coercivenudomain} on $\PW_2$.

		We choose a smooth domain $\Ud\subseteq\R^3$ such that $\PM_{3/2}\subseteq\Ud\subseteq\PM_2$. Note that this implies $\PZ_1\subseteq\partial\Ud$. Furthermore, we set $\Gp_1:=\Ud$ and $\Gp_2:=\PW_{2}\setminus\Ud$ and notice that $\Gp = \{\PW_2,\Gp_1,\Gp_2\}$ is a $\CM{\infty}$-partition in the sense of Definition~\ref{partitiondefsmooth}. Hence, according to \cite[Lem.~4.5]{MaxwellMyself} there exists a decomposition $\solv=\solw+\nabla\varphi$ with
		\begin{align}\label{westimate}
				\norm{\solw}{\VSHM{1}(\Gp_1)}+\norm{\solw}{\VSHM{1}(\Gp_2)}\leq C\norm{\curl\solv}{\VSL(\PW_{1})} \quad {\rm and}\quad \norm{\varphi}{\SHM{1}(\PW_{2})}\leq C\norm{\solv}{\VSL(\PW_1)},
		\end{align}
		where we exploited that $\norm{\curl\solv}{\VSL(\PW_{2})} = \norm{\curl\solv}{\VSL(\PW_{1})}$ and $\norm{\solv}{\VSL(\PW_{2})} = \norm{\solv}{\VSL(\PW_1)}$. Furthermore, $\solw$ satisfies $\SCP{\wt{\nu}\solw}{\nabla\xi}{\VSL(\PW_2)}=0$ for all $\xi\in\SHM{1}(\PW_2)$, that is, $\diverg\wt{\nu}\solw=0$ in $\PW_2$ and $\wt{\nu}\solw\cdot\soln=0$ on $\partial\PW_2$. 

For this reason we have 
		\begin{alignat*}{2}
				-\diverg\wt{\nu}\nabla\varphi &= \diverg\wt{\nu}\solv\quad &&{\rm on}\ \PP_{2}\cup\PM_{2},\\
				\njump{\wt{\nu}\nabla\varphi} &= \widetilde{h}  \quad && {\rm on}\ \partial\Ud, \\
				\wt{\nu}\nabla\varphi\cdot\soln &=0 \quad && {\rm on}\ \partial\PW_{2},
		\end{alignat*}
		where $\widetilde{h}$ is the continuation by zero of $h$ to $\partial\Ud$. Without loss of generality we may assume that $\SCP{\varphi}{1}{\SL(\PW_2)}=0$, hence \cite[Lem.~4.2]{MaxwellMyself} implies
		\begin{align}\label{phiestimate}
				\begin{split}
						\norm{\varphi}{\SHM{2}(\Gp_1)}+\norm{\varphi}{\SHM{2}(\Gp_2)}&\leq C\left(\norm{\diverg\wt{\nu}\solv}{\SL(\Gp_1)}+\norm{\diverg\wt{\nu}\solv}{\SL(\Gp_2)}+\norm{\widetilde{h}}{\SHM{1/2}(\partial\Ud)} \right) \\
																					 &\leq C\left(\norm{\diverg\wt{\nu}\solv}{\SL(\PM_1)}+\norm{\diverg\wt{\nu}\solv}{\SL(\PP_1)}+\norm{h}{\SHMZ{1/2}(\PZ_1)} \right),
				\end{split}
					\end{align}
		where we exploited that $\operatorname{supp}\solv\cap\Gp_1\subseteq\PM_1$, $\operatorname{supp}\solv\cap\Gp_2\subseteq\PP_1$ and 
\begin{align*}
\norm{\widetilde{h}}{\SHM{1/2}(\partial\Ud)} \leq C \norm{h}{\SHMZ{1/2}(\PZ_1)},
\end{align*}
see \cite[Prop.~3.3]{PaperFernandesGilardi} for a thorough proof of this estimate.
Finally, due to $\PM_1\subseteq\Gp_1$ and $\PP_1\subseteq\Gp_2$ the estimates \eqref{westimate}-\eqref{phiestimate} yield
\begin{align*}
				\norm{\solv}{\VSHM{1}(\PP_1)}+\norm{\solv}{\VSHM{1}(\PM_1)}\leq C\left(\norm{\curl\solv}{\VSL(\PW_1)}+\norm{\diverg\nu\solv}{\SL(\PP_1)}+\norm{\diverg\nu\solv}{\SL(\PM_1)}+\norm{h}{\SHMZ{1/2}(\PZ_{1})}\right),	
\end{align*}
thus the proof is complete.
\end{fatproof}
}

\subsection{Anisotropic Sobolev-Morrey seminorms on half-balls}\label{adaptnotation4}

On the half-balls $\PM_R$ and $\PP_R$ it is natural to distinguish between derivatives in tangential direction to the disc $\PZ_R$ and derivatives normal to $\PZ_R$. Indeed, when considering a solution $\wu$ of the transformed Maxwell equations from Lemma~\ref{trafolemma}, tangential derivatives of $\wu$ can be estimated by similar techniques as in the proof of Theorem~\ref{Mainresult1}. 
In contrast to this, the normal derivatives can only be estimated by exploiting properties of the underlying Maxwell equations. 
Hence, in order to account for the different nature of tangential and normal derivatives we follow \cite{GrandLivre, MaxwellTomezyk} and work with subsequently defined anisotropic Sobolev-Morrey seminorms on half-balls. 

\medskip 

We recall that for any smooth vector field $\solv = (v_1, v_2, v_3)^T$ and any multiindex $\beta\in\N_0^3$ we use the notation $\D^{\beta}(\solv):=(\D^{\beta}(v_1), \D^{\beta}(v_2), \D^{\beta}(v_3))^T$. For any smooth vector field $\solv$ defined on a half-ball $\PP_R$ or $\PM_R$ and a multiindex $\beta\in\N_0^3$ we call $\D^{\beta}(\solv)$ a tangential derivative, if the third component of $\beta$ is zero, i.e., if $\beta\in\N_0^2\times\{0\}$. Conversely, if the first two components of $\beta$ are zero we call $\D^{\beta}(\solv)$ a normal derivative of $\solv$.

\begin{remark}\label{indexdistinction}
Henceforth, we mark tangential derivatives by an apostrophe and normal derivatives by an asteriks. That is, a multiindex $\beta'\in\N_0^3$ always denotes tangential derivatives, and $\beta^*\in\N_0^3$ indicates a normal derivative. 
\end{remark}

Similarly as before, for any $R>0$  we abbreviate $\norm{\solv}{\VSL(\PP_R)}$ by $\norm{\solv}{\PP_R}$ and $\norm{\solv}{\VSL(\PM_R)}$ by $\norm{\solv}{\PM_R}$. \revision{Moreover, we set $\norm{\solv}{\PS_R}:=\norm{\solv}{\PP_R}+\norm{\solv}{\PM_R}$.} With this convention, we we follow \cite{GrandLivre,MaxwellTomezyk} and define wavenumber-dependent anisotropic Sobolev-Morrey seminorms on half-balls $\PM_R$ and $\PP_R$ as follows:

\begin{definition}\label{anismorrey}
		Let $R>0$, $\ell,n\in\N_0$, let a wavenumber $k\in\Co\setminus\{0\}$ be given, and let $\solv$ be a smooth vector field on $\PM_R$. Then, with $R(\ell,n,k):=\frac{R}{2(\ell+n+|k|)}$ we define
	 \begin{align*}
	 \morrey{\solv}_{\ell,n,\PM_R}:=\max_{0<\rho\leq R(\ell,n,k)}\max_{\substack{|\delta'|=\ell \\ |\alpha^*| = n}}\rho^{\ell+n}\norm{\D^{\delta'+\alpha^*}(\solv)}{\PM_{R-(\ell+n+|k|)\rho}},
	\end{align*}
	as well as 
	 \begin{align*}
	\morreyone{\solv}_{\ell,n,\PM_R}:=\max_{0<\rho\leq R(\ell,n,k)}\max_{\substack{|\delta'|=\ell \\ |\alpha^*| = n}}\rho^{\ell+n+1}\norm{\D^{\delta'+\alpha^*}(\solv)}{\PM_{R-(\ell+n+|k|)\rho}}.
	\end{align*}
	
	In addition, for a smooth tensor field $\nu$ on $\PM_R$ we define the cumulative quantity
	\begin{align*}
		\summorreynu{\solv}{\ell,n,\PM_R} := \morrey{\solv}_{\ell,n,\PM_R}+|k|^{-1}\morrey{\nu\curl\solv}_{\ell,n,\PM_R}.
	\end{align*}
	If $\solv$ is a smooth vector field on $\PP_R$, we define the quantities $\morrey{\solv}_{\ell,n,\PP_R},\ \morreyone{\solv}_{\ell,n,\PP_R}$ and $\summorreynu{\solv}{\ell,n,\PP_R}$ analogously by replacing $\PM_R$ and $\PM_{R-(\ell+n+|k|)\rho}$ by $\PP_R$ and $\PP_{R-(\ell+n+|k|)\rho}$ in all instances.
	
	Finally, we define 
	\begin{align*}
		\morrey{\solv}_{\ell,n,\PS_R}&:=\morrey{\solv}_{\ell,n,\PM_R}+\morrey{\solv}_{\ell,n,\PP_R},\\
		\morreyone{\solv}_{\ell,n,\PS_R}&:=\morreyone{\solv}_{\ell,n,\PM_R}+\morreyone{\solv}_{\ell,n,\PP_R},
	\end{align*}
    and 
    \begin{align*}
		\summorreynu{\solv}{\ell,n,\PS_R}:=\summorreynu{\solv}{\ell,n,\PM_R}+ \summorreynu{\solv}{\ell,n,\PP_R}
	\end{align*}
	for all vector fields $\solv$ and tensor fields $\nu$ that are piecewise smooth on $\PW_R$.
\end{definition}

\medskip

We recall that for a smooth vector field $\solv$, a smooth tensor field $\nu$ and a multiindex $\beta$, the commutator field $\commute(\nu,\solv,\beta)$ is defined as 
\begin{align*}
	\commute(\nu,\solv,\beta) := \D^{\beta}(\nu\solv)-\nu\D^{\beta}(\solv).
\end{align*}

The following lemma generalizes the commutator estimate provided by Lemma \ref{commutatorestimate} to anisotropic Sobolev-Morrey seminorms. The proof is postponed to Appendix \ref{appendix:commest}.

\begin{lemma}\label{commutatoranis}
		Let $0<R\leq \Rmax$ and assume that $\nu$ is a tensor field satisfying $\nu\in\anav(\PM_R)$. Moreover, let $k\in\Co\setminus\{0\}$ and  multiindices $\alpha,\beta\in\N_0^3$ be given and set $m:=|\alpha|+|\beta|$.

		We decompose $\alpha$ and $\beta$ into its normal and tangential components, i.e., $\alpha = \alpha'+\alpha^*$ and $\beta = \beta'+\beta^*$, and we set $m':=|\alpha'|+|\beta'|$, as well as $m^*:=|\alpha^*|+|\beta^*|$.  Then, for all $0<\rho\leq \frac{R}{2(m+|k|)}$ there holds 
	\begin{align}\label{firstineq}
			\begin{split}
					\rho^{m}\norm{\D^{\alpha}(\commute(\nu,\solv,\beta))}{\PM_{R-(m-1+|k|)\rho}} \leq C\sum_{d=0}^{m'}&\sum_{b=0}^{m^*-1}K^{m-1-d-b}\morrey{\solv}_{d,b,\PM_R}\\
																													  &+C\sum_{d=0}^{m'-1}K^{m'-1-d}\morrey{\solv}_{d,m^*,\PM_R}
\end{split}
	\end{align}
    for all smooth vector fields $\solv$, where the constants $C>0$ and $K\geq 1$ depend only on $\nu$ and $\Rmax$. Moreover, for all $\ell,n\in\N_0$ there holds
    \begin{align}\label{secondineq}
    	\morrey{\nu\solv}_{\ell,n,\PM_R}\leq C\sum_{d=0}^{\ell}\sum_{b=0}^nK^{\ell+n-d-b}\morrey{\solv}_{d,b,\PM_R}
    \end{align}
	for all smooth vector fields $\solv$, where the constants $C>0$ and $K\geq 1$ again depend only on $\nu$ and $\Rmax$. Both inequalities remain valid when $\PM_R$ is replaced by $\PP_R$.
\end{lemma}

There is a double sum occurring in the statement of Lemma~\ref{commutatoranis}, and since in subsequent proofs we will need to apply Lemma~\ref{commutatoranis} iteratively, there will be occurrences of quadruple sums. The following auxiliary result helps us to ease the presentation in the presence of such sums. 
Its proof follows from elementary manipulations and standard techniques; for the sake of brevity we omit it.

\begin{lemma}\label{vierersumme}
	Let $(\alpha_{i,j})_{(i,j)\in\N_0\times\N_0}$ and $(\beta_{i,j})_{(i,j)\in\N_0\times\N_0}$ be families of non-negative real numbers and assume that there exist constants $C>0$ and $K\geq 1$ such that for all $\ell, n\in\N_0$ there holds 
	\begin{align*}
		 \alpha_{\ell,n+1}\leq C\sum_{i=0}^{\ell+1}\sum_{j=0}^nK^{\ell+1+n-i-j}\beta_{i,j}.
	\end{align*}
	Then, for any constant $L\geq 1$ there holds
	\begin{align*}
		\sum_{i=0}^{\ell}\sum_{j=1}^{n+1}L^{\ell+n+1-i-j}\alpha_{i,j}\leq C\sum_{p=0}^{\ell+1}\sum_{q=0}^nA^{\ell+1+n-p-q}\beta_{p,q}
	\end{align*}
	for all $\ell,n\in\N_0$, where $A:=2\max\{K,L\}$.
\end{lemma}
%\begin{fatproof}
%	For all $\ell,n\in\N_0$ and $P:=\operatorname{max}\{K,L\}$ there holds 
%	\begin{align*}
%		\sum_{i=0}^{\ell}\sum_{j=1}^{n+1}L^{\ell+n+1-i-j}\alpha_{i,j}&\leq C\sum_{i=0}^{\ell}\sum_{j=1}^{n+1}L^{\ell+n+1-i-j}\sum_{d=0}^{i+1}\sum_{m=0}^{j-1}K^{i+j-d-m}\beta_{d,m} \\
%		&\leq C\sum_{i=0}^{\ell}\sum_{j=1}^{n+1}\sum_{d=0}^{i+1}\sum_{m=0}^{j-1}P^{\ell+n+1-d-m}\beta_{d,m} \\
%		&= C\sum_{p=0}^{\ell+1}\sum_{q=0}^nC_{p,q}\beta_{p,q}.
%	\end{align*}
%	Let $p\in\{0,...,\ell+1\}$ and $q\in\{0,...,n\}$. Then the coefficient $C_{p,q}$ satisfies
%	\begin{align*}
%		C_{p,q} = \sum_{i=\operatorname{max}\{0,p-1\}}^{\ell}\sum_{j=q+1}^{n+1}P^{\ell+n+1-p-q}&\leq (\ell-p+2)(n-q+1)P^{\ell+n+1-p-q} \\
%		&\leq A^{\ell+n+1-p-q},
%	\end{align*}
%	where the last estimate is due to $(n+1)\leq 2^n$ for all $n\in\N_0$.
%\end{fatproof}

\subsection{Sobolev-Morrey seminorms on discs}\label{adaptnotation5}

In order to capture the influence of the boundary conditions \eqref{transformedbdry} on the transformed solution $\wu$, it will be necessary to estimate tangential derivatives on discs $\PZ_R$. The subsequent definition of Sobolev-Morrey seminorms on $\PZ_R$ provides a useful tool for that.

To shorten notation we define $\norm{\solv}{\frac{1}{2},\PZ_R}:=\norm{\solv}{\VSHM{1/2}_T(\PZ_R)}$ and recall that multiindices with an apostrophe always denote tangential derivatives, i.e., a multiindex $\delta'\in\N^3_0$ is required to be zero in its last component, see Remark~\ref{indexdistinction}. 

\begin{definition}\label{anismorreybdr}
		Let $R>0$, $\ell\in\N_0$ and $k\in\Co\setminus\{0\}$ be given. Then, with $R(\ell,k):=\frac{R}{2(\ell+|k|)}$, we define
	\begin{align*}
	\morreyhalf{\solv}_{\ell,\PZ_R}:=\max_{0<\rho\leq R(\ell,k)}\max_{|\delta'|=\ell}\rho^{\ell+1/2}\norm{\D^{\delta'}(\solv)}{\frac{1}{2},\PZ_{R-(\ell+|k|)\rho}}
	\end{align*}
	for all functions and vector fields $\solv$ that are smooth on $\PZ_R$.
\end{definition}

The following result is the analog to the commutator estimate in Lemma \ref{commutatorestimate} for Sobolev-Morrey seminorms on $\PZ_R$. The proof is analogous to the one of Lemma \ref{commutatorestimate} stated in Appendix \ref{appendix:commest}.  

\begin{lemma}\label{commutatorestimatetangential}
		Let $0<R\leq \Rmax$ and assume that $\nu:\PZ_R\rightarrow\Co^{3\times3}$ is a tensor field which satisfies $\nu\in\anagammadisc{\PZ_R}$, that is, there exist $\omega\geq 0, M>0$ such that
		\begin{align}\label{discanalytic}
				\forall\ell\in\N_0:\ \sum_{i,j=1}^3\sum_{|\alpha|=\ell}\norm{\D^{\alpha}\nu_{i,j}}{\SL(\PZ_R)}\leq \omega M^{\ell}\ell^{\ell}.
	\end{align}
	Moreover, let $\alpha',\beta'\in\N_0^2\times\{0\}$ be arbitrary, set $b:=|\alpha'|+|\beta'|$ and assume that $k\in\Co\setminus\{0\}$ is given. Then, for all $0<\rho\leq \frac{R}{2(b+|k|)}$ there holds 
	\begin{align*}
			\rho^{b+1/2}\norm{\D^{\alpha'}(\commute(\nu,\solv,\beta'))}{\PZ_{R-(b-1+|k|)\rho}} \leq  C\sum_{d=0}^{b-1}K^{b-1-d}\morreyhalf{\solv}_{d,R},
	\end{align*}
	for all smooth scalar functions or vector fields $\solv$, where the constants $C,K>0$ depend only on $M$ and $\Rmax$. 
\end{lemma}

%% file: analytic/analytic_normal.tex
\section{Analytic regularity near boundaries and subdomain interfaces}\label{auxiliarsec}

Our next goal is to prove analytic regularity results for the transformed Maxwell equations \eqref{trafosystem}. This will provide the necessary tools for the subsequent proof of Theorem~\ref{Mainresult2}. Henceforth, we suppose that $\solu$ is a piecewise smooth solution of Maxwell's equations \eqref{Maxwellorig} and consider the transformed solution $\wu$ from Definition~\ref{trafoquantitiesdef}; our aim is to discuss the piecewise analyticity properties of $\wu$ on the ball $\PW_R$. 

The first result of this section, Lemma~\ref{linklemmalem}, provides a link between derivatives of $\wu$ in normal direction to the plane $\{x_3=0\}$ and derivatives in tangential direction to $\{x_3=0\}$. Subsequently, rather technical proofs lead to Corollary~\ref{analogontangential}, which provides analytic estimates for tangential derivatives of $\wu$. Finally, we combine Lemma~\ref{linklemmalem} and Corollary~\ref{analogontangential} to obtain the main result of this section, Lemma~\ref{bdryanalyticprefinal}, which provides the sought analytic estimates on $\wu$.

\subsection{Estimates on derivatives in normal direction}

The subsequent lemma is crucial, since it establishes a connection between the normal and tangential derivatives of a transformed solution $\wu$. The idea is quite simple: For a smooth vector field $\solw = (w_1, w_2, w_3)^T$ the definition of $\curl\solw$ reads as
\begin{align*}
\curl\solw = \left(\frac{\partial w_3}{\partial x_2}-\frac{\partial w_2}{\partial x_3}, \frac{\partial w_1}{\partial x_3}-\frac{\partial w_3}{\partial x_1},\frac{\partial w_2}{\partial x_1}-\frac{\partial w_1}{\partial x_2}\right)^T.
\end{align*}
From this, we observe that we can bound the derivatives of $w_1$ and $w_2$ in $x_3$-direction by $\curl\solw$ and the derivatives of $\solw$ in $x_1$ and $x_2$-directions. Moreover, if the divergence of $\solw$ is known, the derivative of $w_3$ in $x_3$-direction can be estimated by $\diverg\solw$ and the tangential derivatives of $w_1$ and $w_2$. 
In total, the normal derivative of $\solw$ can be estimated by $\curl\solw$, $\diverg\solw$ and the tangential derivatives of $\solw$.

The following result exploits this to provide a link between normal and tangential derivatives of $\wu$ in terms of anisotropic Sobolev-Morrey seminorms on half-balls, cf. Definition~\ref{anismorrey}.

\begin{lemma}\label{linklemmalem}
		Let $\Gp$ be an $\ana$-partition and let $\mu^{-1}, \varepsilon\in\anatenspw{\Gp}$ satisfy \eqref{coercivedomain}. In addition, suppose that $\solu\in\Hcurl$ is a piecewise smooth solution of Maxwell's equation \eqref{Maxwellorig} corresponding to a right-hand side $\solf$. 
		Furthermore, let $0<R\leq \Rmax$ and let $\Upsilon:B_R\rightarrow\R^3$ be a local flattening of the boundary $\Gamma$ or of an interface component $\I_j$. Under these assumptions, consider the transformed quantities $\tildf{\solu}$, $\hatf{\solf}$, $\wmu^{-1}$ and $\weps$ from Definition~\ref{trafoquantitiesdef}.
	 
	Then, there exist constants $C>0$ and $A\geq 1$ depending only on $\Rmax$, $\weps$ and $\wmu^{-1}$ such that $\wu$ and $\wf$ satisfy
	\begin{align*}
		\summorrey{\wu}{\ell,n+1,\PM_R}\leq  C\sum_{d=0}^{\ell+1}\sum_{b=0}^nA^{\ell+n+1-d-b}\left[\summorrey{\wu}{d,b,\PM_R}+|k|^{-2}\morrey{\wf}_{d,b,\PM_R}+|k|^{-2}\morreyone{\diverg\wf}_{d,b,\PM_R}\right]
	\end{align*}
	if $\Upsilon$ is a local flattening of $\Gamma$, and 
	\begin{align*}
	\summorrey{\wu}{\ell,n+1,\PS_R}\leq  C\sum_{d=0}^{\ell+1}\sum_{b=0}^nA^{\ell+n+1-d-b}\left[\summorrey{\wu}{d,b,\PS_R}+|k|^{-2}\morrey{\wf}_{d,b,\PS_R}+|k|^{-2}\morreyone{\diverg\wf}_{d,b,\PS_R}\right]
	\end{align*}
	if $\Upsilon$ is a local flattening of $\I_j$.
\end{lemma}

\begin{fatproof}
		We only consider the case that $\Upsilon:B_R\rightarrow\R^3$ is a local flattening of $\Gamma$; if $\Upsilon$ is a local flattening of an interface component $\I_j$, the proof follows the same lines by additionally considering terms involving $\PP_R$. 
	The remainder of the proof is divided into two steps.

	\medskip
	
	{\bf Step 1:} The first step is to find an appropriate estimate on $\morrey{\tildf{\solu}}_{\ell,n+1,\PM_R}$. To that end, we define a smooth tensor field $\lambda_{\weps}$ by 
	\begin{align*}
		\lambda_{\weps} = \begin{pmatrix}
		1 & 0 & 0 \\
		0 & 1 & 0 \\
		\weps_{31} & \weps_{32} & \weps_{33}
		\end{pmatrix}
	\end{align*}
	and set $\tildf\solw := (\wucomp_1, \wucomp_2, (\weps\wu)_3)^T$, where $\wucomp_1$ and $\wucomp_2$ are the first two components of $\wu$. We notice that there holds $\tildf\solw = \lambda_{\weps}\wu$ and, since $\weps$ is coercive, the tensor field $\lambda_{\weps}$ is invertible. Hence, Lemma \ref{commutatoranis} implies
	\begin{align}\label{linklemma1}
		\morrey{\tildf{\solu}}_{\ell,n+1,\PM_R} = \morrey{\lambda_{\weps}^{-1}\lambda_{\weps}\tildf{\solu}}_{\ell,n+1,\PM_R}\leq C\sum_{d=0}^{\ell}\sum_{b=0}^{n+1}K^{\ell+n+1-d-b}\morrey{\tildf{\solw}}_{d,b,\PM_R}
	\end{align}
    for appropriate constants $C>0$ and $K\geq 1$.

	For $d,b\in\N_0$, let $\alpha^* = (0,0,1)^T$, $\beta^* = (0,0,b)^T$ and $\gamma^* = (0,0,b+1)^T$. Moreover, let $\delta'\in\N_0^2\times\{0\}$ be a multiindex with $|\delta'|=d$ and $\tau:=d+b+|k|$.
	
	\medskip 
	
	In view of 
	\begin{align*}
	\curl\solw = \left(\frac{\partial w_3}{\partial x_2}-\frac{\partial w_2}{\partial x_3}, \frac{\partial w_1}{\partial x_3}-\frac{\partial w_3}{\partial x_1},\frac{\partial w_2}{\partial x_1}-\frac{\partial w_1}{\partial x_2}\right)^T
	\end{align*}
	for any smooth vector field $\solw = (w_1, w_2, w_3)^T$, we infer that the derivatives of $w_1$ and $w_2$ in $x_3$-direction can be estimated in terms of $\curl\solw$ and the derivatives of $\solw$ in the $x_1$ and $x_2$-directions.
	
	Applying this fact to $\D^{\delta'+\gamma^*}(\tildf{\solw})$ shows that there exists a $C>0$ independent of $\ell$ and $n$ such that 
	\begin{align}\label{normaltmpfund}
	\norm{\D^{\delta'+\gamma^*}(\tildf{w}_i)}{\PM_{R-(\tau+1)\rho}} \leq C\norm{\D^{\delta'+\beta^*}(\curl\tildf{\solu})}{\PM_{R-(\tau+1)\rho}}+C\sum_{|\alpha'|=1}\norm{\D^{\delta'+\alpha'+\beta^*}(\tildf{\solu})}{\PM_{R-(\tau+1)\rho}} 
	\end{align}
 for $i=1,2$, where we exploited $\tildf{w}_i = \tildf{u}_i$ for $i=1,2$. 
	
	\medskip 
	
	From the transformed Maxwell equations \eqref{trafosystem} we infer $\diverg \D^{\delta'+\beta^*}(\weps\wu) = -k^{-2}\diverg\D^{\delta'+\beta^*}(\hatf{\solf})$ in $\PM_R$, thus from the definition of $\tildf{\solw}$ we get
	\begin{align}\label{normaltmp1}
	\begin{split}
	\norm{\D^{\delta'+\gamma^*}(\tildf{w}_3)}{\PM_{R-(\tau+1)\rho}}\leq |k|^{-2}\norm{\diverg\D^{\delta'+\beta^*}(\hatf{\solf})}{\PM_{R-(\tau+1)\rho}}+\sum_{|\alpha'|=1}	\norm{\D^{\delta'+\alpha'+\beta^*}(\weps\tildf{\solu})}{\PM_{R-(\tau+1)\rho}}.
	\end{split}
	\end{align}
	
	Combining \eqref{normaltmpfund} and \eqref{normaltmp1} yields
	\begin{align*}
		\morrey{\tildf{\solw}}_{d,b+1,\PM_R}\leq C\left(|k|^{-1}\morrey{\curl\wu}_{d,b,\PM_R}+\morrey{\wu}_{d+1,b,\PM_R}+\morrey{\weps\wu}_{d+1,b,\PM_R}+|k|^{-2}\morreyone{\diverg\hatf{\solf}}_{d,b,\PM_R}\right),
	\end{align*}
	where we exploited $\rho\leq R|k|^{-1}$. Hence, some manipulations and Lemma \ref{commutatoranis} show 
	\begin{align*}
		\morrey{\tildf{\solw}}_{d,b+1,\PM_R}\leq C&|k|^{-2}\morreyone{\diverg\hatf{\solf}}_{d,b,\PM_R}+C\sum_{i=0}^{d+1}\sum_{j=0}^bL^{d+b+1-i-j}\summorrey{\wu}{i,j,\PM_R}
	\end{align*}
    for appropriate constants $C>0$ and $L\geq 1$.
	
	In total, Lemma \ref{vierersumme} proves 
	\begin{align*}
		\sum_{d=0}^{\ell}\sum_{b=1}^{n+1}K^{\ell+n+1-d-b}\morrey{\tildf{\solw}}_{d,b,\PM_R}\leq C\sum_{d=0}^{\ell+1}\sum_{b=0}^nB^{\ell+n+1-d-b}\left[\summorrey{\wu}{d,b,\PM_R}+|k|^{-2}\morreyone{\diverg\hatf{\solf}}_{d,b,\PM_R}\right]
	\end{align*}
	for some constants $C>0$ and $B\geq 1$. Furthermore, in combination with Lemma~\ref{commutatoranis} \revision{and Lemma~\ref{vierersumme} we get}
    \begin{align*}
    	\sum_{d=0}^{\ell}K^{\ell+n+1-d}\morrey{\tildf{\solw}}_{d,0,\PM_R}\leq C\sum_{d=0}^{\ell}G^{\ell+n+1-d}\morrey{\wu}_{d,0,\PM_R}
    \end{align*}
    for constants $C>0$ and $G\geq K$, hence
    \begin{align*}
    \sum_{d=0}^{\ell}\sum_{b=0}^{n+1}K^{\ell+n+1-d-b}\morrey{\tildf{\solw}}_{d,b,\PM_R}\leq C\sum_{d=0}^{\ell+1}\sum_{b=0}^nP^{\ell+n+1-d-b}\left[\summorrey{\wu}{d,b,\PM_R}+|k|^{-2}\morreyone{\diverg\hatf{\solf}}_{d,b,\PM_R}\right]
    \end{align*}
    for appropriate constants $C>0$ and $P\geq 1$.
	In combination with \eqref{linklemma1} this proves 
	\begin{align}\label{linklemmafinal}
		\morrey{\tildf{\solu}}_{\ell,n+1,\PM_R}\leq C\sum_{d=0}^{\ell+1}\sum_{b=0}^nA^{\ell+n+1-d-b}\left[\summorrey{\wu}{d,b,\PM_R}+|k|^{-2}\morreyone{\diverg\hatf{\solf}}_{d,b,\PM_R}\right]
	\end{align}
	for appropriate $C>0$ and $A\geq 1$.

	\bigskip 
	
	{\bf Step 2:} In the second step we consider the quantity $|k|^{-1}\morrey{\wmu^{-1}\curl\tildf{\solu}}_{\ell,n+1,\PM_R}$. Again, we define a tensor field $\lambda_{\wmu}$ by
	\begin{align*}
	\lambda_{\wmu} = \begin{pmatrix}
	1 & 0 & 0 \\
	0 & 1 & 0 \\
	\wmu_{31} & \wmu_{32} & \wmu_{33}
	\end{pmatrix}
	\end{align*}
	and set $\tildf{\solv}:=((\wmu^{-1}\curl\tildf{\solu})_1,(\wmu^{-1}\curl\tildf{\solu})_2,(\curl\tildf{\solu})_3)^T$. Note that $\wmu_{31}$, $\wmu_{32}$ and $\wmu_{31}$ are the entries of the third row of $\wmu$, where $\wmu$ is the inverse of the coercive tensor field $\wmu^{-1}$.
	
	We notice that $\tildf{\solv} = \lambda_{\wmu}(\wmu^{-1}\curl\tildf{\solu})$ and similarly as in Step 1, Lemma \ref{commutatoranis} shows
	\begin{align}\label{linklemma2}
	\morrey{\wmu^{-1}\curl\tildf{\solu}}_{\ell,n+1,\PM_R} \leq C\sum_{d=0}^{\ell}\sum_{b=0}^{n+1}K^{\ell+n+1-d-b}\morrey{\tildf{\solv}}_{d,b,\PM_R}.
	\end{align}
For $d,b\in\N_0$, let $\alpha^* = (0,0,1)^T$, $\beta^* = (0,0,b)^T$ and $\gamma^* = (0,0,b+1)^T$. Moreover, let $\delta'\in\N_0^2\times\{0\}$ be a multiindex with $|\delta'|=d$ and let $\tau:=d+b+|k|$. We write $\tildf{\solv} = (\tildf{v}_1, \tildf{v}_2, \tildf{v}_3)^T$ and by the same reasoning as in Step 1 we infer 
	\begin{align*}
	\norm{\D^{\delta'+\gamma^*}(\tildf{v}_m)}{\PM_{R-(\tau+1)\rho}} \leq C\norm{\D^{\delta'+\beta^*}(\curl\wmu^{-1}\curl\tildf{\solu})}{\PM_{R-(\tau+1)\rho}}+C\sum_{|\alpha'|=1}\norm{\D^{\delta'+\alpha'+\beta^*}(\wmu^{-1}\curl\tildf{\solu})}{\PM_{R-(\tau+1)\rho}} 
	\end{align*}
	 for $m=1,2$, where we exploited $\tildf{v}_m = (\wmu^{-1}\curl\tildf{\solu})_m$ for $m=1,2$. We note that $\curl\wmu^{-1}\curl\wu = \hatf{\solf}+k^2\weps\wu$, thus Lemma~\ref{commutatoranis} and $\rho \leq R|k|^{-1}$ show
	\begin{align*}
		|k|^{-1}\rho^{d+b+1}\norm{\D^{\delta'+\gamma^*}(\tildf{v}_m)}{\PM_{R-(\tau+1)\rho}}&\leq C|k|^{-2}\morrey{\hatf{\solf}}_{d,b,\PM_R}+C\morrey{\weps\wu}_{d,b,\PM_R}+C|k|^{-1}\morrey{\wmu^{-1}\curl\wu}_{d+1,b,\PM_R} \\
		&\leq C|k|^{-2}\morrey{\hatf{\solf}}_{d,b,\PM_R}+C\sum_{i=0}^{d+1}\sum_{j=0}^bK^{d+b+1-i-j}\summorrey{\wu}{i,j,\PM_R}
		\end{align*}
	for appropriate constants $C>0$ and $K\geq 1$, and $m=1,2$.
	
	\medskip 
	
	Moreover, due to $\diverg\curl\tildf{\solu}=0$ and Lemma \ref{commutatoranis} we get 
	\begin{align*}
	\begin{split}
	|k|^{-1}\rho^{d+b+1}\norm{\D^{\delta'+\gamma^*}(\tildf{v}_3)}{\PM_{R-(\tau+1)\rho}}&\leq C|k|^{-1}\rho^{d+b+1}\sum_{|\alpha'|=1}	\norm{\D^{\delta'+\alpha'+\beta^*}(\curl\tildf{\solu})}{\PM_{R-(\tau+1)\rho}} \\
	&\leq C\sum_{i=0}^{d+1}\sum_{j=0}^bK^{d+b+1-i-j}|k|^{-1}\morrey{\wmu^{-1}\curl\wu}_{i,j,\PM_R}.
	\end{split}
	\end{align*}
	In total, this proves 
	\begin{align*}
		|k|^{-1}\morrey{\tildf{\solv}}_{d,b+1,\PM_R}\leq C|k|^{-2}\morrey{\hatf{\solf}}_{d,b,\PM_R}+C\sum_{i=0}^{d+1}\sum_{j=0}^bK^{d+b+1-i-j}\summorrey{\wu}{i,j,\PM_R}.
	\end{align*}
	Together with \eqref{linklemma2} and Lemma \ref{vierersumme} \revision{we get}
	\begin{align}\label{linklemmafinal2}
		|k|^{-1}\morrey{\wmu^{-1}\curl\tildf{\solu}}_{\ell,n+1,\PM_R}\leq C\sum_{d=0}^{\ell+1}\sum_{b=0}^nB^{\ell+n+1-d-b}\left[\summorrey{\wu}{d,b,\PM_R}+|k|^{-2}\morrey{\hatf{\solf}}_{d,b,\PM_R}\right]
	\end{align}
    for appropriate constants $C>0$ and $B\geq 1$.
	Adding \eqref{linklemmafinal} and \eqref{linklemmafinal2} completes the proof.
	
\end{fatproof}

As a consequence of Lemma~\ref{linklemmalem} and Lemma~\ref{vierersumme} we obtain the subsequent corollary. 

\begin{corollary}\label{auxiliartrace}
		Let $\Gp$ be an $\ana$-partition, let the coefficients $\mu^{-1}, \varepsilon\in\anatenspw{\Gp}$ satisfy \eqref{coercivedomain} and suppose that $\solu\in\Hcurl$ is a piecewise smooth solution of Maxwell's equation \eqref{Maxwellorig} corresponding to a right-hand side $\solf$.  
		Furthermore, let $0<R\leq \Rmax$ and let $\Upsilon:B_R\rightarrow\R^3$ be a local flattening of the boundary $\Gamma$. Under these assumptions we consider the transformed quantities $\wmu^{-1},\weps, \tildf{\solu}$ and $\hatf{\solf}$ from Definition~\ref{trafoquantitiesdef}.
	 
		\medskip

		Moreover, suppose that $\nu:\PM_R\rightarrow\Co^{3\times3}$ is a tensor field which satisfies $\nu\in\anatens{\PM_R}$, that is, there exist $\omega\geq 0, M>0$ such that
	\begin{align*}
			\forall \ell\in\N_0:\ \sum_{i,j=1}^3\sum_{|\alpha|=\ell}\norm{\D^{\alpha}\nu_{i,j}}{\SL(\PM_R)}\leq \omega M^{\ell}\ell^{\ell}.
	\end{align*}

	Finally, let $\beta'\in\N_0^2\times\{0\}$ be arbitrary, set $\ell:=|\beta'|$ and suppose that $k\in\Co\setminus\{0\}$ is given. Then, for all $0<\rho\leq \frac{R}{2(\ell+1+|k|)}$ there holds 
	\begin{multline*}
			\rho^{\ell+1}\sum_{|\alpha|=1}\norm{\D^{\alpha}\commute(\nu,\wu,\beta')}{\PM_{R-(\ell+|k|)\rho}}\leq C\sum_{d=0}^{\ell}A^{\ell-d}\left[\summorrey{\wu}{d,0,\PM_R}+|k|^{-2}\morrey{\wf}_{d,0,\PM_R}\right]\\
																																						 +C|k|^{-2}\sum_{d=0}^{\ell}A^{\ell-d}\morreyone{\diverg\wf}_{d,0,\PM_R},
	\end{multline*}
	as well as 
	\begin{multline*}
			|k|^{-1}\rho^{\ell+1}\sum_{|\alpha|=1}\norm{\D^{\alpha}\commute(\nu,\wmu^{-1}\curl\wu,\beta')}{\PM_{R-(\ell+|k|)\rho}}\leq C\sum_{d=0}^{\ell}A^{\ell-d}\left[\summorrey{\wu}{d,0,\PM_R}+|k|^{-2}\morrey{\wf}_{d,0,\PM_R}\right]\\
																																						 +C|k|^{-2}\sum_{d=0}^{\ell}A^{\ell-d}\morreyone{\diverg\wf}_{d,0,\PM_R}.
	\end{multline*}
   In both cases, the constants $C>0$ and $A\geq 1$ depend only on $\Rmax, \wmu^{-1}, \weps, \omega$ and $M$.
\end{corollary}

\begin{fatproof}
		We prove only the first inequality, the second one follows \revision{from similar arguments}.
		If $\alpha$ is a tangential derivative, the claim follows directly from Lemma~\ref{commutatoranis}. Hence, we only have to discuss the case $\alpha = (0,0,1)^T$ in more detail. According to Lemma~\ref{commutatoranis} we have
		\begin{align}\label{commutatoranisapplied}
		\rho^{\ell+1}\norm{\D^{\alpha}(\commute(\nu,\wu,\beta'))}{\PM_{R-(\ell+|k|)\rho}} \leq C\sum_{d=0}^{\ell}K^{\ell-d}\morrey{\wu}_{d,0,\PM_R}
																													  +C\sum_{d=0}^{\ell-1}K^{\ell-1-d}\morrey{\wu}_{d,1,\PM_R}.
		\end{align}
		Furthermore, due to \revision{the trivial estimate} $\morrey{\wu}_{d,1,\PM_R}\leq \summorrey{\wu}{d,1,\PM_R}$ \revision{we get from Lemma~\ref{linklemmalem} that} 
		\begin{align*}
		\morrey{\wu}_{d,1,\PM_R}\leq C\sum_{m=0}^{d+1}L^{d+1-m}\left[\summorrey{\wu}{m,0,\PM_R}+|k|^{-2}\morrey{\wf}_{m,0,\PM_R}+|k|^{-2}\morreyone{\diverg\wf}_{m,0,\PM_R}\right],
		\end{align*}
		thus Lemma~\ref{vierersumme} yields
		\begin{align*}
		\sum_{d=0}^{\ell-1}K^{\ell-1-d}\morrey{\wu}_{d,1,\PM_R}\leq C\sum_{d=0}^{\ell}A^{\ell-d}\left[\summorrey{\wu}{d,0,\PM_R}+|k|^{-2}\morrey{\wf}_{d,0,\PM_R}+|k|^{-2}\morreyone{\diverg\wf}_{d,0,\PM_R}\right].
		\end{align*}
		Together with \eqref{commutatoranisapplied} this implies the claim.	
\end{fatproof}

%% file: analytic/analytic_tangential.tex
\subsection{Estimates on derivatives in tangential direction}

After having established a connection between normal and tangential derivatives of $\wu$, it remains to provide appropriate estimates for the growth of the tangential derivatives.
Our aim is to prove Lemma~\ref{tangentfinallemma}, which is an analog of Lemma~\ref{analyticinteriorprelemma} and from which we can infer the crucial Corollary~\ref{analogontangential}. The general idea is to employ similar techniques as in Section~\ref{analyticinteriorsec}, however, the presence of boundary terms makes this more involved. In an attempt to ease the presentation we single out the most technical parts of of the proof of Lemma~\ref{tangentfinallemma} and provide them beforehand in the subsections~\ref{sub1}-\ref{sub5}. 

Subsequently, both the Sobolev-Morrey seminorms on half-balls from Definition~\ref{anismorrey} and the Sobolev-Morrey seminorms on discs from Definition~\ref{anismorreybdr} play a prominent role. In addition, let us recall that for $0<\delta<R$, the symbol $\chi_{R,\delta}$ denotes a smooth cutoff-function with the properties stated in Definition~\ref{cutoffdef}.

\medskip

\subsubsection{Auxiliary estimates for volume terms}\label{sub1}

We start by providing inequalities that control certain volume terms on half-balls. 
\begin{lemma}\label{curldivbdry}
		Let $\Gp$ be an $\ana$-partition and assume that $\mu^{-1}, \varepsilon\in\anatenspw{\Gp}$ satisfy \eqref{coercivedomain}.
		Moreover, let $0<R\leq\Rmax$ and let $\Upsilon:\PW_R\rightarrow\R^3$ be a local flattening of the boundary $\Gamma$ or of an interface component $\I_j$. 
	
	Under these assumptions, consider a piecewise smooth solution $\solu\in\Hcurl$ of Maxwell's equations \eqref{Maxwellorig} corresponding to a right-hand side $\solf$, and let $\wu$, $\hatf{\solf}$, $\wmu^{-1}$ and $\weps$ be the transformed quantities from Definition~\ref{trafoquantitiesdef}. Then, for any $\ell\in\N_0$, any $\beta'\in\N_0^2\times\{0\}$ with $|\beta'|=\ell$ and any $0<\rho<\frac{R}{2(\ell+1+|k|)}$ there holds with $\tau:=\ell+|k|$
	\begin{align}\label{divestbdr}
		\rho^{\ell+1}\norm{\curl\chi_{R-\tau\rho,\rho} \D^{\beta'}(\wu)}{\PM_{R-\tau\rho}}\leq C\sum_{d=0}^{\ell}A^{\ell-d}\summorrey{\wu}{d,0,\PM_R},
	\end{align}
	as well as 
	\begin{equation}\label{curlestbdr2}
	\begin{multlined}	
	\rho^{\ell+1}\norm{\diverg\chi_{R-\tau\rho,\rho}\weps\D^{\beta'}(\tildf{\solu})}{\PM_{R-\tau\rho}} \\
	\leq 
	C\sum_{d=0}^{\ell}A^{\ell-d}\left[\summorrey{\wu}{d,0,\PM_R}+|k|^{-2}\morrey{\hatf{\solf}}_{d,0,\PM_R}+|k|^{-2}\morreyone{\diverg\hatf{\solf}}_{d,0,\PM_R}\right],
	\end{multlined}
\end{equation}
	where the constants $C>0$ and $A\geq 1$ depend only on $\Rmax$, $\wmu^{-1}$ and $\weps$. 
	
	Furthermore, we have 
	\begin{align}\label{divcurlbdr}
	|k|^{-1}\rho^{\ell+1}\norm{\curl\chi_{R-\tau\rho,\rho}\D^{\beta'}(\wmu^{-1}\curl\wu)}{\PM_{R-\tau\rho}}\leq C\sum_{d=0}^{\ell}A^{\ell-d}\left[\summorrey{\wu}{d,0,\PM_R}+|k|^{-2}\morrey{\hatf{\solf}}_{d,0,\PM_R}\right],
	\end{align}
	and 
	\begin{equation}\label{curlbdr44}
	\begin{multlined}
			|k|^{-1}\rho^{\ell+1}\norm{\diverg\chi_{R-\tau\rho,\rho}\wmu\D^{\beta'}(\wmu^{-1}\curl\wu)}{\PM_{R-\tau\rho}} \\\leq C\sum_{d=0}^{\ell}A^{\ell-d}\left[\summorrey{\wu}{d,0,\PM_R}+|k|^{-2}\morrey{\hatf{\solf}}_{d,0,\PM_R}
	+|k|^{-2}\morreyone{\diverg\hatf{\solf}}_{d,0,\PM_R}\right],
	\end{multlined}
	\end{equation}
	where, again, the constants $C>0$ and $A\geq 1$ depend only on $\Rmax$, $\wmu^{-1}$ and $\weps$. If $\Upsilon$ is a local transformation of $\I_j$, the inequalities \eqref{divestbdr}-\eqref{curlbdr44} also hold for $\PM_R$ and $\PM_{R-\tau\rho}$ replaced by $\PP_R$ and $\PP_{R-\tau\rho}$, respectively.
\end{lemma}
\begin{fatproof}
		The idea of the proof of any of the estimates \eqref{divestbdr}-\eqref{curlbdr44} \revision{is to combine} the product rule with the commutator estimate provided by Lemma \ref{commutatoranis}.
	We explain this by giving a proof for \eqref{divestbdr} and \eqref{curlestbdr2}; proving \eqref{divcurlbdr} and \eqref{curlbdr44} is done similarly.

	\medskip 

	The product rule, the properties of $\chi_{R-\tau\rho,\rho}$ and Lemma \ref{commutatoranis} show
	\begin{align*}
		\rho^{\ell+1}\norm{\curl\chi_{R-\tau\rho,\rho} \D^{\beta'}(\wu)}{\PM_{R-\tau\rho}}&\leq |k|^{-1}\rho^{\ell}\norm{\curl \D^{\beta'}(\wu)}{\PM_{R-\tau\rho}}+\rho^{\ell}\norm{\D^{\beta'}(\wu)}{\PM_{R-\tau\rho}} \\
		&\leq C\sum_{d=0}^{\ell}A^{\ell-d}\summorrey{\wu}{d,0,\PM_R},
	\end{align*}
	which proves \eqref{divestbdr}.
	
	\bigskip
	
	To prove \eqref{curlestbdr2}, we note that the product rule and $\diverg\weps\wu = -k^{-2}\diverg\wf$ yield
	\begin{align*}
	\begin{split}
			\rho^{\ell+1}\norm{\diverg\chi_{R-\tau\rho,\rho}\weps\D^{\beta'}(\tildf{\solu})}{\PM_{R-\tau\rho}} &\leq \rho^{\ell+1}\norm{\diverg\weps\D^{\beta'}(\tildf{\solu})}{\PM_{R-\tau\rho}}+C\rho^{\ell}\norm{\weps\D^{\beta'}(\tildf{\solu})}{\PM_{R-\tau\rho}}\\
	&\leq \rho^{\ell+1}\norm{\diverg\commute(\weps,\wu,\beta')}{\PM_{R-\tau\rho}}+|k|^{-2}\morreyone{\diverg\wf}_{\ell,0,\PM_R}\\
	&\hskip 6cm +C\rho^{\ell}\norm{\weps\D^{\beta'}(\tildf{\solu})}{\PM_{R-\tau\rho}}.
	\end{split}
	\end{align*}
	Moreover, Lemma~\ref{commutatoranis} proves 
	\begin{align*}
		\rho^{\ell+1}\norm{\diverg\commute(\weps,\wu,\beta')}{\PM_{R-\tau\rho}}+C\rho^{\ell}\norm{\weps\D^{\beta'}(\tildf{\solu})}{\PM_{R-\tau\rho}}\leq C\sum_{d=0}^{\ell}K^{\ell-d}\morrey{\wu}_{d,0,\PM_R}+C\sum_{d=0}^{\ell-1}K^{\ell-1-d}\morrey{\wu}_{d,1,\PM_R}
	\end{align*}
	for appropriate constants $C>0$ and $K\geq 1$. As in the proof of Corollary~\ref{auxiliartrace} we obtain
	\begin{align*}
		\sum_{d=0}^{\ell-1}K^{\ell-1-d}\morrey{\wu}_{d,1,\PM_R}\leq C\sum_{d=0}^{\ell}B^{\ell-d}\left[\summorrey{\wu}{d,0,\PM_R}+|k|^{-2}\morrey{\wf}_{d,0,\PM_R}+|k|^{-2}\morreyone{\diverg\wf}_{d,0,\PM_R}\right]
	\end{align*}
	for some constants $C>0$ and $B\geq 1$.
	In total, we end up with 
	\begin{multline*}
			\rho^{\ell+1}\norm{\diverg\chi_{R-\tau\rho,\rho}\weps\D^{\beta'}(\tildf{\solu})}{\PM_{R-\tau\rho}}\\
			\leq C\sum_{d=0}^{\ell}A^{\ell-d}\left[\summorrey{\wu}{d,0,\PM_R}+|k|^{-2}\morrey{\wf}_{d,0,\PM_R}+|k|^{-2}\morreyone{\diverg\wf}_{d,0,\PM_R}\right]
	\end{multline*}
	for appropriate constants $C>0$ and $A\geq 1$, and this proves \eqref{curlestbdr2}.

\end{fatproof}

\subsubsection{Auxiliary estimates for boundary terms}\label{sub2}

Our next aim is to provide further auxiliary estimates for boundary terms. Henceforth, in order to ease our notation we abbreviate

\begin{align*}
		\norm{\solv}{\thalfm,\PZ_R}:=\norm{\solv}{\VSHMZ{-1/2}_T(\PZ_R)}\quad\text{and}\quad  \norm{\solw}{\thalf,\PZ_R}:=\norm{\solw}{\VSHMZ{1/2}_T(\PZ_R)}
\end{align*}
for all $\solv\in\VSHMZ{-1/2}_T(\PZ_R)$ and all $\solw\in\VSHMZ{1/2}_T(\PZ_R)$, respectively.
Analogously, for all $f\in\SHMZ{-1/2}(\PZ_R)$ and all $g\in\SHMZ{1/2}(\PZ_R)$ we set
\begin{align*}
		\norm{f}{\thalfm,\PZ_R}:=\norm{f}{\SHMZ{-1/2}(\PZ_R)}\quad\text{and}\quad  \norm{g}{\thalf,\PZ_R}:=\norm{g}{\SHMZ{1/2}(\PZ_R)}.
\end{align*}

\begin{remark}\label{rem:wzetaextension}
		According to Lemma~\ref{trafotensors}, if $\zeta\in\anagamma$ then $\wzeta\in\anagammadisc{\PZ_R}$. As a consequence, there exists a constant $R'>0$ depending only on $\wzeta$ such that $\wzeta$ can be extended to an analytic and coercive tensor field on $\PW_{R'}$. Extending $\wzeta$ to $\PW_{R'}$ will be necessary in the proof of Lemma~\ref{devialemma} below. Subsequently, we will apply Lemma~\ref{devialemma} only in the case $0<\Rmin\leq R\leq \Rmax$ for sufficiently small $\Rmax>0$, and therefore we may assume without loss of generality that $R'\geq \Rmax$. This implies that for $R\leq\Rmax$ the impedance tensor $\wzeta$ can be extended to a coercive $\wzeta\in\anatens{\PW_R}$.
\end{remark}

The following lemma provides auxiliary estimates for boundary terms. 

\begin{lemma}\label{devialemma}
		Let $\Gp$ be an $\ana$-partition, assume that the coefficients $\mu^{-1}, \varepsilon\in\anatenspw{\Gp}$ satisfy \eqref{coercivedomain} and that $\zeta\in\anagamma$ satisfies \eqref{coercivebdr} and \eqref{zetaproperty}.
		Moreover, let $0<\Rmin\leq R\leq\Rmax$ and let $\Upsilon:\PW_R\rightarrow\R^3$ be a local flattening of the boundary $\Gamma$. Finally, let $\solu\in\Hcurl$ be a piecewise smooth solution of Maxwell's equations \eqref{Maxwellorig} corresponding to a right-hand side $\solf$ and impedance data $\solgi$. 

		Under these assumptions we consider the transformed quantities $\wu$, $\hatf{\solf}$, $\solgicheck$, $\wmu^{-1}$, $\weps$ and $\wzeta$ from Definition~\ref{trafoquantitiesdef}.
		Furthermore, let $\ell\in\N_0$ and $\beta'\in\N_0^2\times\{0\}$ with $|\beta'|=\ell$ be arbitrary, define the quantity $\tau:=\ell+|k|$ and let $\rho\in\left(0,\frac{R}{2(\ell+1+|k|)}\right)$.
		Then, with $\wv:=\chi_{R-\tau\rho,\rho}\D^{\beta'}(\wu)$ and $\ww:=\chi_{R-\tau\rho,\rho}\D^{\beta'}(\wmu^{-1}\curl\wu)$ there holds 
	
		\begin{align}\label{divbdrdevia}
	\begin{split}
			\norm{\divergt\wzeta\wv_T}{\revision{\thalfm,\PZ_{R-\tau\rho}}}+|k|^{-1}\norm{\curlt\wzeta^{-1}\ww_t}{\revision{\thalfm,\PZ_{R-\tau\rho}}}\leq &C\rho^{-\ell-1}\sum_{d=0}^{\ell}K^{\ell-d}\summorrey{\wu}{d,0,\PM_R}\\
																																											 &+C|k|^{-2}\rho^{-\ell-1}\sum_{d=0}^{\ell}K^{\ell-d}M_{k,d,R}(\wf,\solgicheck),
	\end{split}
	\end{align}
	where 
	\begin{align*}
			M_{k,d,R}(\wf,\solgicheck):=\morrey{\hatf{\solf}}_{d,0,\PM_R}+\morreyone{\diverg\hatf{\solf}}_{d,0,\PM_R}+|k|\morreyhalf{\solgicheck}_{d,\PZ_R}.
	\end{align*}
	The constants $C>0$ and $K\geq 1$ depend only on $\Rmin$, $\Rmax$, $\wmu^{-1},\weps$ and $\wzeta$.
\end{lemma}

\begin{fatproof}
		We notice that the impedance boundary condition $\wmu^{-1}\curl\wu\times\sole_3-ik\wzeta\wu_T = \solgicheck$ from Lemma~\ref{trafolemma} \revision{implies}
	\begin{align}\label{lemma44i0}
	\ww_t-ik\wzeta\wv_T +ik\chi_{R-\tau\rho,\rho}\commute(\wzeta,\wu_T,\beta')= \chi_{R-\tau\rho,\rho}\D^{\beta'}(\solgicheck),
	\end{align}
	where $\ww_t:=\ww\times\sole_3$ and $\wu_T:=\sole_3\times(\wu\times\sole_3)$ denote the tangential trace of $\ww$ and the tangential component of $\wu$, respectively.
	Hence, \revision{by using the mapping properties \eqref{divtmapping1}-\eqref{divtmapping2} we arrive at}
	\begin{align}\label{lemma44i1}
			\begin{split}
					\rho^{\ell+1}\norm{\divergt\wzeta\wv_T}{\revision{\thalfm,\PZ_{R-\tau\rho}}}\leq &\rho^{\ell+1}\left(|k|^{-1}\norm{\divergt\ww_t}{\revision{\thalfm,\PZ_{R-\tau\rho}}}\right.+\norm{\chi_{R-\tau\rho,\rho}\commute(\wzeta,\solu_T,\beta')}{\revision{\thalf,\PZ_{R-\tau\rho}}}\\ 
																									 &\hskip 4.7cm +\left.|k|^{-1}\rho^{\ell+1}\norm{\chi_{R-\tau\rho,\rho}\D^{\beta'}(\solgicheck)}{\thalf,\PZ_{R-\tau\rho}}\right).		  %&\hskip 7cm +C|k|^{-1}\morreyhalf{\solgicheck}_{\ell,\PZ_R}.
	\end{split}
	\end{align}
	\revision{Lemma~\ref{interpolation} shows
			\begin{align*}
					\rho^{\ell+1}\norm{\chi_{R-\tau\rho,\rho}\D^{\beta'}(\solgicheck)}{\thalf,\PZ_{R-\tau\rho}}\leq C\rho^{\ell+1/2}\norm{\D^{\beta'}(\solgicheck)}{\frac{1}{2},\PZ_{R-\tau\rho}}\leq \morreyhalf{\solgicheck}_{\ell,\PZ_R}.
			\end{align*}
	}

	As a consequence of Lemma~\ref{lem:trbdr}, the mapping properties \eqref{divtmapping1}-\eqref{divtmapping2} and the product rule we get
	\begin{align}\label{lemma44i2}
			\begin{split}
					\rho^{\ell+1}|k|^{-1}\norm{\divergt\ww_t}{\thalfm,\PZ_{R-\tau\rho}}&\leq C\left(\summorrey{\wu}{\ell,0,\PM_R}+|k|^{-1}\rho^{\ell+1}\norm{\D^{\beta'}(\curl\wmu^{-1}\curl\wu)}{\PM_{R-\tau\rho}}\right) \\
																						   &\leq C\left(\summorrey{\wu}{\ell,0,\PM_R}+|k|^{-2}\morrey{\hatf{\solf}}_{\ell,0,\PM_R}\right).
	\end{split}
	\end{align}

	We aim to apply the trace inequality to the commutator term. To that end, we notice that there holds the identity $\commute(\wzeta,\solu_T,\beta')=\commute(\wzeta,\solu,\beta')_T$. A priori, the transformed impedance tensor $\wzeta$ is only defined on $\PZ_R$ and not on $\PM_R$. However, as explained in Remark~\ref{rem:wzetaextension} we may extend $\wzeta$ to a coercive tensor field $\wzeta\in\anatens{\PW_R}$, then a trace inequality and $\rho\leq R$ show
	\begin{align*}
			\rho^{\ell+1}\norm{\chi_{R-\tau\rho,\rho}\commute(\wzeta,\solu_T,\beta')}{\revision{\thalf,\PZ_{R-\tau\rho}}}\leq C\rho^{\ell} \norm{\commute(\wzeta,\solu,\beta')}{\PM_{R-\tau\rho}}+C\rho^{\ell+1}\sum_{|\alpha|=1}\norm{D^{\alpha}\left(\commute(\wzeta,\solu,\beta')\right)}{\PM_{R-\tau\rho}}.
	\end{align*}
	Consequently, by Lemma~\ref{commutatoranis} and Corollary~\ref{auxiliartrace} we conclude
	\begin{align}\label{lemma44i3}
		\rho^{\ell+1}\norm{\chi_{R-\tau\rho,\rho}\commute(\wzeta,\solu_T,\beta')}{\thalf,\PZ_{R-\tau\rho}}\leq C\sum_{d=0}^{\ell}K^{\ell-d}\summorrey{\wu}{d,0,\PM_R}+C|k|^{-2}\sum_{d=0}^{\ell}K^{\ell-d}M_{k,d,R}(\wf,\solgicheck),
\end{align}
and combining \eqref{lemma44i1}-\eqref{lemma44i3} then proves the claim concerning $\divergt\wzeta\wv_T$.

\medskip

As for the statement concerning $\curlt\wzeta^{-1}\ww_t$, we notice that \eqref{lemma44i0} yields 
\begin{align}\label{lemma44i4}
		\begin{split}
		|k|^{-1}\rho^{\ell+1}\norm{\curlt\wzeta^{-1}\ww_t}{\thalfm,\PZ_{R-\tau\rho}}\leq\rho^{\ell+1}&\left(\norm{\curlt\wv_T}{\thalfm,\PZ_{R-\tau\rho}}+\norm{\chi_{R-\tau\rho,\rho}\commute(\wzeta,\solu_T,\beta')}{\thalf,\PZ_{R-\tau\rho}}\right)\\
		&\hskip 5cm +C|k|^{-1}\morreyhalf{\solgicheck}_{\ell,\PZ_R}
\end{split}
\end{align}
The commutator term was already estimated in \eqref{lemma44i3}, and for the term involving $\curlt\wv_T$ we notice that the identity \eqref{surfacetraceprop}, the product rule and the fact that $\rho\leq R|k|^{-1}$ show
\begin{align}\label{lemma44tmp}
		\rho^{\ell+1}\norm{\curlt\wv_T}{\thalfm,\PZ_{R-\tau\rho}}\leq C\left(|k|^{-1}\morrey{\curl\wu}_{\ell,0,\PM_R}+\morrey{\wu}_{\ell,0,\PM_R}\right).
\end{align}
Finally, according to Lemma~\ref{commutatoranis} we have
\begin{align*}
		\morrey{\curl\wu}_{\ell,0,\PM_R} \leq\morrey{\wmu\ \wmu^{-1}\curl\wu}_{\ell,0,\PM_R}\leq C\sum_{d=0}^{\ell}K^{\ell-d}\morrey{\wmu^{-1}\curl\wu}_{d,0,\PM_R}, 
\end{align*}
hence \eqref{lemma44tmp} implies
\begin{align}\label{lemma44i5}
		\rho^{\ell+1}\norm{\curlt\wv_T}{\thalfm,\PZ_{R-\tau\rho}}\leq C\sum_{d=0}^{\ell}K^{\ell-d}\summorrey{\wu}{d,0,\PM_R}.
\end{align}
Combining \eqref{lemma44i5} with \eqref{lemma44i4} and \eqref{lemma44i3} proves the statement concerning $\curlt\wzeta^{-1}\ww_t$ and thus concludes the proof of \eqref{divbdrdevia}. 

\end{fatproof}

%% file: analytic/analytic_tangential_2.tex
\subsubsection{Auxiliary estimates for interface terms}\label{sub5}

So far, we derived auxiliary estimates in the case of $\Upsilon$ being a local flattening of the boundary $\Gamma$. It remains to consider the case of $\Upsilon$ being a local flattening of an interface component $\interf_j$. For a piecewise smooth vector field $\solv$ on $\PW_R$ we recall the normal and tangential jumps $\njump{\solv}$ and $\tjump{\solv}$ from Section~\ref{interfacesubsec}.
With these definitions there holds the following auxiliary result:

\begin{lemma}\label{interfaceest}
		Let $\Gp$ be an $\ana$-partition and assume that $\mu^{-1},\varepsilon\in\anatenspw{\Gp}$ satisfy \eqref{coercivedomain}.
		In addition, let $0<\Rmin\leq R\leq \Rmax$, let $\Upsilon:\PW_R\rightarrow\R^3$ be a local flattening of an interface component $\interf_j$, and let $\solu\in\Hcurl$ be a piecewise smooth solution of Maxwell's equation \eqref{Maxwellorig} corresponding to a right-hand side $\solf\in\Hdiv$. Furthermore, let $\wu$, $\hatf{\solf}$, $\wmu^{-1}$ and $\weps$ be the transformed quantities from Definition~\ref{trafoquantitiesdef}.

		Under these assumptions, let $\ell\in\N_0$ and $\beta'\in\N_0^2\times\{0\}$ with $|\beta'|=\ell$ be arbitrary, define $\tau:=\ell+|k|$ and choose $\rho\in\left(0,\frac{R}{2(\ell+1+|k|)}\right)$.
		Then, there holds 
		\begin{equation}\label{interf1}
				\begin{multlined}
						\rho^{\ell+1}\norm{\njump{\chi_{R-\tau\rho,\rho}\weps\D^{\beta'}(\wu)}}{\thalf,\PZ_{R-\tau\rho}}\leq C\sum_{d=0}^{\ell}K^{\ell-d}\left[\summorrey{\wu}{d,0,\PS_R}+|k|^{-2}\morrey{\wf}_{d,0,\PS_R}\right]\\ 
																																					+C|k|^{-2}\sum_{d=0}^{\ell}K^{\ell-d}\morreyone{\diverg\wf}_{d,0,\PS_R}
		\end{multlined}
	\end{equation}
as well as 
\begin{equation}\label{interf2}
		\begin{multlined}
		|k|^{-1}\rho^{\ell+1}\norm{\njump{\chi_{R-\tau\rho,\rho}\wmu\D^{\beta'}(\wmu^{-1}\curl\wu)}}{\thalf,\PZ_{R-\tau\rho}}\leq C\sum_{d=0}^{\ell}K^{\ell-d}\left[\summorrey{\wu}{d,0,\PS_R}+|k|^{-2}\morrey{\wf}_{d,0,\PS_R}\right] \\
+C|k|^{-2}\sum_{d=0}^{\ell}K^{\ell-d}\morreyone{\diverg\wf}_{d,0,\PS_R},
\end{multlined}
\end{equation}
where the constants $C>0$ and $K\geq 1$ depend only on $\Rmin$, $\Rmax$, $\wmu^{-1}$ and $\weps$.

\end{lemma}

\begin{fatproof}
		From the transformed Maxwell equations in Lemma~\ref{trafolemma} and $\solf\in\Hdiv$ we infer $\weps\wu\in\Hdivdomain{\PW_R}$, and therefore $\njump{\chi_{R-\tau\rho,\rho}\D^{\beta'}(\weps\wu)}=0$ on $\PZ_R$. Hence, a trace inequality yields
	\begin{align*}
		\rho^{\ell+1}\norm{\njump{\chi_{R-\tau\rho,\rho}\weps\D^{\beta'}(\wu)}}{\thalf,\PZ_{R-\tau\rho}}&=\rho^{\ell+1}\norm{\chi_{R-\tau\rho,\rho}\njump{\commute(\weps,\wu,\beta')}}{\thalf,\PZ_{R-\tau\rho}} \\
																											 &\leq C\rho^{\ell}\norm{\commute(\weps,\wu,\beta')}{\PS_{R-\tau\rho}}+C\rho^{\ell+1}
\sum_{|\alpha| =  1}\norm{\D^{\alpha}(\commute(\weps,\wu,\beta'))}{\PS_{R-\tau\rho}}
	\end{align*}
	for a constant $C>0$ depending only on $\Rmin$ and $\Rmax$, where we exploited the product rule and $\rho\leq R$. Analogous arguments as in the proof of \eqref{lemma44i3} then show 
\begin{multline*}
		\rho^{\ell+1}\norm{\njump{\chi_{R-\tau\rho,\rho}\weps\D^{\beta'}(\wu)}}{\thalf,\PZ_{R-\tau\rho}}\leq C\sum_{d=0}^{\ell}K^{\ell-d}\left[\summorrey{\wu}{d,0,\PS_R}+|k|^{-2}\morrey{\wf}_{d,0,\PS_R}\right]\\ 
																																					+C|k|^{-2}\sum_{d=0}^{\ell}K^{\ell-d}\morreyone{\diverg\wf}_{d,0,\PS_R},
	\end{multline*}
	which concludes the proof of \eqref{interf1}.

	\medskip 

	In order to see \eqref{interf2}, we notice that due to $\curl\wu\in\Hdivdomain{\PW_R}$ there holds $\njump{\chi_{R-\tau\rho,\rho}\D^{\beta'}(\curl\wu)}=0$ on $\PZ_R$, hence a trace inequality and the product rule yield 
\begin{align*}
		\rho^{\ell+1}\norm{\njump{\chi_{R-\tau\rho,\rho}\wmu\D^{\beta'}(\wmu^{-1}\curl\wu)}}{\thalf,\PZ_{R-\tau\rho}}&=\rho^{\ell+1}\norm{\chi_{R-\tau\rho,\rho}\njump{\commute(\wmu,\wmu^{-1}\curl\wu,\beta')}}{\thalf,\PZ_{R-\tau\rho}} \\
																														  &\leq C\rho^{\ell} \norm{\commute(\wmu,\wmu^{-1}\curl\wu,\beta')}{\PS_{R-\tau\rho}} \\
																														  & \hskip 2cm +C\rho^{\ell+1}\sum_{|\alpha| = 1}\norm{\D^{\alpha}(\commute(\wmu,\wmu^{-1}\curl\wu,\beta'))}{\PS_{R-\tau\rho}}
	\end{align*}
for a constant $C>0$ depending only on $\Rmin$ and $\Rmax$,	where we again exploited that $\rho\leq R$.
	Similar arguments as in the proof of \eqref{lemma44i3} then show
\begin{multline*}
			\rho^{\ell+1}\norm{\njump{\chi_{R-\tau\rho,\rho}\wmu\D^{\beta'}(\wmu^{-1}\curl\wu)}}{\thalf,\PZ_{R-\tau\rho}}\leq C\sum_{d=0}^{\ell}K^{\ell-d}\left[\summorrey{\wu}{d,0,\PS_R}+|k|^{-2}\morrey{\wf}_{d,0,\PS_R}\right] \\ 
+C|k|^{-2}\sum_{d=0}^{\ell}K^{\ell-d}\morreyone{\diverg\wf}_{d,0,\PS_R},
	\end{multline*}
	which finishes the proof of \eqref{interf2}.
\end{fatproof}

\subsubsection{Estimates on derivatives in tangential direction: The main result}

The subsequent result is the analog to Lemma~\ref{analyticinteriorprelemma} and provides a recursive inequality for tangential derivatives of the transformed solution $\wu$.
For $\ell\in\N_0$ and $R>0$ we abbreviate  
\begin{align*}
		\solF_{k,\ell,R}:= \morrey{\wf}_{\ell,0,\PM_R}+\morreyone{\diverg\wf}_{\ell,0,\PM_R},
\end{align*}
and with this abbreviation, there holds the following result:
%and, depending on the imposed boundary conditions, 
%\begin{align*}
%		&\solG_{k,\ell,R}:= |k|\morreyhalf{\solgicheck}_{\ell,\PZ_R}\quad && {\rm in\ the\ case\ of\ } \eqref{Maxwellimpedance}, \\
%		&\solG_{k,\ell,R}:= |k|\morreyhalf{\solgncheck}_{\ell,\PZ_R}+\morreyhalf{\diverg_{\PZ_R}\solgncheck}_{\ell,\PZ_R}+\morreyhalf{\wf\cdot\sole_3}_{\ell,\PZ_R} \quad && {\rm in\ the\ case\ of\ }\eqref{Maxwellnatural},
% \\
%		&\solG_{k,\ell,R}:= 0 \quad && {\rm in\ the\ case\ of\ }\eqref{Maxwellessential}, 
%\end{align*}
%where $\sole_3 := (0,0,1)^T$.

\begin{lemma}\label{tangentfinallemma}
		Let $\Gp$ be an $\ana$-partition, let the tensor fields $\mu^{-1},\varepsilon\in\anatenspw{\Gp}$ satisfy \eqref{coercivedomain} and suppose that $\zeta\in\anagamma$ satisfies \eqref{coercivebdr} and \eqref{zetaproperty}.
	In addition to that, let $0<\Rmin\leq R\leq\Rmax$ and let $\Upsilon:\PW_{R}\rightarrow\R^3$ be a local flattening of either the boundary $\Gamma$ or of an interface component $\interf_j$.

		Furthermore, suppose that $\solu\in\Hcurl$ is a piecewise smooth solution of Maxwell's equation \eqref{Maxwellorig} corresponding to a right-hand side $\solf\in\Hdiv$ and boundary data $\solgi$. 

		Under these assumptions we consider the transformed quantities $\wu$, $\hatf{\solf}$, $\solgicheck$, $\wmu^{-1}$, $\weps$ and $\wzeta$ from Definition~\ref{trafoquantitiesdef}.
		Then, for all $\ell\in\N_0$ there holds 
		\begin{align}\label{firstcase}
				\summorrey{\wu}{\ell+1,0,\PM_R}\leq C\sum_{d=0}^{\ell}A^{\ell-d}\left[\summorrey{\wu}{d,0,\PM_R}+|k|^{-2}\solF_{k,d,R}+|k|^{-3/2}\morreyhalf{\solgicheck}_{d,\PZ_R}\right]
		\end{align}
	if $\Upsilon$ is a local flattening of $\Gamma$, and 
	\begin{align}\label{secondcase}
			\summorrey{\wu}{\ell+1,0,\PS_R}\leq C\sum_{d=0}^{\ell}A^{\ell-d}\left[\summorrey{\wu}{d,0,\PS_R}+|k|^{-2}\morrey{\wf}_{d,0,\PS_R}+|k|^{-2}\morreyone{\diverg\wf}_{d,0,\PS_R}\right]
	\end{align}
	if $\Upsilon$ is a local flattening of an interface component $\interf_j$. In both cases, the constants $C>0$ and $A\geq 1$ depend only on $\Rmin$, $\Rmax$, $\wmu^{-1}$, $\weps$ and $\wzeta$. 
\end{lemma}

\begin{fatproof}
		We treat the two cases of $\Upsilon$ being a local flattening of $\Gamma$ and of $\Upsilon$ being a local flattening of an interface component $\interf_j$ separately. The remainder of the proof is divided into two steps, where the first step deals with $\Upsilon$ being a local flattening of $\Gamma$ and the second step treats the case of $\Upsilon$ being a local flattening of an interface component $\interf_j$.

\medskip 

{\bf Case 1:} Let $\Upsilon$ be a local flattening of $\Gamma$. Let $\ell\in\N_0$, let $\delta'\in\N_0^2\times\{0\}$ with $|\delta'|=\ell+1$ be arbitrary, and write $\delta'=\alpha'+\beta'$ for multiindices $\alpha',\beta'\in\N_0^2$ satisfying $|\alpha'|=1$ and $|\beta'| = \ell$. Moreover, let $\rho\in\left(0,\frac{R}{2(\ell+1+|k|)}\right)$ and define $\tau:=\ell+|k|$. Then, from the properties of the cutoff function $\chi_{R-\tau\rho,\rho}$ and Proposition~\ref{finiteshifttrafo} we get
\begin{align}\label{lemma47i1}
		\begin{split}
			\norm{\D^{\delta'}(\wu)}{\PM_{R-(\tau+1)\rho}}&\leq C\norm{\chi_{R-\tau\rho,\rho}\D^{\beta'}(\wu)}{\VSHM{1}\left(\PM_{R-\tau\rho}\right)} \\
		&\leq C\norm{\curl\wv}{\PM_{R-\tau\rho}}+C\norm{\diverg\weps\wv}{\PM_{R-\tau\rho}} 
		 +C\norm{\wv_T}{\thalf,\PZ_{R-\tau\rho}},
	\end{split}
	\end{align}
	where $\wv := \chi_{R-\tau\rho,\rho}\D^{\beta'}(\wu)$. Furthermore, from Proposition~\ref{shiftprop} and \eqref{surfacetraceprop} we infer
	\begin{align*}
			\norm{\wv_T}{\thalf,\PZ_{R-\tau\rho}}&\leq C\norm{\divergt\wzeta\wv_T}{\thalfm,\PZ_{R-\tau\rho}}+C\norm{\curlt\wv_T}{\thalfm,\PZ_{R-\tau\rho}}\\
		& \leq C\norm{\divergt\wzeta\wv_T}{\thalfm,\PZ_{R-\tau\rho}}+C\norm{\curl\wv}{\PM_{R-\tau\rho}}.
	\end{align*}
	Hence, \eqref{lemma47i1} improves to
\begin{align*}
		\rho^{\ell+1}\norm{\D^{\delta'}(\wu)}{\PM_{R-(\tau+1)\rho}}\leq C\rho^{\ell+1}\left(\norm{\curl\wv}{\PM_{R-\tau\rho}}+\norm{\diverg\weps\wv}{\PM_{R-\tau\rho}}+\norm{\divergt\wzeta\wv_T}{\thalfm,\PZ_{R-\tau\rho}}\right),
\end{align*}
from which the auxiliary estimates provided by Lemma~\ref{curldivbdry} and Lemma~\ref{devialemma} yield
\begin{align}\label{secondcasefirst}
		\morrey{\wu}_{\ell+1,0,\PM_R}\leq C\sum_{d=0}^{\ell}A^{\ell-d}\left[\summorrey{\wu}{d,0,\PM_R}+|k|^{-2}\solF_{k,d,R}+|k|^{-3/2}\morreyhalf{\solgicheck}_{d,\PZ_R}\right]
\end{align}
for appropriate constants $C>0$ and $A\geq 1$.

\medskip 	

It remains to find an estimate for $|k|^{-1}\morrey{\wmu^{-1}\curl\wu}_{\ell+1,0,\PM_R}$. \revision{Because $\chi_{R-\tau\rho}$ is} identically one inside of $\PM_{R-(\tau+1)\rho}$ we notice
\begin{align*}
		\norm{\D^{\delta'}(\wmu^{-1}\curl\wu)}{\PM_{R-(\tau+1)\rho}}\leq\norm{\chi_{R-\tau\rho,\rho}\D^{\beta'}(\wmu^{-1}\curl\wu)}{\VSHM{1}\left(\PM_{R-\tau\rho}\right)}.
	\end{align*}
	From this, Proposition~\ref{finiteshifttrafo} proves
	\begin{align}\label{lemma47i2}
	\begin{split}
	\norm{\D^{\delta'}(\wmu^{-1}\curl\wu)}{\PM_{R-(\tau+1)\rho}}&\leq C\left(\norm{\curl\ww}{\PM_{R-\tau\rho}}
	+\norm{\diverg\wmu\ww}{\PM_{R-\tau\rho}} +\norm{\ww_t}{\thalf,\PZ_{R-\tau\rho}}\right),
	\end{split}
	\end{align}
    where $\ww:=\chi_{R-\tau\rho,\rho}\D^{\beta'}(\wmu^{-1}\curl\wu)$.
	Due to Proposition~\ref{shiftprop} and \eqref{surfacetraceprop} we have
	\begin{align*}
			\norm{\ww_t}{\thalf,\PZ_{R-\tau\rho}}&\leq C\norm{\divergt\ww_t}{-\frac{1}{2},\PZ_{R-\tau\rho}}+C\norm{\curlt\wzeta^{-1}\ww_t}{-\thalf,\PZ_{R-\tau\rho}}\\
		&\leq C\norm{\curl\ww}{\PM_{R-\tau\rho}}+C\norm{\curlt\wzeta^{-1}\ww_t}{-\thalf,\PZ_{R-\tau\rho}},
	\end{align*}
	hence \eqref{lemma47i2} becomes
	\begin{align*}
     \norm{\D^{\delta'}(\wmu^{-1}\curl\wu)}{\PM_{R-(\tau+1)\rho}}&\leq C\left(\norm{\curl\ww}{\PM_{R-\tau\rho}}
	+\norm{\diverg\wmu\ww}{\PM_{R-\tau\rho}} +\norm{\curlt\wzeta^{-1}\ww_t}{-\thalf,\PZ_{R-\tau\rho}}
\right).
	\end{align*}
	From this, the auxiliary estimates provided by Lemma~\ref{curldivbdry} and Lemma~\ref{devialemma} prove
	\begin{align*}
			|k|^{-1}\morrey{\wmu^{-1}\curl\wu}_{\ell+1,0,\PM_R}\leq C\sum_{d=0}^{\ell}A^{\ell-d}\left[\summorrey{\wu}{d,0,\PM_R}+|k|^{-2}\solF_{k,d,R}+||k|^{-3/2}\morreyhalf{\solgicheck}_{d,\PZ_R}\right].
	\end{align*}
	Together with \eqref{secondcasefirst} this proves \eqref{firstcase} if $\Upsilon$ is a local flattening of the boundary $\Gamma$.

	\medskip

	\medskip 

	{\bf Case 2:} In the second case we suppose that $\Upsilon$ is a local flattening of an interface component $\interf_j$. In contrast to the first case we cannot employ Proposition~\ref{finiteshifttrafo} here, instead we have to rely on Proposition~\ref{helmholtzmod}.

	\medskip 

	As before, let $\ell\in\N_0$ and $\delta'\in\N_0^2\times\{0\}$ with $|\delta'|=\ell+1$ and $\delta'=\alpha'+\beta'$ for two multiindices $\alpha',\beta'\in\N_0^2\times\{0\}$ satisfying $|\alpha'|=1$ and $|\beta'|=\ell$. 
	Moreover, let $\rho\in\left(0, \frac{1}{2(\ell+1+|k|)}\right)$ and set $\tau:=\ell+|k|$. 
	Then, due to $\chi_{R-\tau\rho,\rho}=1$ in $\PW_{R-(\tau+1)\rho}$ and $\D^{\delta'}$ being a tangential derivative there holds 
	\begin{align}\label{interffirst}
	\begin{split}
			\rho^{\ell+1}\norm{\D^{\delta'}(\wu)}{\PW_{R-(\tau+1)\rho}}&\leq \rho^{\ell+1}\norm{\D^{\alpha'}\chi_{R-\tau\rho,\rho}\D^{\beta'}(\wu)}{\PW_{R-\tau\rho}} \\ &\leq\rho^{\ell+1}\norm{\chi_{R-\tau\rho,\rho}\D^{\beta'}(\wu)}{\VSHM{1}(\PP_{R-\tau\rho})}+\rho^{\ell+1}\norm{\chi_{R-\tau\rho,\rho}\D^{\beta'}(\wu)}{\VSHM{1}(\PM_{R-\tau\rho})}.
	\end{split}
	\end{align}
	
	Due to $\D^{\beta'}(\wu)$ being a tangential derivative of $\wu$ we have $\chi_{R-\tau\rho,\rho}\D^{\beta'}(\wu)\in\Hcurldomain{\PW_R}$ and we notice that $\solz:=\chi_{R-\tau\rho,\rho}\D^{\beta'}(\wu)$ has compact support in $\PW_{R-\tau\rho}$. \revision{As a consequence we may apply Proposition~\ref{helmholtzmod} to obtain
			\begin{align*}
					\norm{\solz}{\VSHM{1}(\PP_{R-\tau\rho})}+\norm{\solz}{\VSHM{1}(\PM_{R-\tau\rho})}\leq C\left(\norm{\curl\solz}{\PW_{R-\tau\rho}}+\norm{\diverg\weps\solz}{\PS_{R-\tau\rho}}+\norm{\njump{\weps\solz}}{\thalf,\PZ_{R-\tau\rho}}\right).
			\end{align*}
			Together with Lemma~\ref{curldivbdry} and Lemma~\ref{interfaceest} and using the abbreviation 
	$$\norm{\solz}{\VSHM{1}(\PS_{R-\tau\rho})}:= \norm{\solz}{\VSHM{1}(\PP_{R-\tau\rho})}+\norm{\solz}{\VSHM{1}(\PM_{R-\tau\rho})}$$ 
			we get
\begin{align*}
\norm{\solz}{\VSHM{1}(\PS_{R-\tau\rho})}\leq C\rho^{-\ell-1}\sum_{d=0}^{\ell}A^{\ell-d}\left[\summorrey{\wu}{d,0,\PS_R}+|k|^{-2}\morrey{\wf}_{d,0,\PS_R}+|k|^{-2}\morreyone{\diverg\wf}_{d,0,\PS_R}\right].
\end{align*}
Together with \eqref{interffirst} we arrive at
\begin{align}\label{infinal}
		\morrey{\wu}_{\ell+1,0,\PS_R}\leq C\sum_{d=0}^{\ell}A^{\ell-d}\left[\summorrey{\wu}{d,0,\PS_R}+|k|^{-2}\morrey{\wf}_{d,0,\PS_R}+|k|^{-2}\morreyone{\diverg\wf}_{d,0,\PS_R}\right]
\end{align}
for appropriate constants $C>0$ and $A\geq 1$. 
	}

\medskip 

The final part of the proof is to find an estimate on $\D^{\delta'}(\wmu^{-1}\curl\wu)$.
The first observation is that
	\begin{align}\label{curlstarttransm}
		|k|^{-1}\rho^{\ell+1}\norm{\D^{\delta'}(\wmu^{-1}\curl\wu)}{\PS_{R-(\tau+1)\rho}}&\leq C|k|^{-1}\rho^{\ell+1}\norm{\chi_{R-\tau\rho,\rho}\D^{\beta'}(\wmu^{-1}\curl\wu)}{\VSHM{1}(\PS_{R-\tau\rho})}.
	\end{align}
	
	\revision{
			Noting that $\solw:=\chi_{R-\tau\rho,\rho}\D^{\beta'}(\wmu^{-1}\curl\wu)\in\Hcurldomain{\PW_R}$ has compact support in $\PW_{R-\tau}$, we may employ Proposition~\ref{helmholtzmod} to obtain 
	\begin{align*}
			\norm{\solw}{\VSHM{1}(\PS_{R})}\leq C\left(\norm{\curl\solw}{\PW_{R-\tau\rho}}+\norm{\diverg\wmu\solw}{\PS_{R-\tau\rho}}+\norm{\njump{\wmu\solw}}{\thalf,\PZ_{R-\tau\rho}}\right)
	\end{align*}
	Hence, Lemma~\ref{curldivbdry}, Lemma~\ref{interfaceest} and \eqref{curlstarttransm} prove
\begin{align*}
			|k|^{-1}\morrey{\wmu^{-1}\curl\wu}_{\ell+1,0,\PS_R}\leq	C\sum_{d=0}^{\ell}A^{\ell-d}\left[\summorrey{\wu}{d,0,\PS_R}+|k|^{-2}\morrey{\wf}_{d,0,\PS_R}+|k|^{-2}\morreyone{\diverg\wf}_{d,0,\PS_R}\right]
	\end{align*}
	for some $C>0$ and $A\geq 1$. Together with \eqref{infinal} this concludes the proof of \eqref{secondcase}.
}
%	and $\diverg\wmu\solz=0$ in $\PW_R$. In combination with Lemma~\ref{curldivbdry} we obtain
%	\begin{align}\label{curlestinterface2}
%			|k|^{-1}\rho^{\ell+1}\norm{\solw}{\VSHM{1}(\PS_{R-\tau\rho})}\leq C\sum_{d=0}^{\ell}A^{\ell-d}\left[\summorrey{\wu}{d,0,\PS_R}+k^{-2}\morreyone{\wf}_{d,0,\PS_R}\right],
%	\end{align}
%	and from the properties of $\solw$ and $\nabla\psi$ and Proposition~\ref{interfacepoisson} we get 
%	\begin{align*}
%			\norm{\nabla\psi}{\VSHM{1}(\PS_{R-\tau\rho})}\leq C&\norm{\diverg \chi_{R-\tau\rho,\rho}\wmu\D^{\beta'}(\wmu^{-1}\curl\wu)}{\PS_{R-\tau\rho}}+C\norm{\njump{\chi_{R-\tau\rho,\rho}\wmu\D^{\beta'}(\wmu^{-1}\curl\wu)}}{\frac{1}{2},\PZ_{R-\tau\rho}}. 
%	\end{align*}

%	Hence, Lemma~\ref{curldivbdry}, Lemma~\ref{interfaceest} and \eqref{curlstarttransm}-\eqref{curlestinterface2} prove
%	\begin{align*}
%			|k|^{-1}\morrey{\wmu^{-1}\curl\wu}_{\ell+1,0,\PS_R}\leq	C\sum_{d=0}^{\ell}A^{\ell-d}\left[\summorrey{\wu}{d,0,\PS_R}+|k|^{-2}\morrey{\wf}_{d,0,\PS_R}+|k|^{-2}\morreyone{\diverg\wf}_{d,0,\PS_R}\right]
%	\end{align*}
%	for some $C>0$ and $A\geq 1$. Together with \eqref{infinal} this concludes the proof of \eqref{secondcase}.

\end{fatproof}

%% file: analytic/analytic_tangential_3.tex
Lemma~\ref{tangentfinallemma} provides recursive estimates for tangential derivatives of a transformed solution $\wu$. Resolving these recursive estimates leads to the subsequent corollary; its proof mainly consists of an induction argument and some computations and can be found in Appendix~\ref{appendix:commest}.

%In order to appropriately capture the influence of the boundary data we recall from Section~\ref{mainresultssec} that, depending on the imposed boundary conditions, 
%\begin{align*}
%		&\solB_{k}:= |k|^{-3/2}C_{\solgi}\quad && {\rm in\ the\ case\ of\ } \eqref{Maxwellimpedance}, \\
%		&\solB_{k}:= |k|^{-3/2}C_{\solgn}+|k|^{-5/2}C_{\diverg\solgn}+|k|^{-5/2}C_{\solf\cdot\soln} \quad && {\rm in\ the\ case\ of\ }\eqref{Maxwellnatural},
% \\
%		&\solB_{k}:= 0 \quad && {\rm in\ the\ case\ of\ }\eqref{Maxwellessential}, 
%\end{align*}
%where the meaning of the constants $C_{\solgi}, C_{\solgn}, C_{\diverg\solgn}$ and $C_{\solf\cdot\soln}$ is specified e.g. below in Corollary~\ref{analogontangential}.
%With these definitions, there holds the following result:
%

\begin{corollary}\label{analogontangential}
		In addition to the hypotheses of Lemma~\ref{tangentfinallemma} suppose that there holds $\solf\in\anaomegavecpw{\Gp}{\omega_{\solf}}{M}$ and $\diverg\solf\in\anaomegapw{\Gp}{\omega_{\diverg\solf}}{M}$, as well as $\solgi\in\anaomegavec{\Gamma}{\omega_{\solgi}}{M}$. 

Then, if $\Upsilon:\PW_R\rightarrow\R^3$ is a local flattening of the boundary $\Gamma$ there holds
\begin{align}\label{analogontangential1}
		\summorrey{\wu}{\ell,0,\PM_R}&\leq \left(\summorrey{\wu}{0,0,\PM_R}+|k|^{-2}\omega_{\solf}+|k|^{-3}\omega_{\diverg\solf}+|k|^{-3/2}\omega_{\solgi}\right)P^{\ell}, 
\end{align}
 and if $\Upsilon$ is a local flattening of an interface component $\interf_j$ there holds
\begin{align}\label{analogontangential2}
		\summorrey{\wu}{\ell,0,\PS_R}&\leq \left(\summorrey{\wu}{0,0,\PS_R}+|k|^{-2}\omega_{\solf}+|k|^{-3}\omega_{\diverg\solf}\right)P^{\ell}. 
\end{align}
In both cases, the constant $P\geq 0$ depends only on $\Upsilon$, $\mu^{-1}$, $\varepsilon$, $\zeta$, $\Rmin$, $\Rmax$ and $M$.
\end{corollary}

\subsection{Analytic estimates on transformed solutions of Maxwell's equations}

After having derived a link between normal and tangential derivatives of $\wu$ in Lemma~\ref{linklemmalem} and analytic estimates for tangential derivatives of $\wu$ in Corollary~\ref{analogontangential}, it remains to combine these two results to obtain analytic estimates for the growth of all derivatives of $\wu$. Indeed, there holds the following lemma:

\begin{lemma}\label{bdryanalyticprefinal}
		Let $\Gp=\geom$ be an $\ana$-partition in the sense of Definition~\ref{partitiondef}, suppose that the tensor fields $\mu^{-1},\varepsilon\in\anatenspw{\Gp}$ satisfy \eqref{coercivedomain} and that $\zeta\in\anagamma$ satisfies \eqref{coercivebdr} and \eqref{zetaproperty}. 

		In addition, let $0<\Rmin\leq R\leq\Rmax$ and let $\Upsilon:\PW_{R}\rightarrow\R^3$ be a local flattening of the boundary $\Gamma$ or of a subdomain interface $\I_{j}$. Furthermore, let $\solu\in\Hcurl$ be a piecewise smooth solution of Maxwell's equation \eqref{Maxwellorig} \revision{for} a right-hand side $\solf\in\Hdiv\cap\anaomegavecpw{\Gp}{\omega_{\solf}}{M}$ with $\diverg\solf\in\anaomegapw{\Gp}{\omega_{\diverg\solf}}{M}$ and impedance data $\solgi\in\anaomegavec{\Gamma}{\omega_{\solgi}}{M}$. 

Under these conditions, we consider the transformed quantities $\wu$, $\wf$ and $\wmu^{-1}$ from Definition~\ref{trafoquantitiesdef}. Then, there exits a constant $L\geq 1$ such that for $\ell,n\in\N_0$ there holds
		\begin{align*}
				\summorrey{\wu}{\ell,n,\PM_R}&\leq \left(\summorrey{\wu}{0,0,\PM_R}+|k|^{-2}\omega_{\solf}+|k|^{-3}\omega_{\diverg\solf}+|k|^{-3/2}\omega_{\solgi}\right)L^{\ell+n},
		\end{align*}
		if $\Upsilon$ is a local flattening of $\Gamma$, and 
		\begin{align*}
				\summorrey{\wu}{\ell,n,\PS_R}&\leq \left(\summorrey{\wu}{0,0,\PS_R}+|k|^{-2}\omega_{\solf}+|k|^{-3}\omega_{\diverg\solf}\right)L^{\ell+n},
		\end{align*} 
		if $\Upsilon$ is a local flattening of an interface component $\interf_j$.
		Furthermore, $L$  depends only on $M$, $\mu^{-1}$, $\varepsilon$, $\Rmin$, $\Rmax$, $\Upsilon$ and $\zeta$. 
	
\end{lemma}
\begin{fatproof}
		We consider only the case of $\Upsilon$ being a local flattening of the boundary $\Gamma$, the case of $\Upsilon$ being a local flattening of $\interf_j$ follows from similar arguments.
	
	According to Lemma~\ref{linklemmalem} there holds 
	\begin{equation}\label{linklemmaeq}
			\begin{multlined}
					\summorrey{\wu}{\ell,n+1,\PM_R}\leq C\sum_{d=0}^{\ell+1}\sum_{b=0}^{n}A^{\ell+n+1-d-b}\left[\summorrey{\wu}{d,b,\PM_R}+|k|^{-2}\morrey{\wf}_{d,b,\PM_R}\right]\\
			+C|k|^{-2}\sum_{d=0}^{\ell+1}\sum_{b=0}^{n}A^{\ell+n+1-d-b}\morreyone{\diverg\wf}_{d,b,\PM_R}
			\end{multlined}
	\end{equation}
	for all $\ell,n\in\N_0$. The remainder of the proof is split into two steps.
	
	\medskip 
	
	{\bf Step 1:} In the first step we show that there exists a constant $N$ depending only on $C,A,\Rmin,\Rmax$ and $M$ such that for all $\ell,n\in\N_0$ there holds 
	\begin{align}\label{step1}
			C\sum_{d=0}^{\ell+1}\sum_{b=0}^nA^{\ell+n+1-d-b}\left[\morrey{\wf}_{d,b,\PM_R}+\morreyone{\diverg\wf}_{d,b,\PM_R}\right]\leq \left[\morrey{\wf}_{0,0,\PM_R}+\morreyone{\diverg\wf}_{0,0,\PM_R}\right]N^{\ell+n+2},
	\end{align}
	where $C,A$ are the constants from Lemma~\ref{linklemmalem}. For the rest of the proof we use the abbreviation
	\begin{align*}
	\solF_{k,0}:=\left[\morrey{\wf}_{0,0,\PM_R}+\morreyone{\diverg\wf}_{0,0,\PM_R}\right].
	\end{align*}
	
	Due to $\solf\in\anaomegavecpw{\Gp}{\omega_{\solf}}{M}$ and $\diverg\solf\in\anaomegapw{\Gp}{\omega_{\diverg\solf}}{M}$, as well as Lemma~\ref{interplay1} and Lemma~\ref{interplay2} we have that there exists a constant $T\geq 0$ depending only on $M$, $\Upsilon$, $\Rmin$ and $\Rmax$ such that
	\begin{align*}
			\forall\ell,n\in\N_0:\ \morrey{\wf}_{\ell,n,\PM_R}+\morreyone{\diverg\wf}_{\ell,n,\PM_R}\leq \left[\omega_{\solf}+|k|^{-1}\omega_{\diverg\solf}\right]T^{\ell+n}.
	\end{align*}
	We choose constants $Q\geq T$ and $N\geq \revision{A}$ such that 
	\begin{align*}
		\frac{C}{N}\left(\sum_{d=0}^{\infty}\left(\frac{A}{Q}\right)^d\right)^2\leq 1,
	\end{align*}
	and observe 
	\begin{align*}
			C\sum_{d=0}^{\ell+1}\sum_{b=0}^nA^{\ell+n+1-d-b}\left[\morrey{\wf}_{d,b,\PM_R}+\morreyone{\diverg\wf}_{d,b,\PM_R}\right]  
								&\leq C\left[\omega_{\solf}+|k|^{-1}\omega_{\diverg\solf}\right]\sum_{d=0}^{\ell+1}\sum_{b=0}^nA^{\ell+n+1-d-b}\revision{T}^{d+b} \\
								&\leq C\left[\omega_{\solf}+|k|^{-1}\omega_{\diverg\solf}\right]  \sum_{d=0}^{\ell+1}\sum_{b=0}^nA^{\ell+n+1-d-b}Q^{d+b} \\
								&\leq C\left[\omega_{\solf}+|k|^{-1}\omega_{\diverg\solf}\right]\revision{A}^{\ell+n+1}\left(\sum_{d=0}^{\infty}\left(\frac{A}{Q}\right)^d\right)^2\\ & \leq \left[\omega_{\solf}+|k|^{-1}\omega_{\diverg\solf}\right]N^{\ell+n+2}.
	\end{align*}
	This proves \eqref{step1}.
	
	\medskip 
	
	{\bf Step 2:}
	In the second step we prove that there exist constants $B,G\geq 1$ such that for all $\ell,n\in\N_0$ there holds 
	\begin{align}\label{step2bdr}
			\summorrey{\wu}{\ell,n,\PM_R}\leq \left(\summorrey{\wu}{0,0,\PM_R}+|k|^{-2}\omega_{\solf}+|k|^{-3}\omega_{\diverg\solf}+|k|^{-3/2}\omega_{\solgi}\right)G^{\ell} B^n.
	\end{align}
	Let $C>0$ and $A,P\geq 1$ be the constants from Corollary~\ref{analogontangential} and Lemma~\ref{linklemmalem}, and choose the constants $B,G\geq 1$ in such a way that both $N+C+P\leq G\leq\frac{1}{2}\sqrt{B}$ and 
	\begin{align*}
		\frac{C}{\sqrt{B}}\left(\sum_{d=0}^{\infty}\left(\frac{A}{G}\right)^d\right)^2\leq 1.
	\end{align*}
	We prove \eqref{step2bdr} by induction with respect to $n$ and note that the case $n=0$ follows from Corollary~\ref{analogontangential}. For arbitrary $\ell\in\N_0$, the inequality \eqref{linklemmaeq} and Step 1 show 
	\begin{align*}
			\summorrey{\wu}{\ell,n+1,\PM_R}\leq  C\sum_{d=0}^{\ell+1}\sum_{b=0}^{n}A^{\ell+n+1-d-b}\summorrey{\wu}{d,b,\PM_R}+|k|^{-2}\left[\omega_{\solf}+|k|^{-1}\omega_{\diverg\solf}\right]N^{\ell+n+2}.
	\end{align*}
	Thus the induction hypothesis yields
	\begin{align*}
			\summorrey{\wu}{\ell,n+1,\PP_R}&\leq C\boldsymbol{W}\sum_{d=0}^{\ell+1}\sum_{b=0}^nA^{\ell+n+1-d-b}G^{d}B^b + |k|^{-2}\left[\omega_{\solf}+|k|^{-1}\omega_{\diverg\solf}\right]N^{\ell+n+2}, 
	\end{align*}
	where
	\begin{align*}
			\boldsymbol{W}:=\left(\summorrey{\wu}{0,0,\PM_R}+|k|^{-2}\omega_{\solf}+|k|^{-3}\omega_{\diverg\solf}+|k|^{-3/2}\omega_{\solgi}\right).
	\end{align*}
	Hence, 
	\begin{align*}
			\summorrey{\wu}{\ell,n+1,\PM_R}&\leq \frac{C}{2}\boldsymbol{W}G^{\ell}B^{n+\frac{1}{2}}\left(\sum_{d=0}^{\infty}\left(\frac{A}{G}\right)^d\right)^2+\frac{1}{2}|k|^{-2}\left[\omega_{\solf}+|k|^{-1}\omega_{\diverg\solf}\right]G^{\ell}B^{n+1}\\
										   &\leq \left(\summorrey{\wu}{0,0,\PM_R}+|k|^{-2}\omega_{\solf}+|k|^{-3}\omega_{\diverg\solf}+|k|^{-3/2}\omega_{\solgi}\right)G^{\ell}B^{n+1},
	\end{align*}
	which proves \eqref{step2bdr}. Due to $G^{\ell}B^{n+1}\leq B^{\ell+n+1}$, defining $L := B$ concludes the proof.
\end{fatproof}

%% file: appendix/appendix1.tex
\section{Postponed proofs}\label{appendix:commest}

It remains to provide proofs for some technical results used throughout this work. 
The following lemma turns out to be helpful on multiple occasions. 

\begin{lemma}
Let $\ell\in\N_0$ and $k\in\Co$. Then, there holds 
\begin{align}\label{elementaryestimate}
		\frac{\ell^{\ell}}{d^d}\leq \frac{(\ell+|k|)^{\ell}}{(d+|k|)^d}
\end{align}
for all $d=0,1,\ldots,\ell$.
\end{lemma}

\begin{fatproof}
		We note that \eqref{elementaryestimate} is equivalent to 
		\begin{align*}
				\left(1+\frac{|k|}{d}\right)^d\leq \left(1+\frac{|k|}{\ell}\right)^{\ell}.
		\end{align*}
		Since for fixed $k$ the sequence $n\mapsto \left(1+\frac{|k|}{n}\right)^n$ is increasing, the statement follows.
\end{fatproof}

\subsection{Commutator estimates}
Subsequently, we provide proofs for the commutator estimates that were exploited extensively in Section~\ref{analyticinteriorsec} and Section~\ref{auxiliarsec}. We begin with the proof of Lemma \ref{commutatorestimate}, which is similar to the proof of \cite[Lem.~1.6.2]{GrandLivre}. 

\medskip

\begin{fatproofmod}{Lemma \ref{commutatorestimate}}
		By a geometric series argument, the Sobolev embedding theorem and Stirling's formula we notice that $\nu\in\anatens{\PW_R}$ is equivalent to the existence of a constant $M>0$ such that
		\begin{align}\label{analyticnu}
				\sup_{x\in\PW_R}|\D^{\alpha}\nu_{i,j}(x)|\leq M^{|\alpha|+1}|\alpha|!
		\end{align}
		for all $i,j=1,2,3$ and all multiindices $\alpha\in\N_0^3$.

We recall that for any smooth vector field $\solv = (v_1,v_2,v_3)^T$ we write $\D^{\beta}(\solv):=(\D^{\beta}(v_1),\D^{\beta}(v_2),\D^{\beta}(v_3))^T$. Analogously, we define $\D^{\beta}(\nu):=(\D^{\beta}(\oldnu_{i,j}))_{i,j=1,2,3}$ for any smooth tensor field $\nu$. With this notation, the Leibniz formula shows for all multiindices $\alpha,\beta\in\N_0^3$
	\begin{align*}
		\D^{\alpha}(\nu\D^{\beta}(\solv)) = \sum_{\gamma\leq\alpha}\frac{\alpha!}{\gamma!(\alpha-\gamma)!}\D^{\gamma}(\nu)\D^{\alpha-\gamma}\D^{\beta}(\solv).
	\end{align*}
 \revision{We define
		\begin{align*}
				\commutetwo(\nu,\solv,\alpha,\beta):=\D^{\alpha}(\nu\D^{\beta}(\solv))-\nu\D^{\alpha}\D^{\beta}(\solv),
		\end{align*}
		then 
		\begin{align}\label{commutetworelation}
				\D^{\alpha}\commute(\nu,\solv,\beta)=\commutetwo(\nu,\solv,\alpha+\beta,0)-\commutetwo(\nu,\solv,\alpha,\beta).
		\end{align}
	From the combinatorial estimate \cite[p. 169]{ConcreteMaths}
	\begin{align*}
		\frac{\alpha!}{\gamma!(\alpha-\gamma)!}\leq\frac{|\alpha|!}{|\gamma|!(|\alpha|-|\gamma|)!}
	\end{align*}
and \eqref{analyticnu} we infer for all multiindices $\alpha,\beta\in\N_0^3$, $b:=|\alpha|+|\beta|$ and all $0<\rho\leq \frac{R}{2(b+|k|)}$ the estimate}
	\begin{align*}
		\rho^b\norm{\commutetwo(\nu,\solv,\alpha,\beta)}{R-(b-1+|k|)\rho} &= \rho^b\norm{\D^{\alpha}(\nu\D^{\beta}(\solv))-\nu\D^{\alpha}\D^{\beta}(\solv)}{R-(b-1+|k|)\rho} \\
		&\leq \rho^b\sum_{1\leq|\gamma|,\gamma\leq\alpha}M^{|\gamma|+1}\frac{|\alpha|!}{(|\alpha|-|\gamma|)!}\norm{\D^{\alpha-\gamma}\D^{\beta}(\solv)}{R-(b-1+|k|)\rho} \\
		&\leq \rho^b\sum_{g=1}^{|\alpha|}\sum_{|\gamma|=g,\gamma\leq\alpha}M^{g+1}\frac{|\alpha|!}{(|\alpha|-g)!}\max_{|\delta|=|\alpha-\gamma+\beta|}\norm{\D^{\delta}(\solv)}{R-(b-1+|k|)\rho}.
	\end{align*}
	Another combinatorial estimate shows $\sum_{|\gamma|=g}1\leq (g+1)^2$. Thus we can estimate 
	\begin{align}\label{tmpcommute}
	\begin{split}
	\rho^b\norm{\commutetwo(\nu,\solv,\alpha,\beta)}{R-(b-1+|k|)\rho}&\leq \rho^b\sum_{g=1}^{b-|\beta|}(g+1)^2M^{g+1}\frac{(b-|\beta|)!}{(b-|\beta|-g)!}\max_{|\delta|=b-g}\norm{\D^{\delta}(\solv)}{R-(b-1+|k|)\rho} \\
	&\leq \rho^b\sum_{d=|\beta|}^{b-1}(b-d+1)^2M^{b-d+1}\frac{(b-|\beta|)!}{(d-|\beta|)!}\max_{|\delta|=d}\norm{\D^{\delta}(\solv)}{R-(b-1+|k|)\rho}.
	\end{split}
	\end{align}
	We consider the two cases $|\beta|>0$ and $|\beta|=0$ separately and start by assuming $|\beta|>0$. In this case, we have 
	\begin{align*}
		\rho^b\norm{\commutetwo(\nu,\solv,\alpha,\beta)}{R-(b-1+|k|)\rho}&\leq \rho^b\sum_{d=1}^{b-1}(b-d+1)^2M^{b-d+1}\frac{b!}{d!}\max_{|\delta|=d}\norm{\D^{\delta}(\solv)}{R-(b-1+|k|)\rho} \\
		&\leq \sum_{d=1}^{b-1}(b-d+1)^2M^{b-d+1}\rho^{b-d}\frac{b!}{d!}\left(\rho^d\max_{|\delta|=d}\norm{\D^{\delta}(\solv)}{R-(b-1+|k|)\rho}\right).
	\end{align*}
	We define $N := \rho^b\norm{\commutetwo(\nu,\solv,\alpha,\beta)}{R-(b-1+|k|)\rho}$, and by applying Stirling's formula to the factorials $b!$ and $d!$ we see that there exists a universal constant $B>0$ such that 
	\begin{align*}
		N&\leq B \sum_{d=1}^{b-1}(b-d+1)^2M^{b-d+1}\rho^{b-d}e^{d-b}\sqrt{\frac{b}{d}}\frac{b^b}{d^d}\left(\rho^d\max_{|\delta|=d}\norm{\D^{\delta}(\solv)}{R-(b-1+|k|)\rho}\right) \\
		&\leq B \sum_{d=1}^{b-1}(b-d+1)^2M^{b-d+1}\rho^{b-d}e^{d-b}\sqrt{\frac{b}{d}}\frac{(b+|k|)^b}{(d+|k|)^d}\left(\rho^d\max_{|\delta|=d}\norm{\D^{\delta}(\solv)}{R-(b-1+|k|)\rho}\right),
	\end{align*}
	where the last estimate is due to \eqref{elementaryestimate}. 
	Moreover, by exploiting $\rho\leq\frac{R}{2(b+|k|)}$ we get 
	\begin{align*}
		N\leq B\sum_{d=1}^{b-1}(b-d+1)^2 M \left(\frac{RM}{2e}\right)^{b-d}\sqrt{\frac{b}{d}}\frac{(b+|k|)^d}{(d+|k|)^d}\left(\rho^d\max_{|\delta|=d}\norm{\D^{\delta}(\solv)}{R-(b-1+|k|)\rho}\right).
	\end{align*}
	 For every value of $d$ we define the quantity $\rho_d:=\frac{(b-1+|k|)\rho}{d+|k|}$; then we have $R-(b-1+|k|)\rho = R-(d+|k|)\rho_d$ and the above inequality becomes
	\begin{align*}
		N\leq B \sum_{d=1}^{b-1}(b-d+1)^2M\left(\frac{RM}{2e}\right)^{b-d}\sqrt{\frac{b}{d}}\frac{(b+|k|)^d}{(b-1+|k|)^d}\left(\rho_d^d\max_{|\delta|=d}\norm{\D^{\delta}(\solv)}{R-(d+|k|)\rho_d}\right).
	\end{align*}
	We note that 
	\begin{align*}
		\left(\frac{b+|k|}{b-1+|k|}\right)^d\leq e, \quad {\rm as\ well\ as}\quad \sqrt{\frac{b}{d}}\leq\frac{b}{d}\leq b-d+1\quad{\rm for\ all}\quad 1\leq d\leq b-1.
	\end{align*}
	and with these inequalities, we get 
	\begin{align*}
		N\leq Be \sum_{d=1}^{b-1}(b-d+1)^3M\left(\frac{RM}{2e}\right)^{b-d}\left(\rho_d^d\max_{|\delta|=d}\norm{\D^{\delta}(\solv)}{R-(d+|k|)\rho_d}\right).
	\end{align*}
	This shows 
	\begin{align*}
		N\leq  C\sum_{d=0}^{b-1}K^{b-1-d}\left(\rho_d^d\max_{|\delta|=d}\norm{\D^{\delta}(\solv)}{R-(d+|k|)\rho_d}\right)
	\end{align*}
	for appropriate constants $C,K\geq 1$ depending on $\Rmax$, $M$ and the universal constant $B$. Taking the maximum over all admissible $\rho_d$ then finishes the proof in the case $|\beta|> 0$. If $|\beta|=0$, then \eqref{tmpcommute} shows $\rho^b\norm{\D^{\alpha}\commutetwo(\nu,\solv,\alpha,\beta)}{R-(b-1+|k|)\rho}\leq N_0+N$, where $N$ is as in the case $|\beta|>0$ and 
	\begin{align*}
		N_0:=\rho^b(b+1)^2M^{b+1}b!\ \norm{\solv}{R-(b-1+|k|)\rho}.
	\end{align*}
	The quantity $N$ has been estimated above, and a similar estimate shows 
	\begin{align*}
		N_0\leq (b+1)^3M\left(\frac{RM}{2e}\right)^b\norm{\solv}{R-|k|\rho}.
	\end{align*}
	In total, we get \revision{
			\begin{align*}
		\rho^b\norm{\commutetwo(\nu,\solv,\alpha,\beta)}{R-(b-1+|k|)\rho}\leq C\sum_{d=0}^{b-1}K^{b-1-d}\left(\rho_d^d\max_{|\delta|=d}\norm{\D^{\delta}(\solv)}{R-(d+|k|)\rho_d}\right)
			\end{align*}
			for all multiindices $\alpha,\beta\in\N_0^3$ and $b:=|\alpha|+|\beta|$.
			This implies
\begin{align*}
		\rho^b\norm{\commutetwo(\nu,\solv,\alpha+\beta,0)}{R-(b-1+|k|)\rho}\leq C\sum_{d=0}^{b-1}K^{b-1-d}\left(\rho_d^d\max_{|\delta|=d}\norm{\D^{\delta}(\solv)}{R-(d+|k|)\rho_d}\right),
			\end{align*}
			hence \eqref{commutetworelation} and the triangle inequality prove
	}
	\begin{align}\label{commutatorstep}
		\rho^b\norm{\D^{\alpha}\commute(\nu,\solv,\beta)}{R-(b-1+|k|)\rho}\leq C\sum_{d=0}^{b-1}K^{b-1-d}\left(\rho_d^d\max_{|\delta|=d}\norm{\D^{\delta}(\solv)}{R-(d+|k|)\rho_d}\right)
	\end{align}
	for appropriate constants $C,K\geq 1$ depending on $\Rmax$, $M$ and the universal constant $B$. 
\end{fatproofmod}

The corresponding commutator estimate for anisotropic seminorms is given by Lemma \ref{commutatoranis}.

\bigskip 

\begin{fatproofmod}{Lemma \ref{commutatoranis}}
	Let $\alpha,\beta\in\N_0^3$ be multiindices and let $m:=|\alpha|+|\beta|$. By repeating the proof of Lemma~\ref{commutatorestimate} with some minor modifications we obtain an analog to \eqref{commutatorstep}, which reads as
	\begin{align*}
			\rho^{m}\norm{\D^{\alpha}(\commute(\nu,\solv,\beta))}{\PM_{R-(m-1+|k|)\rho}} \leq C\sum_{j=0}^{m-1}K^{m-1-j}\rho_j^j\max_{\substack{|\delta|=j \\ \delta\leq\alpha+\beta}}\norm{\D^{\delta}(\solv)}{\PM_{R-(j+|k|)\rho_j}},
	\end{align*}
	where $\rho_j := \frac{(m-1+|k|)}{j+|k|}\rho$.
	Splitting the derivatives into tangential derivatives and normal derivatives yields 
	\begin{align*}
			\rho^{m}\norm{\D^{\alpha}(\commute(\nu,\solv,\beta))}{\PM_{R-(m-1+|k|)\rho}} \leq &C\sum_{d=0}^{m'}\sum_{b=0}^{m^*-1}K^{m-1-d-b}\rho_{d+b}^{d+b}\max_{\substack{|\delta'|=d \\ |\gamma^*| = b}}\norm{\D^{\delta'+\gamma^*}(\solv)}{\PM_{R-(d+b+|k|)\rho_{d+b}}} \\
			&\hskip 1cm +   C\sum_{d=0}^{m'-1}K^{m'-1-d}\rho^{d+m^*}_{d+m^*}\max_{\substack{|\delta'|=d \\ |\gamma^*| = m^*}}\norm{\D^{\delta'+\gamma^*}(\solv)}{\PM_{R-(d+m^*+|k|)\rho_{d+b}}},
	\end{align*}
	and taking the maximum over all admissible $\rho_{d+b}$ and $\rho_{d,m^*}$ then proves 
	\begin{align*}
		\rho^{m}\norm{\D^{\alpha}(\commute(\nu,\solv,\beta))}{\PM_{R-(m-1+|k|)\rho}} \leq C\sum_{d=0}^{m'}\sum_{b=0}^{m^*-1}K^{m-1-d-b}\morrey{\solv}_{d,b,\PM_R}+C\sum_{d=0}^{m'-1}K^{m'-1-d}\morrey{\solv}_{d,m^*,\PM_R},
	\end{align*}
	which is \eqref{firstineq}.
	
	\medskip 
	
	Moreover, we note that for all multiindices $\alpha'$ and $\beta^*$ satisfying $|\alpha'|=\ell$ and $|\beta^*|=n$ and all numbers $\rho\in\left(0,\frac{R}{2(\ell+n+|k|)}\right)$ we have 
	\begin{align*}
		\rho^{\ell+n}\norm{\D^{\alpha'+\beta^*}(\nu\solv)}{\PM_{R-(\ell+n+|k|)\rho}}\leq \rho^{\ell+n}\norm{\commute(\nu,\solv,\alpha'+\beta^*)}{\PM_{R-(\ell+n-1+|k|)\rho}}+\rho^{\ell+n}\norm{\nu\D^{\alpha'+\beta^*}(\solv)}{\PM_{R-(\ell+n+|k|)\rho}}.
	\end{align*}
	We observe that 
	\begin{align*}
		\rho^{\ell+n}\norm{\nu\D^{\alpha'+\beta^*}(\solv)}{\PM_{R-(\ell+n+|k|)\rho}}\leq C'\rho^{\ell+n}\norm{\D^{\alpha'+\beta^*}(\solv)}{\PM_{R-(\ell+n+|k|)\rho}}\leq C'\morrey{\solv}_{\ell,n,\PM_R},
	\end{align*}
    for a constant $C'>0$, and according to \eqref{firstineq} there holds 
		\begin{align*}
	\rho^{\ell+n}\norm{\commute(\nu,\solv,\alpha'+\beta^*)}{\PM_{R-(\ell+n-1+|k|)\rho}}\leq C\sum_{d=0}^{\ell}\sum_{b=0}^{n-1}K^{\ell+n-1-d-b}\morrey{\solv}_{d,b,\PM_R}+C\sum_{d=0}^{\ell-1}K^{\ell-1-d}\morrey{\solv}_{d,n,\PM_R}.
	\end{align*}
	Combining the last two inequalities yields
	\begin{align*}
		\rho^{\ell+n}\norm{\D^{\alpha'+\beta^*}(\nu\solv)}{\PM_{R-(\ell+n+|k|)\rho}}\leq C\sum_{d=0}^{\ell}\sum_{b=0}^{n}K^{\ell+n-d-b}\morrey{\solv}_{d,b,\PM_R}
	\end{align*}
	for redefined constants $C$ and $K$, and thus
	\begin{align*}
		\morrey{\nu\solv}_{\ell,n,\PM_R} \leq C\sum_{d=0}^{\ell}\sum_{b=0}^{n}K^{\ell+n-d-b}\morrey{\solv}_{d,b,\PM_R}.
	\end{align*}
	This proves \eqref{secondineq}. Note that the proof works the same if $\PM_R$ is replaced by $\PP_R$, i.e., both inequalities stay valid for $\PP_R$.

\end{fatproofmod}

\subsection{Resolving recursive inequalities}

Subsequently, we present how to resolve the recursive inequalities derived in Section~\ref{analyticinteriorsec} and Section~\ref{auxiliarsec}. 

\medskip 

\begin{fatproofmod}{Lemma \ref{analyticinterior}}
	The proof of \eqref{estimateomega} is done in a similar way as \cite[Proof of Prop.~1.6.3]{GrandLivre}: The inequality \eqref{summorreyiteration} implies
	\begin{align}\label{newiterateomega}
	\summorrey{\solu}{\ell+1,R}\leq B\left(|k|^{-2}\morrey{\solf}_{\ell,R}+|k|^{-2}\morreyone{\diverg\solf}_{\ell,R}+\sum_{d=0}^{\ell}B^{\ell-d}\summorrey{\solu}{d,R}\right)
	\end{align}
	for a sufficiently large constant $B\geq 1$.
     We choose a constant $A\geq 1$ such that 
	\begin{align}\label{Aestimate}
	\sum_{d=0}^{\infty}\left(\frac{B}{A}\right)^{d+1}\leq 1,
	\end{align}
	and for this $A$ we will prove \eqref{estimateomega} by induction with respect to $\ell$. The case $\ell = 0$ follows from \eqref{newiterateomega}, since \eqref{Aestimate} implies $B\leq A$. 
	To ease notation, we abbreviate  $|k|^{-2}\morrey{\solf}_{\ell,R}+|k|^{-2}\morreyone{\diverg\solf}_{\ell,R}$ by ${\bf F}_{\ell}$ for all $\ell\in\N_0$.
	
	 Suppose that \eqref{estimateomega} is valid for $0,...,\ell-1$. Together with \eqref{newiterateomega} this yields 
	\begin{align}\label{omegaestimateprefinal}
	\summorrey{\solu}{\ell+1,R}\leq B{\bf F}_{\ell}+B\sum_{d=1}^{\ell}B^{\ell-d}\sum_{j=0}^{d-1}A^{d-j}{\bf F}_{j}+B\sum_{d=0}^{\ell}B^{\ell-d}A^d\summorrey{\solu}{0,R}.
	\end{align}
	We note that due to \eqref{Aestimate}, the coefficient $\oldlambda_0$ of $\summorrey{\solu}{0,R}$ satisfies
	\begin{align*}
	\oldlambda_0 = B\sum_{d=0}^{\ell}B^{\ell-d}A^d = BA^{\ell}\sum_{d=0}^{\ell}\left(\frac{B}{A}\right)^{\ell-d}\leq A^{\ell+1},
	\end{align*}
	and, again due to \eqref{Aestimate}, for fixed $j\in\{0,...,\ell-1\}$ the coefficient $\alpha_j$ of ${\bf F}_{j}$ satisfies 
	\begin{align*}
	\alpha_j = B\sum_{d=j+1}^{\ell}B^{\ell-d}A^{d-j} = BA^{\ell-j}\sum_{d=j+1}^{\ell}\left(\frac{B}{A}\right)^{\ell-d}\leq A^{\ell+1-j}.
	\end{align*}
	Together with \eqref{omegaestimateprefinal} these observations prove \eqref{estimateomega}.
	
	In order to prove \eqref{analyticityfinal}, we observe that as a consequence of the definitions there holds 
	\begin{align*}
			\forall\ell\in\N_0:\ \left(\frac{R}{2(\ell+1+|k|)}\right)^{\ell+1}\ \max_{|\delta|=\ell+1}\norm{\D^{\delta}(\solu)}{R/2}&\leq\morrey{\solu}_{\ell+1, R}, \\
			\forall d\in\N_0:\ \morrey{\solf}_{d,R}&\leq \left(\frac{R}{2(d+|k|)}\right)^{d}\ \max_{|\beta|=d}\norm{\D^{\beta}(\solf)}{R}, \\
			\forall d\in\N_0:\ \morreyone{\diverg\solf}_{d,R}&\leq \left(\frac{R}{2(d+|k|)}\right)^{d+1}\ \max_{|\beta|=d}\norm{\D^{\beta}(\diverg\solf)}{R}.
	\end{align*}
	Together with \eqref{estimateomega} these inequalities show for all $\ell\in\N_0$
	\begin{align}\label{analfinalpre}
	\begin{split}
			\left(\frac{R}{2(\ell+1+|k|)}\right)^{\ell+1}\ \max_{|\delta|=\ell+1}\norm{\D^{\delta}(\solu)}{R/2}\leq |k|&^{-2}\sum_{d=0}^{\ell}A^{\ell+1-d}\left(\frac{R}{2(d+|k|)}\right)^{d}\ \max_{|\beta|=d}\norm{\D^{\beta}(\solf)}{R} \\ &+ R|k|^{-3}\sum_{d=0}^{\ell}A^{\ell+1-d}\left(\frac{R}{2(d+|k|)}\right)^{d}\ \max_{|\beta|=d}\norm{\D^{\beta}(\diverg\solf)}{R} \\ 
	& \hskip 1.5cm+A^{\ell+1}\summorrey{\solu}{0,R}.
	\end{split}
	\end{align}
In combination with \eqref{analfinalpre}, some manipulations and straightforward estimates yield \eqref{analyticityfinal} for an appropriate constant $P\geq 1$. 
\end{fatproofmod}

Similar techniques as in the proof of Lemma~\ref{analyticinterior} can be employed to show Corollary~\ref{analogontangential}, as we shall see below.

\medskip 

\begin{fatproofmod}{Corollary~\ref{analogontangential}}
		We only prove \eqref{analogontangential1}, the estimate \eqref{analogontangential2} follows from analogous arguments. 

	As in Lemma~\ref{tangentfinallemma} we employ the abbreviation  
\begin{align*}
		\solF_{k,\ell,R}:= \morrey{\wf}_{\ell,0,\PM_R}+\morreyone{\diverg\wf}_{\ell,0,\PM_R}.
\end{align*}
%\begin{align*}
%		&\solG_{k,\ell,R}:= |k|\morreyhalf{\solgncheck}_{\ell,\PZ_R}+\morreyhalf{\diverg_{\PZ_R}\solgncheck}_{\ell,\PZ_R}+\morreyhalf{\wf\cdot\sole_3}_{\ell,\PZ_R} \hskip 1.59cm\ {\rm in\ the\ case\ of\ }\eqref{Maxwellnatural},
% \\
%		&\solG_{k,\ell,R}:= 0 \hskip 10.17cm\ {\rm in\ the\ case\ of\ }\eqref{Maxwellessential}. 
%\end{align*}
%

	The remainder of the proof is divided into two steps. In the first step we show that there exists a constant  $Q\geq 1$ such that for all $\ell\in\N_0$ there holds
	\begin{align}\label{tangentcontrol}
			C\sum_{d=0}^{\ell}A^{\ell-d}\left[\solF_{k,d,R}+|k|\morreyhalf{\solgicheck}_{d,\PZ_R}\right]\leq \frac{1}{2}\left[\omega_{\solf}+|k|^{-1}\omega_{\diverg\solf}+|k|^{1/2}\omega_{\solgi}\right]Q^{\ell+1},
	\end{align}
	where $C$ and $A$ are the constants from Lemma~\ref{tangentfinallemma}. Afterwards, in the second step we conclude.
	
	{\bf Step 1:}  In combination with some basic computations, Lemma~\ref{interplay1} and Lemma~\ref{interplay2} show that there exists a constant $M'\geq 1$ such that for all $d\in\N_0$ there holds 
\begin{align*}
		\solF_{k,d,R}+|k|\morreyhalf{\solgicheck}_{d,\PZ_R}\leq \left[\omega_{\solf}+|k|^{-1}\omega_{\diverg\solf}+|k|^{1/2}\omega_{\solgi}\right]M'^{d}.
\end{align*}

 Let $Q\geq M'$ be such that 
	\begin{align*}
		2C\sum_{d=0}^{\infty}\left(\frac{A}{Q}\right)^d\leq Q.
	\end{align*}
	With this choice we get
	\begin{align*}
			C\sum_{d=0}^{\ell}A^{\ell-d}\left[\solF_{k,d,R}+|k|\morreyhalf{\solgicheck}_{d,\PZ_R}\right]&\leq C\left[\omega_{\solf}+|k|^{-1}\omega_{\diverg\solf}+|k|^{1/2}\omega_{\solgi}\right]Q^{\ell}\sum_{d=0}^{\ell}\left(\frac{A}{Q}\right)^{\ell-d}\\
		&\leq C\left[\omega_{\solf}+|k|^{-1}\omega_{\diverg\solf}+|k|^{1/2}\omega_{\solgi}\right]Q^{\ell}\sum_{d=0}^{\infty}\left(\frac{A}{Q}\right)^d\\
		&\leq \frac{1}{2}\left[\omega_{\solf}+|k|^{-1}\omega_{\diverg\solf}+|k|^{1/2}\omega_{\solgi}\right]Q^{\ell+1},
	\end{align*}
	which proves \eqref{tangentcontrol}.
	
	{\bf Step 2:}  Let $P\geq Q$ be such that there holds 
	\begin{align*}
		2C^2\sum_{d=0}^{\infty}\left(\frac{A}{P}\right)^d\leq P.
	\end{align*}
	We proceed by induction with respect to $\ell$. For $\ell = 0$ we note that \eqref{analogontangential1} is trivial, and for arbitrary $\ell\in\N_0$ we have that Lemma~\ref{tangentfinallemma}, the inequality \eqref{tangentcontrol} and the induction hypothesis show 
    \begin{align*}
			\summorrey{\wu}{\ell+1,0,\PM_R}&\leq C\sum_{d=0}^{\ell}A^{\ell-d}\left(\summorrey{\wu}{d,0,\PM_R}+|k|^{-2}\solF_{k,d,R}+|k|^{-1}\morreyhalf{\solgicheck}_{d,\PM_R}\right) \\
										   &\leq C^2\left(\summorrey{\wu}{0,0,\PM_R}+|k|^{-2}\omega_{\solf}+|k|^{-3}\omega_{\diverg\solf}+|k|^{-3/2}\omega_{\solgi}\right)\sum_{d=0}^{\ell}A^{\ell-d}P^d\\
										   &\hskip 3cm +\frac{1}{2}\left[|k|^{-2}\omega_{\solf}+|k|^{-3}\omega_{\diverg\solf}+|k|^{-3/2}\omega_{\solgi}\right]Q^{\ell+1}\\ 
										   &\leq \left(\summorrey{\wu}{0,0,\PM_R}+|k|^{-2}\omega_{\solf}+|k|^{-3}\omega_{\diverg\solf}+|k|^{-3/2}\omega_{\solgi}\right)P^{\ell+1}.
    \end{align*} 
    This finishes the proof.

\end{fatproofmod}

%% file: appendix/appendix2.tex
\subsection{Wavenumber-explicit analytic change of variables}

It remains to give a proof for the analytic change of variables estimate from Lemma~\ref{analyticchangevars}. 

\medskip

\begin{fatproofmod}{Lemma~\ref{analyticchangevars}}
		We only consider the case that $\Upsilon$ is a local flattening of $\Gamma$; the case of $\Upsilon$ being a local flattening of a subdomain component $\I_{j}$ follows the same lines. 
		
		\medskip

		From \cite[(1.22)]{GrandLivre} and Stirling's formula we infer
		\begin{align*}
				\norm{\D^{\gamma}(\solv)}{\VSL\left(\mathcal{D}_R\cap\Omega\right)}\leq C^{\ell+n+1}\sum_{d=0}^{\ell}\sum_{b=0}^n\frac{(\ell+n)^{\ell+n}}{(d+b)^{d+b}}\max_{\substack{|\delta'|=d \\ |\alpha^*|=b}}\norm{\D^{\delta'+\alpha^*}(\solv\circ\Upsilon)}{\PM_{R}}.
		\end{align*}
		Let $R(\ell,n,k):=\frac{R}{\ell+n+|k|}$ and choose an arbitrary $0<\rho\leq 2R(\ell,n,k)$. Due to \eqref{elementaryestimate} we get	
\begin{align*}
		\rho^{\ell+n}\norm{\D^{\gamma}(\solv)}{\VSL\left(\mathcal{D}_R\cap\Omega\right)}\leq C^{\ell+n+1}\sum_{d=0}^{\ell}\sum_{b=0}^n\frac{(\ell+n+|k|)^{\ell+n}}{(d+b+|k|)^{d+b}}\rho^{\ell+n}\max_{\substack{|\delta'|=d \\ |\alpha^*|=b}}\norm{\D^{\delta'+\alpha^*}(\solv\circ\Upsilon)}{\PM_{2R-\tau\rho}},
		\end{align*}
		where $\tau:=\ell+n+|k|$.
		By defining $\rho_{d+b}':=\frac{\ell+n+|k|}{d+b+|k|}\rho$ and exploiting $\rho\leq 2R(\ell,n,k)$ we obtain
\begin{align*}
		\rho^{\ell+n}\norm{\D^{\gamma}(\solv)}{\VSL\left(\mathcal{D}_R\cap\Omega\right)}\leq C^{\ell+n+1}\sum_{d=0}^{\ell}\sum_{b=0}^nB^{\ell+n-d-b}\rho_{d+b}'^{d+b}\max_{\substack{|\delta'|=b \\ |\alpha^*|=d}}\norm{\D^{\delta'+\alpha^*}(\solv\circ\Upsilon)}{\PM_{2R-(d+b+|k|)\rho_{d+b}}}
		\end{align*}
		for some constants $B,C\geq 1$.
Taking the maximum over $\rho$ on both sides and rearranging expressions yield 
\begin{align*}
		\frac{1}{(\ell+n+|k|)^{\ell+n}}\norm{\D^{\gamma}(\solv)}{\VSL\left(\mathcal{D}_R\cap\Omega\right)}\leq C^{\ell+n+1}\sum_{d=0}^{\ell}\sum_{b=0}^nB^{\ell+n-d-b}\morrey{\solv\circ\Upsilon}_{d,b,\PM_{2R}}
\end{align*}
for appropriate constants $B,C\geq 1$. Furthermore, in combination with \eqref{secondineq} and Lemma~\ref{vierersumme} a calculation leads to
\begin{align*}
		\sum_{d=0}^{\ell}\sum_{b=0}^nB^{\ell+n-d-b}\morrey{\solv\circ\Upsilon}_{d,b,\PM_{2R}}&= 
		\sum_{d=0}^{\ell}\sum_{b=0}^nB^{\ell+n-d-b}\morrey{F^{-T}F^T(\solv\circ\Upsilon)}_{d,b,\PM_{2R}}\\
																							 &\leq C\sum_{d=0}^{\ell}\sum_{b=0}^nA^{\ell+n-b-d}\morrey{\wv}_{b,d,\PM_{2R}}, 
\end{align*}
thus 
\begin{align*}
		\frac{1}{(\ell+n+|k|)^{\ell+n}}\norm{\D^{\gamma}(\solv)}{\VSL\left(\mathcal{D}_R\cap\Omega\right)}\leq C^{\ell+n+1}\sum_{d=0}^{\ell}\sum_{b=0}^nA^{\ell+n-d-b}\morrey{\wv}_{d,b,\PM_{2R}}
\end{align*}
for appropriate constants $C,A\geq 1$. This concludes the proof.
\end{fatproofmod}